\documentclass[11pt]{article}
\usepackage[english]{babel}
\usepackage[utf8]{inputenc}
\usepackage[OT1]{fontenc}
\usepackage[a4paper, left=2.5cm, right=2.5cm, top=2.5cm, bottom=2.5cm]{geometry}
\usepackage{amsfonts, amsmath, amsthm, amssymb, enumitem, bbm, hyperref, longtable, natbib, graphicx, footnote, xcolor, float,mathtools,authblk}
\usepackage{setspace}
\usepackage[title]{appendix}
\usepackage{endnotes}
\newtheorem{theorem}{Theorem}
\newtheorem{corollary}{Corollary}
\newtheorem{proposition}{Proposition}
\newtheorem{lemma}{Lemma}

\theoremstyle{definition}

\newtheorem{example}{Example}
\theoremstyle{remark}

\DeclareMathOperator*{\esssup}{ess\,sup}
\DeclareMathOperator*{\essinf}{ess\,inf}
\title{Constructing Uncertainty Sets for Robust Risk Measures: A Composition of $\phi$-Divergences Approach to Combat Tail Uncertainty}
\author[a]{Guanyu Jin}
\author[a,c,d]{Roger J.~A.~Laeven}
\author[b]{Dick den Hertog}
\author[e]{Aharon Ben-Tal}
\affil[a]{{\small Dept.~of Quantitative Economics, University of Amsterdam,\protect\\ 
1001 NJ Amsterdam, The Netherlands}}
\affil[b]{{\small Dept.~of Business Analytics, University of Amsterdam,\protect\\  1001 NJ Amsterdam, The Netherlands}}
\affil[c]{{\small EURANDOM, 5600 MB Eindhoven, The Netherlands}}
\affil[d]{{\small CentER, Tilburg University, 5000 LE Tilburg, The Netherlands}}
\affil[e]{{\small Faculty of Industrial Engineering and Management, Technion – Israel Institute of Technology}}
\begin{document}
\maketitle
\onehalfspacing
\begin{abstract}
Risk measures, which typically evaluate the impact of extreme losses, are highly sensitive to misspecification in the tails. This paper studies a robust optimization approach to combat tail uncertainty by proposing a unifying framework to construct uncertainty sets for a broad class of risk measures, given a specified nominal model. Our framework is based on a parametrization of robust risk measures using two (or multiple) $\phi$-divergence functions, which enables us to provide uncertainty sets that are tailored to both the sensitivity of each risk measure to tail losses and the tail behavior of the nominal distribution. In addition, our formulation allows for a tractable computation of robust risk measures, and elicitation of $\phi$-divergences that describe a decision maker's risk and ambiguity preferences.
\end{abstract}

\renewcommand{\endnote}[1]{\footnote{#1}}

\section{Introduction}\label{sec:Intro}
Model misspecification, especially in the tails, is an inevitable issue that underlies the practice of financial and operational risk management. One way to safeguard against model uncertainty is distributionally robust optimization (DRO). The idea is to take a set of plausible models, often referred to as the uncertainty/ambiguity set, and evaluate the worst-case risk among all these models. DRO has its roots in the classical robust optimization pioneered by \citet{RObook}, as well as in economics and decision theory, by the works of \citet{MaxminEU}, \citet{Maccheroni}, and \citet{Hansen_Sargent}.

An essential part of distributional robust optimization is to choose a proper uncertainty set. In data-driven optimization, DRO is often used to address sampling errors of empirical risk minimization, where uncertainty sets are calibrated to offer certain statistical guarantees. For this line of research, we refer the readers to works such as \citet{BHWMR13,EsfahaniKuhn, Lam_burg_entropy,DuchiMOR, Parys_Kuhn}. In risk management, such as the calculation of capital requirement, one typically works with a pre-specified nominal model (e.g., an internal model proposed by the regulator) instead of the empirical distribution. In this case, the goal is to construct uncertainty sets that are sufficiently large to protect against model misspecification.   

In much of the DRO literature, the optimization problem is studied in a risk-neutral environment, where a random loss is evaluated under its expected value. In finance, insurance, and economics, it is more natural to model a decision maker's preferences using a risk measure that is nonlinear in the probability, since humans typically do not perceive changes in probabilities linearly, especially for extreme events. This complicates the choice of an uncertainty set, since in some cases the linearity assumption can be crucial for deriving a tractable reformulation of the primal robust problem. Moreover, the sensitivity of a risk measure to the tails of a distribution also determines the specification of an uncertainty set. For example, an evaluation of the entropic risk measure is only finite for light-tailed distributions. This implies that any uncertainty set containing a heavy-tailed distribution is inadequate for the entropic risk measure. Thus, the natural question arises: given a risk measure and a pre-specified nominal model, how should one properly specify an uncertainty set that is not too conservative (and not too restrictive), and that in addition admits a tractable reformulation?

In this paper, we provide an answer to the above question by proposing a simple, yet unifying framework to specify uncertainty sets for a broad class of risk measures, based on $\phi$-divergences. Introduced by \citet{C75}, $\phi$-divergences are statistical measures for probability models that have received much attention in the DRO literature since the seminal work of \citet{BHWMR13}. Uncertainty sets that are defined as a ball around a nominal model, measured using a $\phi$-divergence, have been extensively studied by \citet{KCS19,KCS21} in the context of model risk assessment. Interestingly, $\phi$-divergences are also strongly connected to risk measures. As shown by \citet{BBT91, BT07}, the four major classes of risk measures: Expected utility, optimized certainty equivalent, shortfall risk measure (a.k.a. u-Mean certainty equivalent), and distortion risk measures, all admit a robust representation with a $\phi$-divergence penalty function.   

Our main contributions can be summarized as follows: by utilizing the pivotal role of $\phi$-divergence in both robust optimization and risk theory, we characterize a wide range of distributionally robust risk measures in terms of two $\phi$-divergences $(\phi_1,\phi_2)$: a $\phi_1$-divergence that specifies the risk measure through an \emph{inner} robust representation, and a $\phi_2$-divergence that controls an \emph{outer} ambiguity set to address model uncertainty. We call this the \emph{composition} approach, since the parametrization with two $\phi$-divergences is formulated as a composite robust optimization problem with two layers of uncertainty. This approach allows us to translate the problem of choosing an uncertainty set for a risk measure into the task of specifying a $\phi_2$ function given $\phi_1$. By deriving a tractable reformulation of the composite robust problem, we show that our $(\phi_1,\phi_2)$ characterization is not only computationally tractable, but also provides a blueprint for how to calibrate the $\phi_2$-divergence uncertainty set to address tail uncertainty, given a risk measure and a nominal model. For risk measures such as the Conditional Value-at-Risk and the entropic risk measure, we provide examples showing the inadequacy of standard divergences to evaluate model risk for certain nominal models, and provide explicit construction of new $\phi_2$-divergences that address these shortcomings, using our suggested framework. In addition, we show that our composition representation also offers other advantages, such as the elicitation of divergences $(\phi_1,\phi_2)$ through queries on decision makers, and natural extension to include a higher level of uncertainty using \textit{globalized robust optimization} (as introduced in \citealt{BenTalGlobal}), that also offers more flexibility in the construction of uncertainty sets. These results serve as extra tools for selecting proper divergences that are customized to a decision maker's risk and ambiguity preferences. 

Finally, we give a brief review of other related literature. \citet{LamMottet} studied the robust optimization approach to address tail uncertainty by constructing uncertainty sets that impose shape constraints on probability densities. However, the tractability of their formulation relies on an explicit representation of the risk measure in terms of density functions, which are not always obtainable for risk measures that are outside the expected utility framework. On the other hand, the tractability of our composite divergences formulation can be easily achieved via convex duality. \citet{BreuerCsiszar} used $\phi$-divergences to measure model risk and provided a tractable reformulation of the corresponding worst-case expectation problem. Building upon this work, 
\citet{KCS19,KCS21} demonstrated how specific $\phi$-divergences can be designed to control the tail behaviors of the distributions in the divergence ball. Our work extends and unifies the results of \citet{BreuerCsiszar} and \citet{KCS19,KCS21} beyond the framework of expected utility models. \citet{GlassermanXU, SchneiderSchweizer} have used the Kullback-Leibler divergence to measure model risk for many financial problems such as option hedging and portfolio credit risk management. \citet{BB09} provided a method for constructing uncertainty sets using coherent risk measures. However, while this approach is theoretically sound, it requires the specification of a decision maker's individual preference for risk, which has to be elicited. Robust risk measures with uncertainty sets that are defined using the Wasserstein distance have also been studied by works such as \cite{Bartl2020} and \citet{Robust_Distortion_Pesenti}, where the focus is more on tractability than calibration of uncertainty sets.

The remaining parts of the paper are organized as follows: Section \ref{sec:Preliminaries} formalizes our two divergences approach in more detail. Section \ref{sec: Divergence_choose} and \ref{sec: Other_Risk_Measures} establish the conditions that we impose on $\phi_2$ for different types of risk measures and provide explicit examples of the construction of new divergences, tailored to a given risk measure and nominal model. Section \ref{sec: OCE_properties} shows how to elicit divergences from robust risk measures, and Section \ref{sec:Compute} discusses the computational aspects of robust risk measures, as well as an extension of our framework using globalized robust optimization to address more uncertainty. Numerical examples of robust options hedging and inventory planning problems are displayed in Section \ref{sec: numerical examples}. A concluding remark is given in Section \ref{sec: conclusion}. Furthermore, all proofs and additional technical details are included in an Electronic Companion.

\section{$\phi$-Divergences and Risk Measures}\label{sec:Preliminaries}

Let $\Phi_0$ denote the set of all non-negative functions $\phi:\mathbb{R}\to [0,\infty)$ that are convex and normalized, i.e., $\phi(1)=0$, with an effective domain $\mathrm{dom}(\phi):=\{t\in \mathbb{R}~|~ \phi(t)<\infty\}\subset [0,\infty)$ which also contains a neighborhood around 1. The $\phi$-divergence between any two measures $\mathbb{P}$ and $\mathbb{Q}$ on an event space  $(\Omega,\mathcal{F})$ is defined as
\begin{equation}\label{def_phi_divergence}
    I_{\phi}(\mathbb{Q},\mathbb{P}):=\begin{cases}
    \int_{\Omega}\phi\left(\frac{\mathrm{d}\mathbb{Q}}{\mathrm{d}\mathbb{P}}\right)\mathrm{d}\mathbb{P}&~\text{if}~\mathbb{Q}\ll \mathbb{P}\\
    +\infty &~\text{else},
    \end{cases}
\end{equation}
where $\mathbb{Q}\ll \mathbb{P}$ denotes absolute continuity with respect to $\mathbb{P}$ (i.e., $\mathbb{P}(A)=0\Rightarrow \mathbb{Q}(A)=0, \forall A\in\mathcal{F}$), and $\frac{\mathrm{d}\mathbb{Q}}{\mathrm{d}\mathbb{P}}$ denotes the Radon-Nikodym derivative of $\mathbb{Q}$ with respect to $\mathbb{P}$. Furthermore, we adopt the conventions: $\phi(0):=\lim_{t\downarrow 0}\phi(t)$, $0\phi \left(0/0\right)=0$, $0\phi \left(a/0\right)=\lim_{\epsilon\downarrow 0}\epsilon\phi \left(a/\epsilon\right)=a\lim_{t\to \infty}\phi(t)/t,~ a>0$. We note that for any probability $\mathbb{Q},\mathbb{P}$, we have $I_{\phi}(\mathbb{Q},\mathbb{P})\geq 0$ and $I_{\phi}(\mathbb{P},\mathbb{P})=0$. 

If $\mathbb{P}$ and $\mathbb{Q}$ have density functions $f,g$, respectively, with respect to some $\sigma$-finite measure $\mathrm{d}x$ on $\mathbb{R}$,\endnote{Such as the Lebesgue measure or the counting measure or a combination of both} then we can also write the $\phi$-divergence in terms of density functions:
\begin{equation*}%\label{phi_div_densities}
I_\phi(g,f)=\int_\mathbb{\mathbb{R}}\phi\left(\frac{g(x)}{f(x)}\right)f(x)\mathrm{d}x.    
\end{equation*}

Let $X: \Omega\to \mathbb{R}$ be a random payoff.\endnote{Negative payoffs are considered as losses} As shown in \citet{ BT87,BT07} and \citet{BBT91}, many types of risk measure/certainty equivalents admit a robust representation with a $\phi$-divergence penalization. This is due to the dual relationship between a utility function and the convex/concave conjugates of $\phi$ (For any $\psi$, we have the convex conjugate $\psi^*(\mathbf{y})\triangleq\sup_{\mathbf{x}\in \mathbb{R}^d}\{\mathbf{y}^T\mathbf{x}-\psi(\mathbf{x})\}$ and the concave conjugate $\psi_*(\mathbf{y})\triangleq-\psi^*(-\mathbf{y})$). For example, the optimized certainty equivalents (OCE) defined with respect to a probability measure $\mathbb{P}$ and a utility function $u$, is equal to the following robust problem where all probability measures $\mathbb{Q}\ll \mathbb{P}$ are penalized by $\phi=-u_*$:
\begin{align}\label{def:OCE}
\begin{split}
   \rho_{\mathrm{oce},\mathbb{P}}(X)&\triangleq\sup_{\{\mathbb{Q}:\mathbb{Q}\ll \mathbb{P}\}}\mathbb{E}_{\mathbb{Q}}[-X]-I_{\phi}(\mathbb{Q},\mathbb{P})\\
   &=\inf_{\eta\in \mathbb{R}}\{\eta - \mathbb{E}_{\mathbb{P}}[u(X+\eta)]\}.
\end{split}
\end{align} 
Therefore, the OCE is also known as the \textit{divergence risk measure}, which includes examples such as:
\begin{itemize}
    \item Conditional Value-at-Risk ($\mathrm{CVaR}_{\alpha}$) at $\alpha\in (0,1)$, where $\phi(t)= \delta_{[0,1/(1-\alpha)]}(t)$ is the indicator function on $[0,1/(1-\alpha)]$ \endnote{For a set $\mathcal{C}\in \mathbb{R}^n$, the indicator function is defined as $\delta_{\mathcal{C}}(x)=\begin{cases}
    0& x\in \mathcal{C}\\
    +\infty& x\notin \mathcal{C}.
\end{cases}$ \label{def:indicator} }. \label{def:CVaR}
    \item Entropic risk measure $\rho_{e,\gamma}(X)= \log \left(\mathbb{E}[e^{-\gamma X}]\right)/\gamma$, $\gamma>0$, where $\phi(t)=(t\log(t)-t+1)/\gamma$, for $t\geq 0$. \label{def:entropy_riskmeasure}
\end{itemize}
The shortfall risk measures, a.k.a. u-Mean certainty equivalents, can also be represented by a $\phi$-divergence:
\begin{align}\label{u-Mean}
\begin{split}
        \rho_{\mathrm{sf},\mathbb{P}}(X)&\triangleq\sup_{\lambda >0}\sup_{\{\mathbb{Q}:\mathbb{Q}\ll \mathbb{P}\}}\mathbb{E}_{\mathbb{Q}}[-X]-I_{\phi^\lambda}(\mathbb{Q},\mathbb{P})\\
    &=\inf_{\eta\in \mathbb{R}}\{\eta~|~\mathbb{E}_{\mathbb{P}}[-u(X+\eta)]\leq 0\},
\end{split}
\end{align}
where $\phi^\lambda(t)=\lambda \phi\left(t/\lambda\right)$ and $\phi(t)=-u_*(t)$. 
%\begin{align}\label{u-Mean}
%    \rho_{\mathrm{sf},\mathbb{P}}(X)=\sup_{\lambda >0}\sup_{\{\mathbb{Q}:\mathbb{Q}\ll \mathbb{P}\}}\mathbb{E}_{\mathbb{Q}}[-X]-I_{\phi^\lambda_1}(\mathbb{Q},\mathbb{P}),
%\end{align}

The expected utility risk measures can also be expressed using a $\phi$-divergence. The only difference is that the maximization is taken over all finite measures, instead of only probability measures:
\begin{equation}\label{EU_dual_representation}
    \mathbb{E}_{\mathbb{P}}[-u(X)]=\sup_{\{\mathbb{Q}:\mathbb{Q}\ll \mathbb{P}, \mathbb{Q}(\Omega)<\infty\}}\mathbb{E}_{\mathbb{Q}}[-X]-I_{\phi}(\mathbb{Q},\mathbb{P}),
\end{equation}
where again $\phi=-u_*$. 

Therefore, we have seen that many types of risk measures admit a robust representation in terms of $\phi$-divergences. Adding a second $\phi_2$-divergence, we can thus parametrize a broad class of robust risk measures by a composition of $\phi$-divergences. We can address ambiguity in two ways: Either by means of the \emph{multiple prior} approach via an uncertainty set (which we often denote using $\rho^s$), or by means of the \emph{multiplier preference} approach through penalization (which we often denote using $\rho^l$). In the case of OCE risk measures, this leads to the following composite robust formulations, where the maximization is taken over two variables $(\mathbb{Q},\bar{\mathbb{Q}})$:
\begin{align}\label{intro_rob_measures}
\begin{split}
        \rho^{s}_{\mathrm{oce},\mathbb{P}_0}(X)
        &\triangleq\sup_{\{\mathbb{Q}:I_{\phi_2}(\mathbb{Q},\mathbb{P}_0)\leq r\}}\sup_{\{\bar{\mathbb{Q}}:\bar{\mathbb{Q}}\ll \mathbb{Q}\}}\mathbb{E}_{\bar{\mathbb{Q}}}[-X]\\
        &\qquad-I_{\phi_1}(\bar{\mathbb{Q}},\mathbb{Q}),
\end{split}
\end{align}
and
\begin{align}\label{intro_rob_measures_pen}
\begin{split}
     \rho^{l}_{\mathrm{oce},\mathbb{P}_0}(X)
     &\triangleq \sup_{\{\mathbb{Q}:\mathbb{Q}\ll \mathbb{P}_0\}}\sup_{\{\bar{\mathbb{Q}}:\bar{\mathbb{Q}}\ll \mathbb{Q}\}}\mathbb{E}_{\bar{\mathbb{Q}}}[-X] \\
     &\qquad-I_{\phi_1}(\bar{\mathbb{Q}},\mathbb{Q})- I_{\phi_2}(\mathbb{Q},\mathbb{P}_0).
\end{split}
\end{align}
Other classes of robust risk measures are defined similarly. We note that for any constant $\lambda>0$, we have $\lambda I_{\phi}(\mathbb{Q},\mathbb{P})=I_{\lambda \phi}(\mathbb{Q},\mathbb{P})$. Hence, the representations in \eqref{intro_rob_measures} and \eqref{intro_rob_measures_pen} also include penalization constants that are not equal to one. Throughout this paper, we fix $\mathbb{P}_0$ as a \emph{given} nominal distribution. This can be a distribution that is already specified in an economic model, or a distribution calibrated from a parametric family. Furthermore, we assume that $X\in L^1(\mathbb{P}_0)$ (i.e., $\int |X|\mathrm{d}\mathbb{P}_0<\infty)$. In the following section, we extensively study the robust OCE risk measures as the prime example of this paper. We show that by reformulating the problems \eqref{intro_rob_measures} and \eqref{intro_rob_measures_pen} into finite-dimensional dual problems, we can obtain a blueprint for constructing $\phi_2$ that is specifically tailored to a given risk measure $\phi_1$ and the nominal model $\mathbb{P}_0$.
\section{Reformulation of Robust OCE Risk Measures and Its Implications for Divergence Choices} \label{sec: Divergence_choose}

To illustrate the importance of choosing a proper $\phi$-divergence for a robust risk measure, we consider an example of a divergence function defined by the Dirac delta function $\phi_{\alpha}(t)=\delta_{[0,1/(1-\alpha)]}(t)$, for some $\alpha\in (0,1)$. Then, any alternative distributions with a density function $g$ that has a heavier tail than the nominal density $f_0$ (i.e., $\lim_{x\to \infty}g(x)/f_0(x)=\infty$) will have an infinite divergence value with respect to $\phi_\alpha$. Indeed, the likelihood ratio $g(x)/f_0(x)$ exceeds $1/(1-\alpha)$ for sufficiently large values of $x$, eventually lying outside the effective domain of $\phi_{\alpha}$. Therefore, the $\phi$-divergence uncertainty set induced by $\phi_\alpha$ does not include any distributions that have a heavier right-side tail than the nominal distribution, while heavy-tailed losses are often the real threats in financial risk management that one wishes to be robust against. 

The choice of $\phi$-divergence also depends on the sensitivity of a risk measure to tail losses. For example, the entropic risk measure $\rho(X)=\log(\mathbb{E}[e^{-X}])$ is only finite for distributions where the tail of their density functions decays faster than $e^{-|x|}$. However, if the modified chi-squared divergence $\phi_2(x)=(x-1)^2$ is used to construct an uncertainty set where the nominal model for $Y:=-X$ is the exponential distribution $f_0(y)= \lambda e^{-\lambda y}$ for $ y\geq 0, \lambda>1$, then this can lead to an infinitely pessimistic robust risk evaluation. Indeed, any exponential distribution $g(y)=\eta e^{-\eta y}$ with $\eta<1<\lambda<2\eta$ has a finite modified chi-squared divergence with respect to $f_0$. Therefore, the modified chi-squared divergence ball will contain the distribution $g$ for sufficiently large radius $r$. However, since $1-\eta >0$, the entropic risk measure is infinite under the distribution of $g$. In fact, for any radius $r>0$, one can construct a density $\tilde{g}$, that has the same asymptotic tail behavior as $e^{-\eta y}$, and such that $I_{\phi_2}(\tilde{g},f_0)<r$ (see \citealt{KCS21}). Therefore, in this particular example, the modified chi-squared divergence ball will always lead to an infinite robust entropic risk measure evaluation, for all radius $r>0$. 

Motivated by these examples, it is thus important to establish the necessary and sufficient conditions that one must impose on $\phi_2$ to ensure a finite evaluation of the robust risk measures $\rho^{s}_{\mathrm{oce},\mathbb{P}_0}$ and $\rho^{l}_{\mathrm{oce},\mathbb{P}_0}$. The key to identifying these conditions is by examining the duals of \eqref{intro_rob_measures} and \eqref{intro_rob_measures_pen}. For any divergence function $\phi$, we define $0\phi^*\left(s/0\right):=0$ if $s\leq 0$ and $0\phi^*\left(s/0\right):=\infty$ if $s>0$.
\begin{theorem}\label{thm:ref_rob_OCE_pen}
Let $X\in L^1(\mathbb{P}_0)$ and $\phi_1,\phi_2\in \Phi_0$ be lower-semicontinuous. Then, we have the following equalities:
\begin{align}
\begin{split}
    \rho^{l}_{\mathrm{oce},\mathbb{P}_0}(X)&=\inf_{\theta_1,\theta_2\in \mathbb{R}}-\theta_1-\theta_2\\
   &\qquad+\mathbb{E}_{\mathbb{P}_0}\left[\phi_2^*(\phi_1^*(\theta_2-X)+\theta_1)\right]\label{dual_OCE_p}
\end{split}
\end{align}
  \begin{align}
  \begin{split}
     \rho^{s}_{\mathrm{oce},\mathbb{P}_0}(X)&=\inf_{\substack{\lambda\geq 0\\\theta_1,\theta_2\in \mathbb
    R}}-\theta_1-\theta_2+\lambda r\\
    &\qquad+\mathbb{E}_{\mathbb{P}_0}\left[\lambda \phi_2^*\left(\frac{\phi_1^*(\theta_2-X)+\theta_1}{\lambda}\right)\right].\label{dual_OCE_b} 
  \end{split}
\end{align}
Furthermore, if the following conditions are satisfied:
\begin{enumerate}%[label=(\arabic*)]
    \item $\phi_2^*$ and $\phi_1^*$ are differentiable on $\mathbb{R}$,
    \item Each domain set $\mathcal{F}(\boldsymbol{\theta})\triangleq\{\boldsymbol{\theta}\in \mathbb{R}^2~|~ \int_\Omega \beta^*(\omega,\boldsymbol{\theta})\mathrm{d}\mathbb{P}_0(\omega)<\infty\}$ and $\tilde{\mathcal{F}}(\boldsymbol{\theta},\lambda)\triangleq\{(\boldsymbol{\theta},\lambda)\in \mathbb{R}^2\times [0,\infty)~|~ \int_\Omega \tilde{\beta}^*(\omega,\boldsymbol{\theta},\lambda)\mathrm{d}\mathbb{P}_0(\omega)<\infty\}$ has non-empty interior, where $\beta^*(\omega,\boldsymbol{\theta})\triangleq \phi_2^*(\phi_1^*(\theta_2-X(\omega))+\theta_1)$ and $\tilde{\beta}^*(\omega,\boldsymbol{\theta},\lambda)\triangleq\lambda \phi_2^*\left(\frac{\phi_1^*(\theta_2-X(\omega))+\theta_1}{\lambda}\right)$.
\end{enumerate}
Then, there exists dual solutions $\boldsymbol{\theta}^*$ and $(\tilde{\boldsymbol{\theta}}^*,\lambda^*)$ that attain respectively the minimum \eqref{dual_OCE_p} and \eqref{dual_OCE_b}. Moreover, if $\boldsymbol{\theta}^*\in \mathrm{int}(\mathcal{F}(\boldsymbol{\theta}))$ and $(\tilde{\boldsymbol{\theta}}^*,\lambda^*)\in \mathrm{int}(\tilde{\mathcal{F}}(\boldsymbol{\theta},\lambda))$, then the partial derivatives of $\beta^*(\omega,.)$ and $\tilde{\beta}^*(\omega,.,\lambda)$ with respect to the variable $\boldsymbol{\theta}$,
evaluated at the dual solutions $\boldsymbol{\theta}^*$ and $(\tilde{\boldsymbol{\theta}}^*,\lambda^*)$, are respectively the worst-case densities of the measure 
 $\mathbb{Q}^*, \bar{\mathbb{Q}}^*$ with respect to $\mathbb{P}_0$ that attain the maximum of the primal problems \eqref{intro_rob_measures_pen} and \eqref{intro_rob_measures}.
\end{theorem}

As a direct consequence of Theorem \ref{thm:ref_rob_OCE_pen}, we derive the following necessary and sufficient conditions for the finiteness of $\rho^{s}_{\mathrm{oce},\mathbb{P}_0}(X)$ and $\rho^{l}_{\mathrm{oce},\mathbb{P}_0}(X)$. In the following, we assume that the distribution of $X$ under the nominal measure $\mathbb{P}_0$ has a real-valued density function $f_0(x)$.
\begin{corollary}\label{cor:finite_oce}
    Let $X\in L^1(\mathbb{P}_0)$ and $\phi_1,\phi_2\in \Phi_0$ be lower-semicontinuous. Then, we have that $\rho^{s}_{\mathrm{oce},\mathbb{P}_0}(X)<\infty$, if and only if there exists $\theta_1,\theta_2\in \mathbb{R}, \lambda\geq 0$, such that
    \begin{equation}\label{oceb_integral}
   \int_{\mathbb{R}} \lambda \phi_2^*\left(\frac{\phi_1^*(\theta_2-x)+\theta_1}{\lambda}\right)f_0(x)\mathrm{d}x<\infty.
    \end{equation}
    Similarly,  we have that $\rho^{l}_{\mathrm{oce},\mathbb{P}_0}(X)<\infty$, if and only if there exists $\theta_1,\theta_2\in \mathbb{R}$, such that
    \begin{equation}\label{ocep_integral}
   \int_{\mathbb{R}} \phi_2^*\left(\phi_1^*(\theta_2-x)+\theta_1\right)f_0(x)\mathrm{d}x<\infty.
    \end{equation}
\end{corollary}
The integral conditions \eqref{oceb_integral} and \eqref{ocep_integral} displayed in Corollary \ref{cor:finite_oce} show that the finiteness of a robust OCE risk measure is completely determined by the tail behaviors of three components $(\phi_2^*, \phi_1^*, f_0)$, where $\phi_2^*, f_0$ control the content of the $\phi_2$-uncertainty set, and $\phi_1^*$ encodes the sensitivity of the risk measure to losses in the tail. This provides a guideline on how to calibrate an uncertainty set to a given risk measure and nominal model. Namely, one must choose the $\phi_2$-divergence function such that its conjugate $\phi_2^*$ leads to a finite integral in \eqref{oceb_integral} or \eqref{ocep_integral}.  

We note that the conditions in Corollary \ref{cor:finite_oce} are only verifiable if we have information on the tail behavior of the nominal density $f_0$. In practice, $X$ might depend on many underlying risk factors $Z_1,\ldots,Z_I$, where only the nominal distributions of the marginals are specified and have an explicit form. In the following proposition, which is adapted from \citet{KCS19}, we provide a sufficient condition under which one can choose the divergences based on the nominal distributions of the marginals $Z_i$'s. 

\begin{proposition}\label{prop: Risk factors}
    Let $X\in L^1(\mathbb{P}_0)$. Suppose there exists a constant $C>0$, such that 
    \begin{equation}\label{bound_risk_factor}
        |X|\leq C\left(1+\sum^m_{i=1}|Z_i|\right),
    \end{equation}
holds $\mathbb{P}_0$-almost surely. 
If there exists $\theta_1,\theta_2\in \mathbb{R}$, such that for all $i=1,\ldots, I$,
\begin{equation*}
   \mathbb{E}_{\mathbb{P}_0}\left[ \phi_2^*\left(\theta_1+\phi_1^*(\theta_2+C(1+m\cdot |Z_i|))\right)\right]<\infty.
\end{equation*}
Then, the integral conditions \eqref{oceb_integral} and \eqref{ocep_integral} are also satisfied.
\end{proposition}

Using the guideline given in Corollary \ref{cor:finite_oce}, we perform several tail analyses and illustrate in a couple of examples how most canonical choices of $\phi$-divergences are inadequate to address tail uncertainty for some standard risk measures and nominal distributions. Throughout this paper, we use the big O notation $f_1(x)=\mathcal{O}(f_2(x))$ if $\limsup_{x\to \infty}|f_1(x)|/|f_2(x)|<\infty$, and the small o notation $f_1(x)=o(f_2(x))$ if $\lim_{x\to \infty}|f_1(x)|/|f_2(x)|=0$.

\begin{example}\label{example:KL_lognormal}
Consider the $\mathrm{CVaR}_{\alpha}$ risk measure with a log-normal nominal model. A common choice of a $\phi$-divergence is the Kullback-Leibler divergence, where $\phi(t)=t\log(t)-t+1$. Its conjugate is given by $\phi^*(s)=e^s-1$. The $\mathrm{CVaR}_{\alpha}$ risk measure \eqref{def:CVaR} is an OCE risk measure with $\phi^*_1(s)=\max\{s/(1-\alpha),0\}$. However, examining the integral condition \eqref{oceb_integral} reveals that for any $\theta_1,\theta_2\in \mathbb{R}, \lambda>0$,
%\begin{align*}
%    \int^0_{-\infty}\left(e^{\frac{1}{\lambda}(\max\{\frac{1}{1-\alpha}(\theta_2-x),0\}+\theta_1)}-1\right)\frac{1}{|x|\sigma \sqrt{2\pi}}e^{-\frac{(\log(|x|)-\mu)^2}{2\sigma^2}}\mathrm{d}x=+\infty,
%\end{align*}
the integrand has a tail behavior of $\mathcal{O}\left(\exp\left\{|x|/(\lambda\alpha)\right\}\right)$, 
%\begin{equation*}
%    \exp\left({\frac{1}{\lambda\alpha}|x|-\log(|x|)-\frac{(\log(|x|)-\mu)^2}{2\sigma^2}}\right),
%\end{equation*}
which diverges to $+\infty$ as $x\to -\infty$. Hence, the uncertainty set induced by the Kullback-Leibler divergence is too conservative for the $\mathrm{CVaR}_\alpha$ risk measure under a log-normal nominal model. In Table \ref{tab:CVaR_finite}, we give an overview of the finiteness status of robust $\mathrm{CVaR}$ with other nominal distributions and $\phi$-divergences.
\end{example}
\begin{example}
    The $\chi^2$-distance is another popular $\phi$-divergence where $\phi(t)=(t-1)^2/t$. However, its conjugate $\phi^*(s)=2-2\sqrt{1-s}, s<1$ is only finite on the domain on $(-\infty,1)$, i.e., $\phi^*(s)=+\infty$, for $s>1$. Hence, the integral condition \eqref{oceb_integral} can never be satisfied if the $\chi^2$-distance is used to define a robust OCE risk measure, for any nominal distribution that has unbounded support towards $-\infty$. 
    
    Similarly, other canonical examples of $\phi$-divergences such as Burg entropy $\phi(t)=-\log(t)+t-1$, total variation distance $\phi(t)=|t-1|$, and polynomial divergence $\phi(t)=(t^p-p(t-1)-1)/(p(p-1))$ with degree $p<1$, all have a conjugate function $\phi^*$ that is infinite on certain interval $[a,\infty)$ (see Table 2 of \citealt{BHWMR13} and page 6 of \citealt{Pardo06}). Hence, they all lead to an infinite robust OCE risk evaluation for a nominal distribution with unbounded support towards $-\infty$. 
\end{example}
\begin{example}\label{example: KLentropic_gaussian}
 The entropic risk measure $\rho_{e,\gamma}(X)=\log \left(\mathbb{E}[e^{-\gamma X}]\right)/\gamma$ is an OCE risk measure with $\phi_1^*(x)=(e^{\gamma x}-1)/\gamma$. Due to its exponential growth behavior, the entropic risk measure has infinite evaluation for heavy-tailed distributions and therefore also its robust counterpart. For the Gaussian distribution, the entropic risk measure is finite and equal to $\mu+\gamma \sigma^2/2$. However, if one chooses the KL-divergence to define a robust entropic risk measure with a Gaussian nominal model, then condition \eqref{oceb_integral} is not satisfied for any $\theta_1,\theta_2\in \mathbb{R}, \lambda>0$, since the integrand has a tail asymptotic of the order
\begin{equation*}
   \mathcal{O}\left(\exp\left\{\frac{1}{\lambda}\left(\frac{1}{\gamma}\left(e^{\gamma(\theta_2-x)}-1\right)+\theta_1\right)-\left(\frac{x-\mu}{\sigma}\right)^2\right\}\right),
\end{equation*}
which diverges to $\infty$ as $x\to -\infty$. Therefore, the KL-divergence is not adequate to address model uncertainty for the entropic risk measure, under many nominal distributions such as the Gaussian and Weibull distribution, as shown in Table \ref{tab:exp_OCE_finite}. 
\end{example}

\begin{table}[t]
    \centering
    \begin{tabular}{cccccc}
 
        Divergence&Gaussian&Weibull&Log-normal&Pareto&Student t\\ \hline    
         Kullback-Leibler&$<\infty$&$*$&$\infty$&$\infty$&$\infty$  \\ 
         Polynomial, $p>1$&$<\infty$&$<\infty$&$<\infty$&$*$ &$*$\\
         Polynomial, $0<p<1$&$\infty$&$\infty$&$\infty$&$\infty$&$\infty$\\
    \end{tabular}
    \caption{Finiteness status of robust $\mathrm{CVaR}_{\alpha}$ risk measures \eqref{intro_rob_measures} and \eqref{intro_rob_measures_pen} for different divergence measures versus nominal distributions. $*$ denotes that the finiteness status depends on specific choice of parameters.}
    \label{tab:CVaR_finite}
\end{table}
 \begin{table}[t]
    \centering
    \begin{tabular}{cccccc}

        Divergence&Gaussian&Weibull&Log-normal&Pareto&Student t\\ \hline    
         Kullback-Leibler&$\infty$&$\infty$&$\infty$&$\infty$ &$\infty$  \\ 
         Polynomial, $p>1$&$ <\infty$&$*$&$\infty$&$\infty$&$\infty$\\
         Polynomial, $0<p<1$&$\infty$&$\infty$&$\infty$&$\infty$&$\infty$\\
    \end{tabular}
    \caption{Finiteness status of robust entropic risk measures \eqref{intro_rob_measures} and \eqref{intro_rob_measures_pen} for different divergence measures versus nominal distributions, with exponential utility OCE. $*$ denotes that the finiteness status depends on specific choice of parameters.}
    \label{tab:exp_OCE_finite}
\end{table}
From these examples, we observe that many standard choices of $\phi$-divergence uncertainty sets often lead to infinitely conservative robust risk evaluations. On the other hand, the polynomial divergence $\phi(t)=(t^p-p(t-1)-1)/(p(p-1))$ can often be made to satisfy the finiteness conditions \eqref{oceb_integral} and \eqref{ocep_integral}, by choosing a sufficiently large $p$. Indeed, an examination of its conjugate function $\phi^*(s)=\max\left\{1+s(p-1),0\right\}^{\frac{p}{p-1}}/p-1/p$ shows that as $p\to \infty$, $\phi^*$ approaches the linear function $\phi^*(s)=s$, which is the least conservative conjugate function of the $\phi$-divergence that only considers the nominal model (i.e., $\phi(1)=0$, $\phi(t)=+\infty$ elsewhere). However, the drawback of a polynomial divergence uncertainty set is that it can be very restrictive. For example, consider two Weibull distributions $g(x)= kx^{k-1}e^{-x^k}$ and $f(x)= lx^{l-1}e^{-x^l}$ with $k<l$. Then, any polynomial divergence $\phi(t)= \mathcal{O}(t^p)$ with $p>1$ would give $I_{\phi}(g,f)=\infty$, since $\phi\left(g(x)/f(x)\right)f(x)=\mathcal{O}(x^{k-l}e^{(p-1)x^l-px^k})$. Therefore, when the nominal distribution belongs to the Weibull class, an uncertainty set induced by the polynomial divergence does not include any other Weibull distribution that has slightly heavier tail. Hence, apart from the conditions in Corollary \ref{cor:finite_oce}, we need to impose additional constraints on the conjugate function $\phi_2^*$, such that the resulting uncertainty set not only guarantees finiteness of a robust OCE risk measure, but is also not too restrictive. We will address this issue in the next subsection.

\subsection{Calibrating The Content of a Divergence Ball Through The Convex Conjugate}\label{sec:phi_construct_conjugate}
As illustrated by previous examples, many classical divergences fail to guarantee the finiteness of a robust risk measure, due to improper matching of tail behaviors between $\phi_2^*$, $\phi_1^*$ and $f_0$. Therefore, to construct an alternative $\phi_2$-divergence that addresses this issue, we can specify a new tail function $\tilde{\psi}(s), s\geq 0$ that is compatible with $(\phi_1^*,f_0)$ according to the conditions established in Corollary \ref{cor:finite_oce}, and then set $\phi_2^*(s)=\tilde{\psi}(s)$, for $s\geq 0$. Since the integral conditions in Corollary \ref{cor:finite_oce} depend only on the tail behavior of $\phi_2^*(s)$ as $s\to \infty$, we may simply set $\phi_2^*(s)=e^s-1$ for $s\leq 0$. Finally, the new tail function $\tilde{\psi}(s)$ must be normalized such that $(\phi_2^*)^*$ becomes a $\phi$-divergence function (i.e., a convex, non-negative function satisfying $\phi(1)=0$ and $\mathrm{dom}(\phi)\subset [0,\infty)$). This can be done by choosing $\tilde{\psi}(s)$ to be increasing and convex with $\tilde{\psi}(0)=0$ and $\tilde{\psi}'(0)=1$. To conclude, we propose the following construction: let $\psi$ be an increasing, convex differentiable function on $[0,\infty)$ such that $\psi''(0)\neq 0$. Set, for $s\geq 0$, 
\begin{equation}\label{def_psi}
    \tilde{\psi}(s)=\frac{1}{\psi''(0)}\left(\psi(s)+\left(\psi''(0)-\psi'(0)\right)s-\psi(0)\right),
\end{equation}
and define the function
\begin{equation}\label{def_phi^*}
    \phi^*_2(s)=\begin{cases}
        \tilde{\psi}(s)& s\geq 0\\
        e^s-1& s\leq 0.
    \end{cases}
\end{equation}
The following proposition, which is an application of Corollary 23.5.1 of \citet{Rockafellar}, verifies that $\phi_2^*$, as defined in \eqref{def_phi^*}, is indeed the conjugate function of a $\phi$-divergence.
\begin{proposition}\label{prop: construct_phi_conj}
    The function $\phi_2:=(\phi_2^*)^*$, where $\phi_2^*$ is defined in \eqref{def_phi^*}, is a convex, non-negative function on $[0,\infty)$ with $\phi_2(1)=0$. If $\psi$ in \eqref{def_psi} is strictly convex and continuously differentiable with $\lim_{s\to \infty}\psi'(s)=\infty$, then $\phi_2$ is finite on $[0,\infty)$ and we have that $\phi_2'(t)=((\phi_2^*)')^{-1}(t)$, for $t>0$. 
\end{proposition}
   
Besides ensuring the finiteness of a robust risk measure, it is also desirable to have a divergence that is not too restrictive; for example, one that allows for distributions with heavier tails than the nominal one in the divergence ball. Intuitively, the degree of penalization of $\phi_2$ is determined by its growth rate as $t\to \infty$, which is governed by the derivative $\phi_2'(t)$. Therefore, one must control $\phi_2'(t)$ to adjust the content of a divergence ball. However, if the divergence is constructed through the conjugate as in \eqref{def_phi^*}, then there is not always an explicit form of $\phi_2'$ available that allows us to examine its growth behavior as $t\to\infty$. Indeed, although the equality $\phi_2'=((\phi_2^*)')^{-1}$ holds, the inverse of $(\phi_2^*)'$ cannot always be calculated explicitly. Therefore, it would seem more convenient to specify $\phi_2$ through its derivative instead of the conjugate, and to derive sufficient conditions for \eqref{oceb_integral} and \eqref{ocep_integral} in terms of $\phi_2'$ (e.g., by bounding the conjugate using $\phi_2^*(s)\leq s(\phi_2')^{-1}(s)$, similarly as in \citealt{KCS19}). This creates a dilemma, since the latter construction does not give an analytical expression for the conjugate function and its derivatives, which is important for the use of optimization techniques to compute the robust risk measures via their dual problems, as we elaborate in Section \ref{sec:Compute}. Nevertheless, we can resolve this dilemma by bounding $(\phi_2^*)'$ by other functions that have similar tail properties, but for which the inverse of the derivatives can be explicitly determined. This provides bounds on $((\phi_2^*)')^{-1}=\phi_2'$, and therefore bypasses the need for an explicit form of $\phi_2'$ to analyze the content of a $\phi_2$-divergence ball. We state this more precisely in the following proposition, which is adapted from \citet{KCS19}, and allows us to control the moments of the distributions in a $\phi_2$-divergence ball, through the derivative of $\phi_2^*$.   

\begin{proposition}\label{prop: phi_bound}
    Let $\phi_2$ be a continuous divergence function with a strictly increasing derivative $\phi_2'$. Suppose $\psi_1'$ is a continuous function such that its inverse satisfies $(\psi_1')^{-1}(y) \leq (\phi_2')^{-1}(y)$ for all $y\geq y_0$, for some $y_0\in \mathbb{R}$. If for some $d>1$, we have that
    \begin{equation}\label{moment_include}
        \limsup_{|x|\to \infty}\psi_1'\left(\frac{1}{f_0(x)}\right)|x|^{-d}<\infty,
    \end{equation} 
    then for any probability density $g$ such that the likelihood $g/f_0$ is bounded on any compact subset of $(-\infty,0]$,  we have that if $g$ has a finite $d$-th moment: $\int^\infty_{-\infty}g(x)|x|^d\mathrm{d}x<\infty$, then we also have $I_{\phi_2}(g,f_0)<\infty$.
    
    On the other hand, if there exists a divergence function $\psi_2\in \Phi_0$ such that $(\phi_2')^{-1}(y)\leq (\psi_2')^{-1}(y)$ for all $y\geq 0$, we have that $I_{\phi_2}(g,f_0)=\infty$, if $I_{\psi_2}(g,f_0)=\infty$, for any density function $g$.
\end{proposition}

\subsection{Examples of Tailored Divergences}\label{sec: examples_tailored_div}

In this subsection, we provide examples of new $\phi_2$-divergences that are tailored to a given risk measure and nominal model. We will only specify $\phi_2^*(s)$ for $s\geq 0$ in these examples, since we always set $\phi_2^*(s)=e^s-1$ as in \eqref{def_phi^*}. 
\begin{example}
Consider the $\mathrm{CVaR}_{\alpha}$ risk measure with the generalized log-normal nominal distribution, where the density function is given by 
\begin{equation*}
    f_0(x)= \frac{1}{C\sigma |x|}\exp\left\{-\frac{1}{p\sigma^p}(\log|x|-\mu)^p\right\},~x\leq 0,
\end{equation*}
where $\mu\in \mathbb{R}, \sigma>0, p\geq 2$, and $C>0$ is a normalization constant. As shown in Example \ref{example:KL_lognormal}, the KL-divergence is not suitable for this combination of nominal model and risk measure. Hence, using the construction method in \eqref{def_phi^*}, we propose the following divergence specified by its conjugate function:
\begin{align}\label{div_gl_cvar}
\begin{split}
     \phi_{2}^*(s)&=
          c_1(s+e)\exp\left\{\frac{1}{p\left(\sigma d\right)^p}\log^p(s+e)\right\}\\
          &\qquad+c_2s+c_3,~ s\geq 0, 
\end{split}
\end{align}
where the constants are $c_1=1/(p^2(a^2+a)e^{a-1}), c_2=1-e^a(ap+1)/c_1, c_3= -e^{a+1}/c_1$ with $a=1/(p\left(\sigma d\right)^p)$ and $d>1$. We note that the proposed divergence \eqref{div_gl_cvar} contains the parameters $(\sigma,p)$, due to its dependence on the log-normal nominal model. It also contains a parameter $d>1$, which controls the moments of the distributions in its divergence ball. Notably, the divergence \eqref{div_gl_cvar} is independent of the parameter $\alpha$ of the Conditional Value-at-Risk, which is due to the linearity of the $\phi_1^*$ function corresponding to the $\mathrm{CVaR}_{\alpha}$ risk measure, that does not determine the asymptotic behavior of the integrand in \eqref{ocep_integral}. 
\end{example}
\begin{lemma}\label{lem: div_gl_cvar}
     The function $\phi_2^*$ defined in \eqref{div_gl_cvar} is strictly convex and twice continuously differentiable on $\mathbb{R}$. Its conjugate $\phi_2:=(\phi_2^*)^*$ belongs to $\Phi_0$. Moreover, $\phi_2^*$ satisfies \eqref{oceb_integral} and \eqref{ocep_integral}, for the $\mathrm{CVaR}_{\alpha}$ risk measure and the generalized log-normal nominal model.

     In addition, for any continuous probability density $g$ defined on $(-\infty,0]$ that has finite $d$-th moment: $\int^0_{-\infty}g(x)|x|^d<\infty$, we have that $I_{\phi_2}(g,f_0)<\infty$. Conversely, for any density function $g$ such that $\liminf_{x\to -\infty}|x|^{t+1}g(x)>0$, where $t\in (1,d)$, we have $I_{\phi_2}(g,f_0)=\infty$.
\end{lemma}

\begin{example}
Consider now the $\mathrm{CVaR}_{\alpha}$ risk measure with a Weibull nominal distribution, where the density $f_0$ is given by:
\begin{equation*}
     f_0(x)= \frac{k}{\lambda}\left(\frac{|x|}{\lambda}\right)^{k-1}\exp\left\{{-\left(\frac{|x|}{\lambda}\right)^k}\right\},~ x\leq 0,
\end{equation*}
where $\lambda, k>0$.
If $k<1$, the Weibull distribution becomes heavy-tailed, and the robust $\mathrm{CVaR}_{\alpha}$ risk measure defined by the KL-divergence will lead to an infinite evaluation. To match the exponential tail of $f_0$, we propose the following divergence:
\begin{align}\label{div:exp_power}
\begin{split}
       \phi_2^*(s)&=c_1(s+1)\exp\{{(s+1)^{k/d}}\}\\
       &\qquad+c_2s+c_3, ~s\geq 0.
\end{split}
\end{align}
where $d>1$, and the normalization constants are $c_1 = 1/(ep(2p+1))$, $c_2 = 1-(1+p)/(2p+1), c_3 =1/(p(2p+1))$ with $p=k/d$.
\end{example}
\begin{lemma}\label{lem: exp_power}
    The function $\phi_2^*$ defined in \eqref{div:exp_power} is strictly convex and twice continuously differentiable on $\mathbb{R}$. Its conjugate $\phi_2:=\phi_{2}^{**}$ belongs to $\Phi_0$. Moreover, $\phi_2^*$ satisfies the finiteness conditions \eqref{oceb_integral} and \eqref{ocep_integral}, for the $\mathrm{CVaR}_{\alpha}$ risk measure and the Weibull nominal model. 

     In addition, for any continuous probability density $g$ defined on $(-\infty,0]$ that has finite $d$-th moment: $\int^0_{-\infty}g(x)|x|^d<\infty$, we have that $I_{\phi_2}(g,f_0)<\infty$. Conversely, for any density function $g$ such that $\liminf_{x\to -\infty}|x|^{t+1}g(x)>0$, where $t\in (1,d)$, we have $I_{\phi_2}(g,f_0)=\infty$.
\end{lemma}

\begin{example}    
Consider the entropic risk measure $\rho_{e,\gamma}$ with the Weibull nominal distribution, where $k>1$. Table \ref{tab:exp_OCE_finite} shows that the KL-divergence is generally not suitable for the entropic risk measure. Moreover, as mentioned previously, the polynomial divergence is too restrictive since its uncertainty set does not include any Weibull distributions with a heavier tail than the nominal one. Therefore, we propose the following divergence:  
\begin{align}\label{div_entr_weibull}
\begin{split}
     \phi_{2}^*(s)&=
          c_1(s+e)\exp\left\{\frac{1}{\left(2\gamma\lambda \right)^k}\log^k(s+e)\right\}\\
          &\qquad+c_2s+c_3, ~s\geq 0,
\end{split}
\end{align}
where the constants are $c_1=1/(k^2(a^2+a)e^{a-1}), c_2=1-e^a(ak+1)/c_1, c_3= -e^{a+1}/c_1$ with $a=1/\left(2\gamma\lambda\right)^k$. 
\end{example}
\begin{lemma}\label{lem: entr_weibull}
    The function $\phi_2^*$ defined in \eqref{div_entr_weibull} is strictly convex and twice continuously differentiable on $\mathbb{R}$. Its conjugate $\phi_2:=\phi_{2}^{**}$ belongs to $\Phi_0$. Moreover, $\phi_2^*$ satisfies the finiteness conditions \eqref{oceb_integral} and \eqref{ocep_integral}, for the entropic risk measure (with parameter $\gamma>0$) and the Weibull nominal distribution $f_0$ with parameters $(\lambda,k)$. Furthermore, for any Weibull density function $g$ with parameters $(\lambda,l)$ such that $1<l<k$, we have $I_{\phi_2}(g,f_0)<\infty$.
\end{lemma}
\section{Other Types of Robust Risk Measures} \label{sec: Other_Risk_Measures}
In this section, we show that the formulation with two $\phi$-divergences $(\phi_1, \phi_2)$ provides a general framework that encompasses other families of robust risk measures. Similar to the previous sections, we establish the sufficient and necessary conditions for finiteness and derive the finite-dimensional dual problems of the corresponding robust risk measures.  

We first note that the robust expected utility risk measure can also be formulated in terms of two $\phi$-divergences $(\phi_1,\phi_2)$. Indeed, as shown in \eqref{EU_dual_representation}, expected utility can be represented by the divergence $\phi_1=-u_*$. The dual reformulation of robust expected utility risk measures has already been derived in \citet{BreuerCsiszar}, and its implication for divergence choices has been studied in \cite{KCS19,KCS21}. Therefore, in this section, we focus on the class of robust shortfall risk measures, which can be characterized as follows:
\begin{align}\label{robust_uMean_b}
\begin{split}
    \rho^s_{\mathrm{sf},\mathbb{P}}(X)&\triangleq \sup_{\lambda >0}\sup_{\{\mathbb{Q}: I_{\phi_2}(\mathbb{Q}, \mathbb{P}_0)\leq r\}}\sup_{\{\bar{\mathbb{Q}}:\bar{\mathbb{Q}}\ll \mathbb{Q}\}}\mathbb{E}_{\bar{\mathbb{Q}}}[-X]\\
    &\qquad-I_{\phi^\lambda_1}(\bar{\mathbb{Q}},\mathbb{Q}).
\end{split}
\end{align}
%and 
%\begin{align}\label{robust_uMean_pen_bad}
%\begin{split}
%     \rho^l_{\mathrm{sf},\mathbb{P}}(X)&\triangleq\sup_{\lambda >0}\sup_{\{\mathbb{Q}: \mathbb{Q}\ll \mathbb{P}_0\}}\sup_{\{\bar{\mathbb{Q}}:\bar{\mathbb{Q}}\ll \mathbb{Q}\}}\mathbb{E}_{\bar{\mathbb{Q}}}[-X]\\
 %    &\qquad-I_{\phi^\lambda_1}(\bar{\mathbb{Q}},\mathbb{Q})- I_{\phi_2}(\mathbb{Q},\mathbb{P}_0).
%\end{split}
%\end{align}
By examining the dual of \eqref{robust_uMean_b}, we can once again establish the finiteness conditions for \eqref{robust_uMean_b} in terms of the tail behavior of $\phi_1^*$, $\phi_2^*$ and $f_0$. 
\begin{theorem}\label{thm: uMean_b}
        Let $X\in L^1(\mathbb{P}_0)$. Let $\phi_1,\phi_2\in \Phi_0$ be lower-semicontinuous. Then, we have that \eqref{robust_uMean_b} is finite, if and only if there exists $\theta_1,\theta_2\in \mathbb{R}, \eta\geq 0$, such that
\begin{equation}\label{finite_uMeanb_condition}
    \int_{\mathbb{R}} \eta \phi_2^*\left(\frac{\phi_1^*(-\theta_2-x)+\theta_1}{\eta}\right)f_0(x)\mathrm{d}x\leq \theta_1- \eta r.
\end{equation}
Moreover, we have the duality relation:
\begin{align}\label{dual_sf_risk_measure}
\begin{split}
       \rho^s_{\mathrm{sf},\mathbb{P}}(X)=&\inf_{\substack{\eta\geq  0\\\theta_1,\theta_2\in \mathbb
    R}}\left\{\theta_2~\middle |~\mathbb{E}_{\mathbb{P}_0}\left[ \beta^*(X,\theta_1,\theta_2,\eta)\right]\right.\\
    &\qquad \left.\leq\theta_1- \eta r\right\},
\end{split}
\end{align}
where $\beta^*(X,\theta_1,\theta_2,\eta)=\eta \phi_2^*\left(\frac{\phi_1^*(-\theta_2-X)+\theta_1}{\eta}\right)$.
\end{theorem} 
However, a similar duality result cannot be obtained if the robust shortfall risk measure is defined by adding a second $\phi_2$ penalization (such as in \eqref{intro_rob_measures_pen}), since the non-convex product term $\lambda \mathbb{Q}$ in $I_{\phi^\lambda_1}(\bar{\mathbb{Q}},\mathbb{Q})$ is hard to handle \endnote{In fact, this is also the case for \eqref{robust_uMean_b}, but the set formulation in \eqref{robust_uMean_b} allows us to circumvent the non-convexity by a change of variable}. Nevertheless, as an alternative, we note that one can instead also define a robust shortfall risk measure by adding the $\phi_2$ penalization to the expectation constraint in \eqref{u-Mean}, and then apply the reformulation in \cite{BreuerCsiszar} for the robust expectation. This would lead to the dual formulation in \eqref{rob_sf_pen}. 

\section{Elicitation of Divergences from Robust Risk Measures} \label{sec: OCE_properties}
In the previous sections, we have shown that the robust risk measures characterized by $(\phi_1,\phi_2)$ can be reformulated as finite-dimensional dual problems in terms of the conjugates $(\phi_1^*,\phi_2^*)$. In this section, we study the inverse problem of recovering the $(\phi_1^*,\phi_2^*)$ from the robust risk measures. In the non-robust case where ambiguity is not present, it is shown by \citet{BB91AOR, BT07} that utility functions which are normalized and strongly risk-averse (i.e., $u(0)=0, u(x)<x, \forall x\in \mathbb{R}$) can be elicited by querying the certainty equivalent values (EU, OCE, and u-Mean) of random variables $X_p$, which for any $x\in \mathbb{R}$ is defined as
\begin{equation}\label{Xp}
    X_p=\begin{cases}
        x&~\text{with probability}~p\\
        0&~\text{with probability}~1-p.
    \end{cases}
\end{equation}
More precisely, for any certainty equivalent $\mathrm{CE}\in \{\mathrm{EU},\mathrm{OCE},\text{u-Mean}\}$, we have the limit relation
\begin{equation}\label{CE_recover}
    \lim_{p\downarrow 0}\frac{\mathrm{CE}(X_p)}{p}=u(x).
\end{equation}
Our goal is to extend the identity \eqref{CE_recover} to robust risk measures $\{\rho^l_{\mathrm{eu},\mathbb{P}}, \rho^l_{\mathrm{oce},\mathbb{P}}, \rho^l_{\mathrm{sf},\mathbb{P}}\}$ that are defined through penalization.
In the case of robust expected utility $\rho^l_{\mathrm{eu},\mathbb{P}}$, such an extension is already evident due to its duality relation to OCE (see \eqref{def:OCE}):
\begin{equation*}
    \rho^l_{\mathrm{eu},\mathbb{P}}(X):=\inf_{\theta\in \mathbb{R}}-\theta +\mathbb{E}_{\mathbb{P}}[\phi_2^*\left(\phi_1^*(-X)+\theta\right)],
\end{equation*}
where the right-hand side is equal to $-\mathrm{OCE}_{v,\mathbb{P}}(u(X))$, for $u(x):=-\phi_1^*(-x)$ and $v(x)=-\phi_2^*(-x)$. Since identity \eqref{CE_recover} holds for the OCE, it follows that if $\phi_2^*(x)<x, \forall x\in \mathbb{R}\setminus \{0\}$, then we have that  
\begin{equation}\label{elicitate_EU}
    \lim_{p\downarrow 0}\frac{\rho^l_{\mathrm{eu},\mathbb{P}}(X_p)}{p}=\phi_2^*(\phi_1^*(-x)).
\end{equation}
In Proposition \ref{prop:OCEMeasure}, we state that the same relation \eqref{elicitate_EU} can also be derived for the robust OCE risk measure $\rho^l_{\mathrm{oce},\mathbb{P}}$ and the robust shortfall risk measure $\rho^l_{\mathrm{sf},\mathbb{P}}$ through their dual formulation. Recall that the dual of $\rho^l_{\mathrm{oce},\mathbb{P}}$ is given in \eqref{dual_OCE_p}. The robust risk measure $\rho^l_{\mathrm{sf},\mathbb{P}}$, which is defined by replacing the expected utility constraint in \eqref{u-Mean} with its robust counterpart, has the following dual formulation:
\begin{align}\label{rob_sf_pen}
\begin{split}
    \rho^l_{\mathrm{sf},\mathbb{P}}(X)=&\inf_{\theta_1,\theta_2\in \mathbb{R}}\left\{\theta_2~\middle |~\mathbb{E}_{\mathbb{P}}[\phi_2^*(\theta_1-u(X+\theta_2))]\right.\\
    &\qquad \left.\leq \theta_1\right\},
\end{split}
\end{align}
 
\begin{proposition}\label{prop:OCEMeasure}
    Let $\phi_1,\phi_2\in \Phi_0$ be lower-semicontinuous and suppose that $\mathrm{dom}(\phi^*_1)=\mathrm{dom}(\phi_2^*)=\mathbb{R}$. Assume furthermore that $\phi_1^*(s)>s$, $\phi_2^*(s)>s$ for all $s\neq 0$. Then, for $X_p$ defined in \eqref{Xp} with any $x\in \mathbb{R}$, we have
    \begin{equation}\label{oce_elicit}
        \lim_{p\downarrow 0}\frac{\rho^l_{\mathrm{oce},\mathbb{P}}(X_p)}{p}=\phi_2^*(\phi_1^*(-x)).
    \end{equation}
    If $\phi_2^*\circ \phi^*_1$ is also differentiable at zero with derivative $(\phi_2^*\circ \phi_1^*)'(0)=1$, then we also have 
    \begin{equation}\label{uMean_elicit}
        \lim_{p\downarrow 0}\frac{\rho^l_{\mathrm{sf},\mathbb{P}}(X_p)}{p}=\phi_2^*(\phi_1^*(-x)).
    \end{equation}
\end{proposition}

Therefore, Proposition \ref{prop:OCEMeasure} suggests a method for measuring the divergences $\phi^*_1,\phi^*_2$, by conducting the following experiment. One can first measure $\phi_1^*$, by retrieving the certainty equivalent for a sequence of $X_p$, where $p$ is known and tends to zero. In this case, the decision maker is not ambiguous and \eqref{CE_recover} gives a measurement of $u(x)=-\phi_1^*(-x)$. In a second experiment, we determine the robust certainty equivalent for $X_p$, where $p$ is now hidden from the decision maker. This will yield $\phi_2^*(\phi_1^*(-x))$. Of course, this only allows us to recover the values of $\phi_2^*$ on the image set of $\phi_1^*$. On the other hand, the values of $\phi_2^*$ outside this image set do not contribute to the values of robust risk measures and can be thus taken arbitrarily.

Finally, we note that besides the elicitation of divergences, the formulation of the robust OCE and robust shortfall risk measures in terms of $(\phi_1^*,\phi_2^*)$ also allows us to characterize risk aversion by simple convexity conditions imposed on $\phi_1^*,\phi_2^*$. We refer the readers to the Electronic Companion for more details.

\section{Computing the Robust OCE Risk Measures}\label{sec:Compute}
In this section, we discuss the computation of a robust OCE risk measure using sample average approximation (SAA).
Theorem \ref{thm:ref_rob_OCE_pen} states that robust OCE risk measures can be reformulated into the finite-dimensional dual problems \eqref{dual_OCE_p} and \eqref{dual_OCE_b}, which are stochastic programming problems. The SAA method approximates these stochastic optimization problems by replacing $\mathbb{P}_0$ with the empirical distribution constructed from samples $X_1,\ldots, X_N$ from $\mathbb{P}_0$. This leads to the following optimization problems:
\begin{equation}\label{SAA: dual_oce_p1}
\min_{\substack{\theta_1,\theta_2\in \mathbb
    R}}-\theta_1-\theta_2 +\frac{1}{N}\sum^N_{i=1} \phi_2^*\left(\phi_1^*(\theta_2-X_i)+\theta_1\right),
\end{equation}
and
\begin{align}\label{SAA: dual_oce_b1}
\begin{split}
    \min_{\substack{\lambda\geq 0\\\theta_1,\theta_2\in \mathbb
    R}}&-\theta_1-\theta_2+\lambda r \\
    &\qquad +\frac{1}{N}\sum^N_{i=1}\lambda \phi_2^*\left(\frac{\phi_1^*(\theta_2-X_i)+\theta_1}{\lambda}\right).
\end{split}
\end{align}
The formulations in \eqref{SAA: dual_oce_p1} and \eqref{SAA: dual_oce_b1} are convex optimization problems with only two or three variables. Therefore, they can be efficiently solved using the ellipsoid method (see, e.g., \citealt{Ellipsoid}). This is particularly useful if $\phi_1^*$ is only subdifferentiable, such as the case of the CVaR function where $\phi_1^*(x)=\max\{x/\alpha ,0\}$. %Although the ellipsoid method requires a starting ellipsoid that contains the minimizer, which is hard to determine. In practice, we can always circumvent this by increasing the radius of the starting ellipsoid, and check a posteriori the optimality of the solution by examining its subgradient. 
If both $\phi_1^*$ and $\phi_2^*$ are conic representable, then \eqref{SAA: dual_oce_p1} and \eqref{SAA: dual_oce_b1} can also be solved using conic optimization software. Otherwise, one can also consider using solvers such as CVXOPT and Ipopt. In practice, we have observed that the ellipsoid method is the most promising one for the problems \eqref{SAA: dual_oce_p1} and \eqref{SAA: dual_oce_b1}. 

The consistency of the SAA estimation method, both in objective value and solutions, is established in Theorem 5.4 of \citet{LectureSP}, which requires the assumption that the optimal solutions of the true stochastic problems \eqref{dual_OCE_p} and \eqref{dual_OCE_b} must be contained in a compact set, for any nominal distribution $\mathbb{P}$. In the following proposition, we show that this is indeed the case.
\begin{proposition}\label{prop: compactOCE}
   Let $X\in L^1(\mathbb{P}_0)$, and $\phi_1,\phi_2\in \Phi_0$ be lower semi-continuous. Suppose that condition \eqref{ocep_integral} holds. Then, the decision variables $\theta_1,\theta_2$ in the dual problems \eqref{dual_OCE_p} and \eqref{dual_OCE_b} can be restricted to a compact set, without changing the optimum.
\end{proposition}

Together with Proposition \ref{prop: compactOCE} and Theorem 5.4 of \citet{LectureSP}, we can now state the consistency of the SAA method for robust OCE risk measures.

\begin{theorem}\label{thm: ConvergenceSAA}
    Let $X\in L^1(\mathbb{P}_0)$ and $\phi_1^*,\phi_2^*$ be proper, non-decreasing convex functions that satisfy the following conditions
    \begin{enumerate}
        \item $\mathrm{dom}(\phi_1^*)=\mathrm{dom}(\phi_2^*)=\mathbb{R}$
        \item $\phi_2^*(s)>s$, $\phi_1^*(s)>s$, for all $s\neq 0$ and $\phi_1^*(0)=\phi_2^*(0)=0$.
    \end{enumerate}
    Assume furthermore that there exists a $(\theta_1^0,\theta_2^0,\lambda_0)$ such that for all $(\theta_1,\theta_2,\lambda)$ in a neighborhood of $(\theta_1^0,\theta_2^0,\lambda_0)$, we have
    \begin{equation*}
    \mathbb{E}_{\mathbb{P}_0}\left[ \lambda\phi_2^*\left(\frac{\phi_1^*(\theta_2-X)+\theta_1}{\lambda}\right)\right]<\infty.    
    \end{equation*}
    Then, as the sample size $N\to\infty$, the SAA approximations \eqref{SAA: dual_oce_p1} and \eqref{SAA: dual_oce_b1} converge both in optimal objective value and in optimal solutions respectively to those of the true problems \eqref{dual_OCE_p} and \eqref{dual_OCE_b}, with probability 1.
\end{theorem}

\subsection{Importance Sampling}
Although theoretically the SAA method converges, it might be very slow in practice due to a lack of data in the tails of the nominal distribution $\mathbb{P}_0$ of $X$. As noted by \citet{Opt_Curse}, empirical risk minimization often underestimates the true risk. If we have a closed-form of the nominal density function $f_0$, then importance sampling can also be employed to improve the performance of SAA through a change of measure. For example, let $g$ be the density function of a measure $\tilde{\mathbb{P}}_0$ such that $\mathbb{P}_0\ll \tilde{\mathbb{P}}_0$, and let $Y_1,\ldots, Y_N$ be i.i.d. samples of $\tilde{\mathbb{P}}_0$. Then, through a multiplication of a likelihood factor, we can also express the SAA of \eqref{SAA: dual_oce_b1} as
\begin{align}\label{ImpSAA: dual_oce_b1}
\begin{split}
    &\min_{\substack{\lambda\geq 0\\\theta_1,\theta_2\in \mathbb
    R}}-\theta_1-\theta_2+\lambda r \\
    &\qquad+\frac{1}{N}\sum^N_{i=1}\lambda \phi_2^*\left(\frac{\phi_1^*(\theta_2-Y_i)+\theta_1}{\lambda}\right)\frac{f(Y_i)}{g(Y_i)}.
\end{split}
\end{align}
Note that the likelihood factor $f(Y_i)/g(Y_i)$ does not depend on the decision variables, and thus does not affect the tractability of \eqref{ImpSAA: dual_oce_b1}.

\subsection{Addressing Sampling Errors and Model Misspecification Using Globalized Robust Optimization}
Another way to address the sampling errors of SAA is to add an extra layer of robustification.  Therefore, we formulate the following \emph{globalized} robust OCE risk measure, by introducing a third divergence $\phi_3$:   
\begin{align}\label{GRO_risk_measure_discrete}
\begin{split}
    \rho^g_{\mathrm{oce},\mathbb{P}_0}(X)&= \sup_{\{\mathbb{P}:I_{\phi_3}(\mathbb{P},\mathbb{P}_0)\leq r\}}\sup_{\{\mathbb{Q}:\mathbb{Q}\ll \mathbb{P}\}}\sup_{\substack{\{\bar{\mathbb{Q}}:\bar{\mathbb{Q}}\ll \mathbb{Q}\}}}\\
    &\qquad\mathbb{E}_{\bar{\mathbb{Q}}}[-X]-I_{\phi_1}(\bar{\mathbb{Q}},\mathbb{Q})-I_{\phi_2}(\mathbb{Q},\mathbb{P}).
\end{split} 
\end{align}
We call \eqref{GRO_risk_measure_discrete} the globalized robust OCE risk measure since it utilizes the idea of globalized robust optimization introduced by \citet{BenTalGlobal}. We note that \eqref{GRO_risk_measure_discrete} is a robust version of \eqref{intro_rob_measures}, where extra robustness is added to address the sampling error of SAA approximation for computing a robust risk measure. In the following theorem, we show that \eqref{GRO_risk_measure_discrete} can also be reformulated as a finite-dimensional dual optimization problem. 
\begin{theorem}\label{thm:GRO_ref_discrete}
    Let $X\in L^1(\mathbb{P}_0)$ and $\phi_1,\phi_2, \phi_3\in \Phi_0$ be lower-semicontinuous. Then, we have the duality:
\begin{align}\label{dual_GRO_discrete}
\begin{split}
    \rho^g_{\mathrm{oce},\mathbb{P}_0}(X) 
   &=\inf_{\substack{\lambda\geq 0\\\theta_1,\theta_2,\theta_3\in \mathbb
    R}}-\theta_1-\theta_2-\theta_3\\
    &\qquad+\lambda r+\mathbb{E}_{\mathbb{P}_0}\left[\gamma^*(X,\theta_1,\theta_2,\theta_3,\lambda)\right],  
\end{split}
\end{align}
where $\gamma^*(X,\theta_1,\theta_2,\theta_3,\lambda)=\lambda \phi_3^*\left(\frac{\phi_2^*(\phi_1^*(\theta_3-X)+\theta_2)+\theta_1}{\lambda}\right)$.
\end{theorem}
We note that the globalized robust OCE risk measure also offers a new alternative for constructing uncertainty sets to address tail misspecification. Indeed, the duality relation \eqref{dual_GRO_discrete} also leads to the following integral condition that guarantees the finiteness of the globalized robust OCE risk measure under the nominal $\mathbb{P}_0$, namely that there must exist $\theta_1,\theta_2,\theta_3\in \mathbb{R}$ and $\lambda\geq 0$ such that
    \begin{equation}\label{GROintegral_condition}
   \int_{\mathbb{R}} \lambda\phi_3^*\left(\frac{\phi_2^*(\phi_1^*(\theta_3-x)+\theta_2)+\theta_1}{\lambda}\right)f_0(x)\mathrm{d}x<\infty.
    \end{equation}
The presence of an extra divergence $\phi_3^*$ in condition \eqref{GROintegral_condition}, compared to the non-globalized case \eqref{oceb_integral}, offers certain flexibility advantages. Instead of defining a single, more complicated divergence to tailor the tails of the conjugate functions and the nominal density, we can also use a composition of two elementary divergences $\phi_3^*\circ \phi_2^*$. For example, a composition of the KL-divergence $\phi_3^*(x)= \mathcal{O}(e^x)$ and the polynomial divergence with $\phi_2^*(x)=  \mathcal{O}(x^{\tilde{p}})$ for some $\tilde{p}>1$, would result in a tail behavior $\phi_3^*(\phi_2^*(-x))= \mathcal{O}(e^{|x|^{\tilde{p}}})$ as $x\to -\infty$, which could be a suitable divergence for $\mathrm{CVaR}_{\alpha}$ risk measures and a Weibull nominal distribution with shape parameter $k>1$. The advantage of using a combination of KL-divergence $\phi_3$ and a polynomial divergence $\phi_2$ is that both divergences have an explicit form of $\phi$ and $\phi^*$ functions.
\section{Numerical Examples}\label{sec: numerical examples}
\subsection{Toy Example}
We investigate a toy example where a risky claim $X$ follows a Pareto distribution on $[-\infty,-1)$ with density $f(x)=2/|x|^3$ for $x\leq -1$. We aim to measure its $\mathrm{CVaR}_{0.975}$ risk at $0.975$ level, aligning with the Basel III regulation. Theoretically, the exact $\mathrm{CVaR}_{0.975}(X)$ value can be calculated to be 12.649 (see \citealt{CVaR_Calc}). 

We investigate the robust $\mathrm{CVaR}_{0.975}$, which we define as in \eqref{intro_rob_measures} by means of a $\phi$-divergence uncertainty set  with some radius $r>0$. The sample average approximation of the robust $\mathrm{CVaR}_{0.975}$ is given by:
\begin{align}\label{robust_SAA}
\begin{split}
        &\min_{\substack{\lambda\geq 0\\\theta_1,\theta_2\in \mathbb
    R}}-\theta_1-\theta_2+\lambda r \\
    &\qquad+\frac{1}{N}\sum^{N}_{i=1}\lambda \phi_2^*\left(\frac{\max\left\{\frac{1}{1-\alpha}(\theta_2-X_i),0\right\}+\theta_1}{\lambda}\right).
\end{split}
\end{align}
Since $\mathbb{P}_0$ is the Pareto distribution, the polynomial divergence with $p=3$ would ensure a finite robust risk evaluation, according to Corollary \ref{cor:finite_oce}. We then solve the problem \eqref{robust_SAA} for a range of the radius $r$, with $N=1000$, and examine the optimal objective values. We note that \eqref{robust_SAA} is a conic optimization problem and thus can be solved efficiently. The results are given in Table \ref{tab:poly_cvar_empirical}. We see that when $r=0$, the SAA yields the value 9.792, which severely underestimates the true risk 12.649. However, the robust model is able to provide an upper bound on the true risk, already at a very small radius $r\geq 0.007$.

\begin{table}[H]\label{tab:poly_cvar_empirical}
\centering
\begin{tabular}{cc}
\hline
Radius & Robust $\mathrm{CVaR}_{0.975}$ \\ \hline
0&9.792\\ 
0.001  & 11.225                       \\ 
0.003   & 12.101                        \\ 
0.005   & 12.622                       \\ 
0.007   & 12.996                        \\ 
0.01   & 13.425                        \\ 
0.03  & 15.035                        \\ 
0.1      &  17.084                       \\ 
0.537    &   19.907                      \\ 
1& 20.940 \\ 
\end{tabular}
\caption{Robust $\mathrm{CVaR}_{0.975}$ risk measure with the polynomial divergence of degree $p=3$, calculated for 1000 samples of a Pareto loss $X$, for each radius $r$. }
\end{table}
We also illustrate the effect of choosing a ``wrong" divergence, by comparing the difference in optimal objective value of \eqref{robust_SAA} between the choice of a polynomial divergence of $p=3$, and the KL-divergence. As shown by condition \eqref{oceb_integral} in Corollary \ref{cor:finite_oce}, if the nominal distribution is Pareto $f(x)=2/{|x|^3}$ and the risk measure is $\mathrm{CVaR}_{\alpha}$, then the choice of KL-divergence will lead to an infinite conservative robust evaluation. To illustrate its effect on SAA, we generate samples from the Pareto distribution of different sizes and compute the optimal objective values of \eqref{robust_SAA} for both divergences. We fix radius at $r=0.02$. The results are shown in Table \ref{tab: poly_KL_compare}. As we can see, the robust $\mathrm{CVaR}_{0.975}$ value given by the KL-divergence is significantly larger and thus more conservative than that of the polynomial divergence. To make this difference even more visible, we use importance sampling to draw more data from the tail by sampling under the Pareto distribution $g(x)=1/x^2$. The results are given in Table \ref{tab: poly_KL_compare_Imp}, where the robust $\mathrm{CVaR}_{0.975}$ calculated using the KL-divergence now gives unreasonably large values with high variance, as reasoned from Corollary \ref{cor:finite_oce}. On the other hand, the robust $\mathrm{CVaR}_{0.975}$ value using the polynomial divergence exhibits more stability, with values that are not much different from the ones in Table \ref{tab: poly_KL_compare}. This shows the importance of choosing a proper divergence, when using robust optimization for risk assessment.

\begin{table}[H]
\centering
\begin{tabular}{ccc}
\hline
Sample Size & Robust $\mathrm{CVaR}_{0.975}$ Polynomial & Robust $\mathrm{CVaR}_{0.975}$ KL \\
\hline
500 & 8.667 & 9.584 \\
1000 & 17.433 & 25.893 \\
1500 & 17.545 & 24.439 \\
2000 & 15.314 & 18.947 \\
2500 & 26.614 & 41.370 \\
3000 & 24.072 & 44.025 \\
3500 & 21.727 & 33.579 \\
4000 & 17.038 & 22.733 \\
4500 & 24.152 & 35.767 \\
5000 & 22.413 & 40.402 \\
5500 & 21.493 & 32.968 \\
6000 & 19.433 & 31.693 \\
\end{tabular}
\caption{The robust $\mathrm{CVaR}_{0.975}$ value, equal to the optimal objective value of \eqref{robust_SAA} with radius $r=0.02$, is computed for a range of sample sizes drawn from a Pareto nominal distribution with $f(x)=2/|x|^3, x\leq -1$. The differences are compared for $\phi_2$ in \eqref{robust_SAA} chosen to be the polynomial divergence with $p=3$ and the KL-divergence with $\phi_2^*(x)=e^x-1$. } \label{tab: poly_KL_compare}
\end{table}

\begin{table}[H]
\centering
\begin{tabular}{ccc}
\hline
Sample Size & Robust $\mathrm{CVaR}_{0.975}$ Polynomial & Robust $\mathrm{CVaR}_{0.975}$ KL \\
\hline
500 & 23.758 & 142.732 \\
1000 & 28.121 & 107.253 \\
1500 & 26.280 & 544.011 \\
2000 & 26.962 & 250.208 \\
2500 & 24.948 & 487.849 \\
3000 & 24.906 & 219.724 \\
3500 & 27.996 & 457.070 \\
4000 & 29.142 & 2564.104 \\
4500 & 28.156 & 1006.002 \\
5000 & 33.983 & 71322.726 \\
5500 & 27.492 & 780.026 \\
6000 & 27.606 & 1262.779 \\
\end{tabular}
\caption{Robust $\mathrm{CVaR}_{0.975}$ with radius $r=0.02$, computed for samples drawn from a Pareto nominal distribution with $f(x)=1/x^2, x\leq -1$, using importance sampling. The differences are compared for $\phi_2$ being the polynomial divergence with $p=3$ and the KL-divergence with $\phi_2^*(x)=e^x-1$.}\label{tab: poly_KL_compare_Imp}
\end{table}
\subsection{Discrete Delta Hedging of Black-Scholes Model}
We consider a discrete delta hedging strategy for a call option under the Black-Scholes framework. In this model, the stock price $S$ follows a geometric Brownian motion with drift $\mu_S$ and volatility $\sigma_S>0$. The call option $C(S_T)$ with a strike price $K$ and a maturity time $T$, has a final payoff given by $C(S_T)=\max\{S_T-K,0\}$.

The discrete delta hedging strategy is implemented by constructing a self-financing portfolio that aims to replicate the payoff of a call option. At each time step $t_i=iT/{n}$, the value of the portfolio is the sum of the stock value $\mathrm{Stock}(i)$ and the cash value $\mathrm{Cash}(i)$ that grows with an interest rate $r_f$. The portfolio is self-financing, which means that any rebalancing of the shares in one asset is only achieved by buying or selling from the other asset. Let $\delta_{S,t_i}$ denote the number of shares of stocks held at time $[t_{i-1},t_i)$, which we calculate using the Black-Scholes formula as in the continuous setting. The dynamics of the stock value and the cash value are determined by the following recursive relation:
\begin{align}\label{dynamics_port}
\begin{split}
  \mathrm{Stock}(t_i)&= S_{t_i}\cdot \delta_{S,t_i}\\
   \mathrm{Cash}(t_i)&=e^{r_f\cdot \frac{1}{n}}-S_{t_i}\cdot (\delta_{S,t_i}-\delta_{S,t_{i-1}})\\
   &\qquad -k_0-k\cdot | \delta_{S,t_i}-\delta_{S,t_{i-1}}|\cdot S_{t_i},  
\end{split} 
\end{align}
where we have also added the transaction cost, $k_0+k\cdot | \delta_{S,t_i}-\delta_{S,t_{i-1}}|\cdot S_{t_i}$, which consists of a fixed amount $k_0$ and a cost that is proportional to the volumes of the trade executed. We define the initial investment of $\mathrm{Stock}(0), \mathrm{Cash}(0)$ to be the portfolio that determines the price of the call option by the Black-Scholes formula, without the transaction cost. The hedging error, which we aim to measure, is defined as the absolute difference between the call option payoff and the payoff of this hedging strategy, evaluated at the maturity date $T$:
\begin{equation*}
    X_{\mathrm{error}}=\left|C(S_T)-\mathrm{Stock}(T)-\mathrm{Cash}(T)\right|.
\end{equation*}
In this experiment, we will measure the risk of the hedging error $X_{\mathrm{error}}$ both in the nominal (when $r=0$) and the robust setting (when $r>0$). Measurements of model risk of hedging error have been previously studied by \citet{KCS19} and \citet{GlassermanXU}. In contrast to these studies, where they examined the expectation of $X_{\mathrm{error}}$, we aim to measure the nominal and robust Conditional Value-at-Risk of $X_{\mathrm{error}}$ at confidence level $\alpha=0.95$. 

To choose a proper $\phi_2$-divergence for measuring ambiguity, we first note that $0<\delta_{S_{t_i}}<1$ since it is chosen according to the Black-Scholes formula. Hence, according to the dynamics given in \eqref{dynamics_port}, one can show that there exists a constant $C>0$, such that $X_{\mathrm{error}}$ satisfies \eqref{bound_risk_factor}:
\begin{equation*}
    |X_{\mathrm{error}}|\leq C\left(1+\sum^n_{i=1}S_{t_i}\right).
\end{equation*}
Therefore, we can choose the divergence according to the nominal distributions of the risk factors $S_{t_i}$, which, for each $i$, follow the log-normal distribution with $\sigma_i= i\sigma_S/\sqrt{n}$. Hence, we can set the volatility parameter of the divergence in \eqref{div_gl_cvar} to be $\sigma= \max_{i}\sigma_i=\sigma_S$. Furthermore, we choose $d=2$, so that our $\phi_2$-divergence set may contain distributions with tails as heavy as $|x|^{-2}$.

\begin{figure}[t]
    \centering
    \caption{Absolute hedging error measured by the nominal and robust Conditional Value-at-Risk at confidence level $\alpha=0.95$, as a function of the hedging frequency.}
    \label{fig:blackscholes}
    \includegraphics[width = 0.8\linewidth]{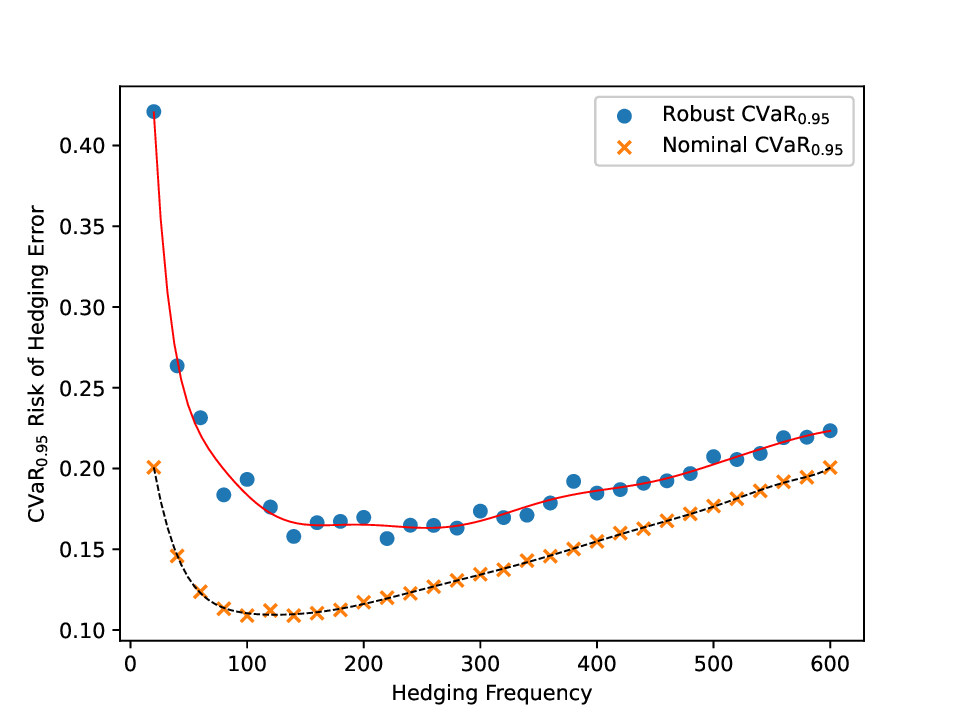}
\end{figure}

As a comparison, we choose the same parameters as in \citet{KCS19}, where $\mu_S= 0.05, \sigma_S=0.3, r_f=0.01, T=1, S_0=1, K=1, k_0 = 0.0002$ and $k=0.005$. We vary the hedging frequency $n$ and simulate for each frequency $8000$ hedging error data points. Then, we calculate the SAA approximation of the $\mathrm{CVaR}_{0.95}$ value of the hedging error, as well as that of the robust $\mathrm{CVaR}_{0.95}$ value, with the radius in \eqref{intro_rob_measures} set to be $r=0.1$. The results are given in Figure \ref{fig:blackscholes}. As we can see, both the nominal and the robust tail risk of the hedging error follow a trend of a $U$-shape as the hedging frequency increases. Indeed, even though theoretically the hedging error would decrease as the hedging frequency increases, in practice, the transaction cost would increase the hedging error. We see that the optimal hedging frequency in the nominal case is around $100$, slightly larger than the case investigated in \citet{KCS19}, where the expected hedging error was measured. In the robust case, the optimal frequency is increased to around $220$, which is substantially different from the nominal case. Moreover, the robust risk value is significantly larger than the nominal risk, showing that model uncertainty increases the risk premium.

\subsection{Risk-Averse Newsvendor Minimization}
We study the risk-averse newsvendor minimization problem, where we aim to minimize the loss in profit of selling a product under uncertain demand. The stock owner needs to decide a priori the amount of products that must be ordered, before the demand is realized. Let $d$ be a realized demand value of the uncertain demand variable $D$, $c$ be the cost of one unit of order, $v>c$ be the selling price, $s<c$ be the salvage value per unsold item returned to the factory, and $l$ be the loss per unit of unmet demand. If the stock owner decides to order $y$ number of items of the product, then the profit function is given by:
\begin{align*}
    \pi(y,d)&\triangleq v\min\{d,y\}+s\max\{y-d,0\}\\
    &\qquad -l\max\{d-y,0\}-cy.
\end{align*}
In the classical newsvendor problem, one minimizes the expected loss $\mathbb{E}[-\pi(y,D)]$ with respect to $y$. In the risk-averse setting, we aim to minimize the Conditional Value-at-Risk $\mathrm{CVaR}_{\alpha}(\pi(y,D))$ in $y$, where $\alpha$ controls the risk-aversion of the decision maker. Without ambiguity, the risk-averse newsvendor $\mathrm{CVaR}$ minimization problem has a closed-form solution (see \citealt{CVaRNews_European}), which is given by: 
\begin{align}\label{news_y_opt}
\begin{split}
    y^*_{\mathrm{nom}}&=\frac{E+V}{E+U}F^{-1}\left(\frac{U(1-\alpha)}{E+U}\right)\\
    &\qquad+\frac{U-V}{E+U}F^{-1}\left(\frac{E\alpha +U}{E+U}\right)
\end{split}
\end{align}
where $F^{-1}$ is the quantile function of the demand, and $E= c-s, U= v+l-c, V=v-c$.

\begin{figure}[htb]
    \centering
    \caption{The effect of increasing radius $r$ of the ambiguity set (thus increasing model uncertainty) on the optimal ordering decision. }
    \label{fig:newsvendor}
    \includegraphics[width = 0.8\linewidth]{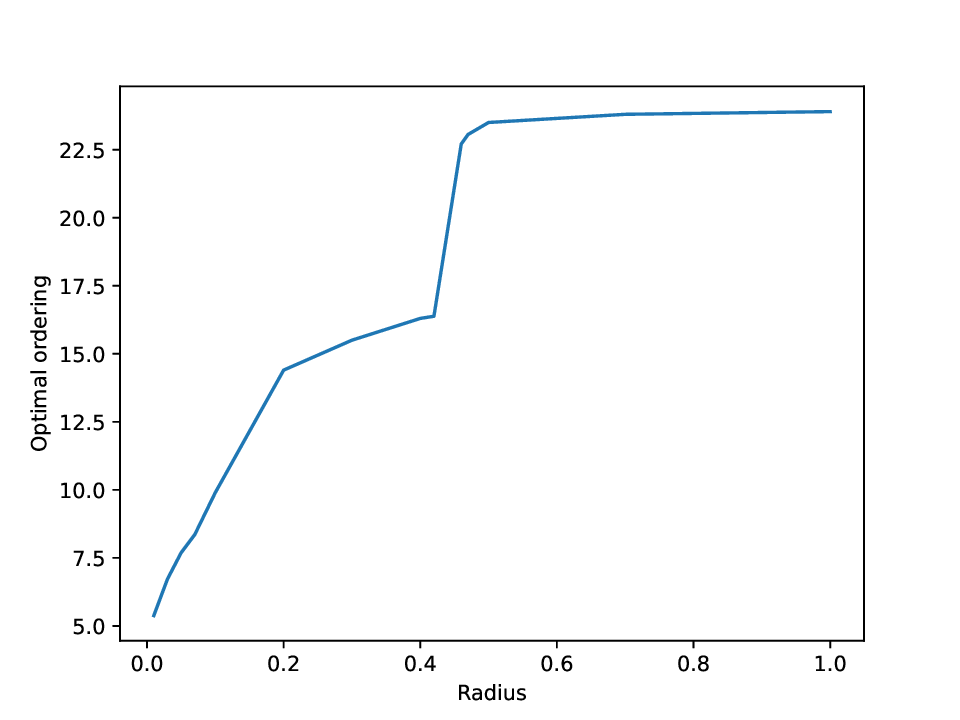}   
\end{figure}

We investigate the effect of model uncertainty on the optimal number of orderings. We do this by minimizing the robust $\mathrm{CVaR}_{0.95}(\pi(y,D))$ over $y$, for a range of radii $r>0$. Since robust $\mathrm{CVaR}_{0.95}(\pi(y,D))$ is a convex function of $y$ and $y$ is one-dimensional, we can search for a local minimum for a range of $y$ that is also a global optimum.
We draw $D_1,\ldots, D_{10000}$ demand samples from the heavy-tailed log-normal distribution with mean $\mu_D=0$ and standard deviation $\sigma_D=1$. Note that since $\pi(y,D)$ is a piecewise-linear concave function in $D$, we have that $|\pi(y,D)|$ also satisfies \eqref{bound_risk_factor} for some constant $C>0$. Hence, we may choose the divergence \eqref{div_gl_cvar} with $\sigma=\sigma_D$. For the remaining parameters, we choose $v = 8, c= 4, s= 2, l =4$. Then, one can calculate using \eqref{news_y_opt} that the optimal decision without model uncertainty is given by $y^*_{\mathrm{nom}}= 4.2$. As shown in Figure \ref{fig:newsvendor}, when model uncertainty increases, the optimal number of orderings also increases and becomes much larger than the optimal ordering $y^*_{\mathrm{nom}}$ without ambiguity. We also observe a jump that occurs around radius $0.45$, which could be due to the piecewise-linear property of the profit function. It is perhaps not very surprising that model uncertainty leads to an increase in the optimal ordering, since each excess ordered items incurs a loss of $s-c=2$ money unit, which is less than the $l=4$ money unit loss due to unmet demand. Therefore, in the presence of model uncertainty and heavy-tailed demand, the decision maker would increase the ordering to minimize the loss of unmet demand.   

\newpage
\section{Concluding Remarks}\label{sec: conclusion}
In this paper, we characterized a broad class of robust risk measures in terms of two $\phi$-divergences $(\phi_1,\phi_2)$, where $\phi_1$ specifies the risk measures and $\phi_2$ specifies the ambiguity set. We have shown that this framework is not only computationally tractable, but also offers a blueprint to systematically construct $\phi$-divergence uncertainty sets that are tailored to a given risk measure and nominal model. In addition, this characterization extends naturally to include higher level of uncertainty using globalized robust optimization. For many robust risk measures, we can also elicit the divergence functions $(\phi_1,\phi_2)$, which serve as an extra tool for calibrating the uncertainty set.
\newpage
\bibliographystyle{plainnat}
\bibliography{compositebiblio}
\newpage

\section*{Appendix}
\appendix
\renewcommand{\thesection}{EC.\arabic{section}}
\counterwithout{equation}{section}
\renewcommand{\theequation}{EC.\arabic{equation}}
\theoremstyle{plain}
\newtheorem{mylemma}{Lemma}[section]
\renewcommand{\themylemma}{EC.\arabic{section}.\arabic{mylemma}}
\theoremstyle{plain}
\newtheorem{mytheorem}{Theorem}[section]
\renewcommand{\themytheorem}{EC.\arabic{section}.\arabic{mytheorem}}
\renewcommand{\theproposition}{EC.\arabic{section}}
\section{Proofs}

\begin{proof}[\textbf{Proof of Theorem \ref{thm:ref_rob_OCE_pen}}]
Due to the absolute continuity relations $\bar{\mathbb{Q}}\ll\mathbb{Q}\ll \mathbb{P}_0$, there exists probability density functions $g,\bar{g}$ with respect to the measure $\mathbb{P}_0$, such that the optimization problem:
\begin{align}
\nonumber
    \sup_{\substack{\{\mathbb{Q}: \mathbb{Q}\ll \mathbb{P}_0\}\\\{\bar{\mathbb{Q}}:\bar{\mathbb{Q}}\ll \mathbb{Q}\}}}\mathbb{E}_{\bar{\mathbb{Q}}}[-X]-I_{\phi_1}(\bar{\mathbb{Q}},\mathbb{Q})-I_{\phi_2}(\mathbb{Q},\mathbb{P}_0),
\end{align}
can be formulated as an optimization problem over density functions:
\begin{align*}
    -\inf_{\substack{g,\Bar{g}\geq 0,\\ \int_{\Omega} g \mathrm{d}\mathbb{P}_0=1, \int_{\Omega} \bar{g}\mathrm{d}\mathbb{P}_0=1}}\int_{\Omega}\beta(\omega,g(\omega),\bar{g}(\omega))\mathrm{d}\mathbb{P}_0(\omega),
\end{align*}
where 
\begin{align*}
\beta(\omega,s,t)=X(\omega)t+s\phi_1\left(\frac{t}{s}\right)+\phi_2 \left(s\right), \qquad s,t\geq 0.
\end{align*}
Note that $\beta(\omega,s,t)=+\infty$ for $s<0$ or $t<0$ since both $\phi_1$ and $\phi_2$ take the value $+\infty$ on the negative real line. By assumption, $\phi_1$ and $\phi_2$ are convex lower-semicontinuous functions that have effective domain with non-empty interior. Therefore, $\beta(\omega,s,t)$ is a normal convex integrand (in the sense of Lemma 2  of \citealt{RockafellarIntegrals} and Proposition 14.39 of \citealt{Rockafellar&Wets}). Furthermore, there exists a feasible point, namely the constant function $(g,\bar{g})\equiv (1,1)$: $\int_{\Omega}\beta(\omega,1,1)\mathrm{d}\mathbb{P}_0(\omega)=\mathbb{E}_{\mathbb{P}_0}[X]+\phi_1(1)+\phi_2(1)<\infty$. Therefore, applying Theorem \ref{thm:J_beta} yields
\begin{align*}
    \inf_{\substack{g,\Bar{g}\geq 0,\\ \int_{\Omega} g \mathrm{d}\mathbb{P}_0=1, \int_{\Omega} \bar{g}\mathrm{d}\mathbb{P}_0=1}}\int_{\Omega}\beta(\omega,g(\omega),\bar{g}(\omega))\mathrm{d}\mathbb{P}_0(\omega)=\sup_{\theta_1,\theta_2\in \mathbb{R}}\theta_1+\theta_2-\int_\Omega\beta^*(\omega,\boldsymbol{\theta})\mathrm{d}\mathbb{P}_0(\omega),
\end{align*}
where
\begin{align*}
    \beta^*(\omega,\boldsymbol{\theta})&=\sup_{s,t\geq 0}\theta_1s+\theta_2t-X(\omega)t-s\phi_1\left(\frac{t}{s}\right)-\phi_2\left(s\right)\\
    &=\sup_{s\geq 0}\theta_1s-\phi_2\left(s\right)+\sup_{t\geq 0}(\theta_2-X(\omega))t-s\phi_1\left(\frac{t}{s}\right)\\
    &=\sup_{s\geq 0}\theta_1s-\phi_2\left(s\right)+s\phi_1^*(\theta_2-X(\omega))\\
    &=\phi_2^*(\theta_1+\phi_1^*(\theta_2-X(\omega))).
\end{align*}
Finally, to obtain the worst-case density, we have that $\mathbb{E}_{\mathbb{P}_0}\left[\phi_2^*(\phi_1^*(\theta_2-X)+\theta_1)\right]<\infty$ for some $(\theta_1,\theta_2)$. Hence, the dual problem \eqref{dual_OCE_p} is finite and there exists a dual solution $(\theta_1^*,\theta_2^*)$ by Theorem \ref{thm:Luenberger}. If furthermore that $(\theta_1^*,\theta_2^*)\in \mathrm{int}(\mathcal{F}(\boldsymbol{\theta}))$, then $\mathbb{E}_{\mathbb{P}_0}\left[\phi_2^*(\phi_1^*(\theta_2-X)+\theta_1)\right]$ is differentiable at $(\theta_1^*,\theta_2^*)$ by Theorem \ref{thm:diff_kbeta}. Therefore, the partial derivatives $(g^*(\omega),\bar{g}^*(\omega)):=(\partial_{\theta_1} \beta^*(\omega,\boldsymbol{\theta}), \partial_{\theta_2}\beta^*(\omega,\boldsymbol{\theta}))$
evaluated at $\boldsymbol{\theta}^*$, for each $\omega\in \Omega$,
satisfies the first-order condition of the dual problem \eqref{dual_OCE_p}
\begin{align}
\nonumber
    \int_\Omega g^*(\omega)\mathrm{d}\mathbb{P}_0(\omega)=1,\quad \int_\Omega \bar{g}^*(\omega)\mathrm{d}\mathbb{P}_0(\omega)=1.
\end{align}
Hence, it follows from Theorem \ref{thm:wc_density} that the worst-case densities are the partial derivatives.

Similarly, the optimization problem
\begin{align*}
\sup_{\substack{\{\mathbb{Q}:I_{\phi_2}(\mathbb{Q},\mathbb{P}_0)\leq r\}}}\sup_{\{\bar{\mathbb{Q}}:\bar{\mathbb{Q}}\ll \mathbb{Q}\}}\mathbb{E}_{\bar{\mathbb{Q}}}[-X]-I_{\phi_1}(\bar{\mathbb{Q}},\mathbb{Q})
\end{align*}
can be rewritten to the following optimization problem over density functions with respect to the measure $\mathbb{P}_0$:
\begin{align}\label{primal_oceb}
    -\inf_{\substack{g, \bar{g}\geq 0,\\ \int_{\Omega} g\mathrm{d}\mathbb{P}_0=1,\int_{\Omega} \bar{g}\mathrm{d}\mathbb{P}_0=1\\
\int_{\Omega}\phi_2\left(g(\omega)\right)\mathrm{d}\mathbb{P}_0(\omega)\leq r}}\int_{\Omega}X(\omega)\bar{g}(\omega)+g(\omega)\phi_1\left(\frac{\bar{g}(\omega)}{g(\omega)}\right)\mathrm{d}\mathbb{P}_0(\omega)=:-J_0.
\end{align}
Assume without loss of generality that $J_0>-\infty$. Then, we also have that $J_0<\infty$, since $(g,\Bar{g})\equiv (1,1)$ is a feasible solution. Then, Theorem \ref{thm:Luenberger} implies that
\begin{align*}
    J_0=\sup_{\lambda \geq 0}-\lambda r+\inf_{\substack{g, \bar{g}\geq 0,\\ \int_{\Omega} g\mathrm{d}\mathbb{P}_0=1,\int_{\Omega} \bar{g}\mathrm{d}\mathbb{P}_0=1}}\int_{\Omega}X(\omega)\bar{g}(\omega)+g(\omega)\phi_1\left(\frac{\bar{g}(\omega)}{g(\omega)}\right)+\lambda \phi_2(g(\omega))\mathrm{d}\mathbb{P}_0(\omega).
\end{align*}
We first examine the case when $\lambda =0$, for which the above would become
\begin{align*}
    J_0&=\inf_{\substack{g, \bar{g}\geq 0,\\ \int_{\Omega} g\mathrm{d}\mathbb{P}_0=1,\int_{\Omega} \bar{g}\mathrm{d}\mathbb{P}_0=1}}\int_{\Omega}X(\omega)\bar{g}(\omega)+g(\omega)\phi_1\left(\frac{\bar{g}(\omega)}{g(\omega)}\right)\mathrm{d}\mathbb{P}_0(\omega)\\
    &=\essinf (X).
\end{align*}
On the other hand, since $0\phi_2^*(s/0)=0$ for $s\leq 0$ and $+\infty$ otherwise, we have
\begin{align*}
    &\sup_{\theta_1,\theta_2\in \mathbb{R}}\theta_1+\theta_2-\mathbb{E}_{\mathbb{P}_0}\left[0 \phi_2^*\left(\frac{\phi_1^*(\theta_2-X)+\theta_1}{0}\right)\right]\\
    &=\sup_{\theta_1,\theta_2\in \mathbb{R}}\left\{\theta_1+\theta_2~\middle |~ \esssup \phi_1^*(\theta_2-X)+\theta_1\leq 0\right\}\\
    &=\sup_{\theta_2\in \mathbb{R}}\theta_2-\phi_1^*(\theta_2-\essinf(X))\\
    &=\essinf(X).
\end{align*}
Hence, the duality formula holds in the case of $\lambda=0$. Assume now $\lambda>0$. Define $\beta_{\lambda}(\omega,s,t):= X(\omega)t+s\phi_1\left(\frac{t}{s}\right)+\lambda \phi_2\left(s\right)$. Note that $\beta_{\lambda}$ is only finite when $s,t\geq 0$. Then, applying Theorem \ref{thm:Luenberger} yields
\begin{align*}
  J_0=\sup_{\substack{\lambda \geq 0\\ \theta_1,\theta_2\in \mathbb{R}}}-\lambda r+\theta_1+\theta_2-\int_\Omega \beta^*_{\lambda}(\omega,\theta_1,\theta_2)\mathrm{d}\mathbb{P}_0(\omega).
\end{align*}
It remains to compute $\beta^*_{\lambda}(\omega,\theta_1,\theta_2)$, which is equal to  $\beta^*_\lambda(\omega,\theta_1,\theta_2)=\lambda \phi_2^*\left(\frac{\theta_1+\phi_1^*(\theta_2-X(\omega))}{\lambda}\right)$. 

Now for the worst-case densities, we note that the domain set $\mathcal{F}(\boldsymbol{\theta},\lambda)$ is assumed to be non-empty, thus a dual solution $(\theta_1^*,\theta_2^*,\lambda^*)$ exists by Theorem \ref{thm:Luenberger}. If the dual solution also lies in the interior of $\mathcal{F}(\boldsymbol{\theta},\lambda)$, then we have differentiability by Theorem \ref{thm:diff_kbeta}, and thus the dual solution must satisfy the gradient conditions
\begin{align*}
    &\int_\Omega (\phi_2^*)'\left(\frac{\theta_1^*+\phi_1^*(\theta^*_2-X(\omega))}{\lambda^*}\right) \mathrm{d}\mathbb{P}_0(\omega)=1\\
    &\int_\Omega (\phi_2^*)'\left(\frac{\theta_1^*+\phi_1^*(\theta^*_2-X(\omega))}{\lambda^*}\right)(\phi_1^*)'(\theta_2^*-X(\omega)) \mathrm{d}\mathbb{P}_0(\omega)=1\\
    &\int_\Omega \frac{\theta_1^*+\phi_1^*(\theta^*_2-X(\omega))}{\lambda^*}(\phi_2^*)'\left(\frac{\theta_1^*+\phi_1^*(\theta^*_2-X(\omega))}{\lambda^*}\right)- \phi_2^*\left(\frac{\theta_1^*+\phi_1^*(\theta^*_2-X(\omega))}{\lambda^*}\right)\mathrm{d}\mathbb{P}_0(\omega)=r,
\end{align*}
where the latter condition is equivalent to $\int_\Omega \phi_2\left((\phi_2^*)'\left(\frac{\theta_1^*+\phi_1^*(\theta^*_2-X(\omega))}{\lambda^*}\right)\right) \mathrm{d}\mathbb{P}_0(\omega)=r$, due to the relation $\phi_2=\phi_2^{**}$. Therefore, the partial derivatives $(g^*(\omega),\bar{g}^*(\omega)):=(\partial_{\theta_1} \beta^*_{\lambda}(\omega,\boldsymbol{\theta}), \partial_{\theta_2} \beta^*_{\lambda}(\omega,\boldsymbol{\theta}))$ evaluated at $(\boldsymbol{\theta}^*,\lambda^*)$, are indeed feasible density functions of the primal problem \eqref{primal_oceb}.

Finally, the partial derivatives are indeed also the worst-case densities, since
\begin{align*}
J_0&=-\lambda^* r+\theta^*_1+\theta^*_2-\int_\Omega \beta^*_{\lambda^*}(\omega,\theta^*_1,\theta^*_2)\mathrm{d}\mathbb{P}_0(\omega)\\
&\stackrel{(*)}{=}-\lambda^*r+\int_\Omega X(\omega)\bar{g}^*(\omega)+g^*(\omega)\phi_1\left(\frac{\bar{g}^*(\omega)}{g^*(\omega)}\right)+\lambda^*\phi_2(g^*(\omega))\mathrm{d}\mathbb{P}_0(\omega)\\
&=\int_\Omega X(\omega)\bar{g}^*(\omega)+g^*(\omega)\phi_1\left(\frac{\bar{g}^*(\omega)}{g^*(\omega)}\right)\mathrm{d}\mathbb{P}_0(\omega),
\end{align*}
where $(*)$ follows from the worst-case densities of \eqref{intro_rob_measures_pen}.
\end{proof}

\begin{proof}[\textbf{Proof of Corollary \ref{cor:finite_oce}}]
Recall that from Theorem \ref{thm:ref_rob_OCE_pen} we have,
\begin{align}\label{roce_dual_proof}
    \sup_{\substack{\{\mathbb{Q}: \mathbb{Q}\ll \mathbb{P}_0\}\\}}\mathbb{E}_{\bar{\mathbb{Q}}}[-X]-I_{\phi_1}(\bar{\mathbb{Q}},\mathbb{Q})-I_{\phi_2}(\mathbb{Q},\mathbb{P}_0)= \inf_{\theta_1,\theta_2\in \mathbb{R}}-\theta_1-\theta_2+\mathbb{E}_{\mathbb{P}_0}\left[\phi_2^*(\phi_1^*(\theta_2-X)+\theta_1)\right].
\end{align}
Hence, if  $\mathbb{E}_{\mathbb{P}_0}\left[\phi_2^*(\phi_1^*(\theta_2-X)+\theta_1)\right]=\infty$ holds for all $\theta_1,\theta_2$, then the infimum on the right-hand side of \eqref{roce_dual_proof} is indeed always $\infty$. If $\mathbb{E}_{\mathbb{P}_0}\left[\phi_2^*(\phi_1^*(\theta_2-X)+\theta_1)\right]<\infty$ for some $\theta_1,\theta_2\in \mathbb{R}$, then this infimum is bounded away from $-\infty$. Similarly, \eqref{oceb_integral} follows from \eqref{dual_OCE_b}. 
\end{proof}

\begin{proof}[\textbf{Proof of Proposition \ref{prop: Risk factors}}]
Apply Corollary \ref{cor:finite_oce} and note that
    \begin{align*}
       \mathbb{E}_{\mathbb{P}_0}\left[ \phi_2^*\left(\theta_1+\phi_1^*(\theta_2-X)\right)\right] &\leq \mathbb{E}_{\mathbb{P}_0}\left[ \phi_2^*\left(\theta_1+\phi_1^*\left(\theta_2+C(1+\sum^m_{i=1}|Z_i|)\right)\right)\right]\\
       &= \mathbb{E}_{\mathbb{P}_0}\left[ \phi_2^*\left(\theta_1+\phi_1^*\left(\frac{1}{m}\sum^m_{i=1}\theta_2+C(1+m|Z_i|)\right)\right)\right]\\
       &\leq \mathbb{E}_{\mathbb{P}_0}\left[ \phi_2^*\left(\theta_1+\frac{1}{m}\sum^m_{i=1}\phi_1^*\left(\theta_2+C(1+m|Z_i|)\right)\right)\right]\\
       &=\mathbb{E}_{\mathbb{P}_0}\left[ \phi_2^*\left(\frac{1}{m}\sum^m_{i=1}\theta_1+\phi_1^*\left(\theta_2+C(1+m|Z_i|)\right)\right)\right]\\
       &\leq \frac{1}{m}\sum^m_{i=1}\mathbb{E}_{\mathbb{P}_0}\left[ \phi_2^*\left(\theta_1+\phi_1^*\left(\theta_2+C(1+m|Z_i|)\right)\right)\right],
    \end{align*}
    where we used Jensen's inequality and the monotonicity of the conjugate functions $\phi^*_1, \phi^*_2$.
\end{proof}

\begin{proof}[\textbf{Proof of Proposition \ref{prop: construct_phi_conj}}]
    By the assumptions on $\psi$, we have that $\phi^*(s)$ as defined in \eqref{def_phi^*} is an increasing, continuous convex function on $\mathbb{R}$ with $\phi^*(0)=0$ and $(\phi^*)'(0)=1$. By definition, we have that $\phi(t):=\sup_{s\in \mathbb{R}}st-\phi^*(s)$. Therefore, $\phi$ is convex. We also have that $\phi(t)=+\infty$ for $t<0$, since $\phi^*(s)<0$ for $s<0$ and thus $\phi(t)\geq \sup_{s\leq 0}st= \infty$, if $t<0$. Moreover, $\phi$ is non-negative on $[0,\infty)$ with $\phi(1)=0$. Indeed, $\phi^*(0)=0$ ensures that $\phi(t)\geq 0\cdot t-\phi^*(0)=0$, and $(\phi^*)'(0)=1$ implies that $\phi(1)=0$ (Theorem 23.5, \citealt{Rockafellar}). Therefore, $\phi\in \Phi_0$ and is thus a $\phi$-divergence function.
    
    If $\psi$ is strictly convex and differentiable, then $\phi^*$ is by construction strictly convex and differentiable on $\mathbb{R}$. Theorem 26.3 of \citet{Rockafellar} implies that its conjugate $\phi:=(\phi^*)^*$ is essentially smooth. Therefore, $\phi$ is differentiable on the interior of its 
 effective domain, denoted as $\mathrm{dom}(\phi)$. By Corollary 23.5.1 and Theorem 23.4 of \citet{Rockafellar}, we have that $\mathrm{dom}(\phi)$ contains the image of the derivative $(\phi^*)'$, which is the set $(0,\infty)$, since by assumption $\psi'$ is continuous and strictly increasing with $\lim_{s\to \infty}\psi'(s)=\infty$, and that $(\phi^*)'(s)=e^s$, for $s\leq 0$. Furthermore, $\phi(0)=1<\infty$, since $\phi(t)=t\log t-t+1$ on $t\leq 1$. Hence, $\mathrm{dom}(\phi)=[0,\infty)$. Moreover, Corollary 23.5.1 of \citet{Rockafellar} implies that $\phi'=((\phi^*)')^{-1}$.
\end{proof}

The proof of Proposition \ref{prop: phi_bound} relies on the following Lemma.

\begin{mylemma}\label{lem:inverse}
    Let $f,g$ be two strictly increasing functions. If for some $y_0\in \mathbb{R}$ that $f^{-1}(y)\leq g^{-1}(y)$ for all $y\geq y_0$, then $f(x)\geq g(x)$, for all $x\geq f^{-1}(y_0)$.
\end{mylemma}
\begin{proof}[\textbf{Proof of Lemma \ref{lem:inverse}}]
    We prove by contradiction. Suppose there is a $x_0\geq f^{-1}(y_0)$ such that $f(x_0)<g(x_0)$. Then, $f(x_0)\geq y_0$. Since $f,g$ are increasing, we have that $f^{-1},g^{-1}$ are also increasing. Therefore, 
    \begin{align*}
        x_0=f^{-1}(f(x_0))\leq g^{-1}(f(x_0))<g^{-1}(g(x_0))=x_0,
    \end{align*}
    which is a contradiction.
\end{proof}

\begin{proof}[\textbf{Proof of Proposition \ref{prop: phi_bound}}]
For notational convenience, we temporarily denote the nominal density $f:= f_0$. Since $\phi_2'$ is strictly increasing and non-negative for $t\geq 1$, we have that for $t\geq 1$, 
  \begin{align*}
    \phi_2(t)=\int^t_1\phi_2'(s)\mathrm{d}s\leq (t-1)\phi_2'(t)\leq t\phi_2'(t).
  \end{align*}
  Since $(\psi_1')^{-1}(y)\leq (\phi_2')^{-1}(y)$ for all $y\geq y_0$, Lemma \ref{lem:inverse} implies that for all $t\geq t_0$, where $t_0:=\max\{(\psi_1')^{-1}(y_0),1\}$, we have that $\psi_1'(t)\geq \phi_2'(t)\geq 0$. Therefore, we have
  \begin{align*}
      &\int^\infty_{-\infty}\phi_2\left(\frac{g(x)}{f(x)}\right)f(x)\mathrm{d}x\\
      &=\int_{x: \frac{g(x)}{f(x)}\leq 1}\phi_2\left(\frac{g(x)}{f(x)}\right)f(x)\mathrm{d}x+ \int_{x: \frac{g(x)}{f(x)}> 1}\phi_2\left(\frac{g(x)}{f(x)}\right)f(x)\mathrm{d}x\\
      &\leq \int_{x: \frac{g(x)}{f(x)}\leq 1}\phi_2\left(\frac{g(x)}{f(x)}\right)f(x)\mathrm{d}x+ \int_{x: \frac{g(x)}{f(x)}> 1}g(x)\phi_2'\left(\frac{g(x)}{f(x)}\right)\mathrm{d}x
  \end{align*}
  Since $\phi_2$ is bounded on $[0,1]$, we only need to bound the second integral. By the limit assumption \eqref{moment_include}, there exists a $u_0$, such that for all $x\in \mathbb{R}$ with $|x|\geq u_0$, we have $\psi_2'\left(\frac{1}{f(x)}\right)\leq M |x|^d$, for some constant $M>0$. We have that $f(x)\to 0$ as $|x|\to \infty$. Hence, there exists a $x_0\geq u_0$ such that for all $x\in \mathbb{R}$ with $|x|\geq x_0$, we have that $\frac{1}{f(x)}\geq t_0$. We then split the following integral by $x_0$:
  \begin{align*}
      &\int_{x: \frac{g(x)}{f(x)}> 1}g(x)\phi_2'\left(\frac{g(x)}{f(x)}\right)\mathrm{d}x\\
      &=\int_{x: \frac{g(x)}{f(x)}> 1\cap \{x: |x|\leq x_0\}}g(x)\phi_2'\left(\frac{g(x)}{f(x)}\right)\mathrm{d}x + \int_{x: \frac{g(x)}{f(x)}> 1\cap \{x: |x|\geq x_0\}}g(x)\phi_2'\left(\frac{g(x)}{f(x)}\right)\mathrm{d}x.
  \end{align*}
  The first term is bounded by the assumption that $\frac{g}{f}$ is bounded on any compact interval, and that $\phi_2'$ is bounded on any compact interval as well, due to its continuity. Therefore, we only need to examine the second term. Due to the non-negativity and monotonicity of $\phi_2',\psi_2'$ on $[1,\infty)$, we have, 
  \begin{align*}
      \int_{x: \frac{g(x)}{f(x)}> 1\cap \{x: |x|\geq x_0\}}g(x)\phi_2'\left(\frac{g(x)}{f(x)}\right)\mathrm{d}x&\leq \int_{x: |x|\geq x_0} g(x)\phi_2'\left(\frac{g(x)}{f(x)}\right)\mathrm{d}x\\
      &\leq \int_{x: |x|\geq x_0} g(x)\phi_2'\left(\frac{1}{f(x)}\right)\mathrm{d}x\\
      &\leq \int_{x: |x|\geq x_0} g(x)\psi_1'\left(\frac{1}{f(x)}\right)\mathrm{d}x\\
      &\leq \int_{x: |x|\geq x_0} g(x)|x|^d\mathrm{d}x<\infty.
  \end{align*}
  Conversely, if there exists a divergence $\psi_2$ such that $(\phi_2^*)'(y)\leq (\psi_2')^{-1}(y)$ for all $y\geq 0$, then we have that $\phi_2' (t)\geq \psi_2'(t)$ for all $t\geq 1$. Since $\phi_2(1)=\psi_2(1)=0$, this also means that 
  \begin{align*}
\phi_2(t)=\int^t_{1}\phi_2'(s)\mathrm{d}s\geq \int^t_{1}\psi_2'(s)\mathrm{d}s=\psi_2(t),~\forall t\geq 1. 
  \end{align*}
  Therefore, we have
  \begin{align*}
       \int_{x: \frac{g(x)}{f(x)}> 1}\phi_2\left(\frac{g(x)}{f(x)}\right)f(x)\mathrm{d}x\geq \int_{x: \frac{g(x)}{f(x)}> 1}\psi_2\left(\frac{g(x)}{f(x)}\right)f(x)\mathrm{d}x. 
  \end{align*}
Since $\psi_2$ is bounded on $[0,1]$, we have that $I_{\psi_2}(g,f)=\infty$ if and only if $\int_{x: \frac{g(x)}{f(x)}> 1}\psi_2\left(\frac{g(x)}{f(x)}\right)f(x)\mathrm{d}x=\infty$, which implies that $I_{\phi_2}(g,f)=\infty$.
\end{proof}

\begin{proof}[\textbf{Proof of Lemma \ref{lem: div_gl_cvar}}]
    Following the construction in \eqref{def_phi^*}, we choose $\psi(s)=(s+e)e^{a\log^p(s+e)}$, where $a = \frac{1}{p\left(\sigma d\right)^p}$. Clearly, $\psi(s)$ is increasing and twice continuously differentiable on $s\geq 0$. Moreover, it is strictly convex, since its second derivative is given by 
    \begin{align*}
        \psi''(s)=\frac{ap\log^{p-2}(s+e)e^{a\log^p(s+e)}(ap\log^p(s+e)+p+\log(s+e)-1)}{s+e},
    \end{align*}
    which is strictly positive for all $s\geq 0$. Hence, $\phi_2^*$ as defined in \eqref{div_gl_cvar} satisfies all properties given in Proposition \ref{prop: construct_phi_conj}. 

    We now show that the conjugate function $\phi_2^*$ satisfies the condition \eqref{ocep_integral}, which also implies \eqref{oceb_integral}. For any $\theta_1\geq 0,\theta_2\in \mathbb{R}$, we have that $\phi_1^*(\theta_2-x)=\max\left\{\frac{1}{\alpha}(\theta_2-x),0\right\}$. Since $f_0(x)$ has support $(-\infty,0]$, we have that $\frac{1}{\alpha}(\theta_2-x)\geq 0$, if $x\leq x_0$ for some $x_0$. Therefore, for all $x\leq x_0$, we have
    \begin{align*}
        &\phi_2^*\left(\phi_1^*(\theta_2-x)+\theta_1)\right)f_0(x)=\phi_2^*\left(\frac{1}{\alpha}(\theta_2-x)+\theta_1\right)\frac{1}{C\sigma |x|}\exp\left\{-\frac{1}{p\sigma^p}(\log(|x|)-\mu)^p\right\}.
    \end{align*}
    We analyze the asymptotic behavior of the above function. We have that
    \begin{align*}
        &\phi_2^*\left(\frac{1}{\alpha}(\theta_2-x)+\theta_1\right)\frac{1}{C\sigma |x|}\exp\left\{-\frac{1}{p\sigma^p}(\log|x|-\mu)^p\right\}\\
        &=\frac{c_1\left(\frac{1}{\alpha}(\theta_2+|x|)+\theta_1+e\right)}{C\sigma|x|}\cdot \exp \left(\frac{1}{p(\sigma d)^p}\log^p\left(\frac{1}{\alpha}(\theta_2+|x|)+\theta_1+e\right)-\frac{1}{p\sigma^p}(\log|x|-\mu)^p\right)\\
        &\qquad + \left(\frac{c_2\left(\frac{1}{\alpha}(\theta_2+|x|)+\theta_1\right)}{C\sigma|x|}+c_3\right)\cdot \exp \left(-\frac{1}{p\sigma^p}(\log|x|-\mu)^p\right)\\
        &=\mathcal{O}(1)\cdot \exp \left(\frac{1}{p(\sigma d)^p}\log^p\left(\frac{1}{\alpha}(\theta_2+|x|)+\theta_1+e\right)-\frac{1}{p\sigma^p}(\log|x|-\mu)^p\right)+\\
        &\qquad  \mathcal{O}(1)\cdot \exp \left(-\frac{1}{p\sigma^p}(\log(|x|)-\mu)^p\right).
    \end{align*}
We only need to examine that, 
\begin{align*}
    &\exp \left(\frac{1}{p(\sigma d)^p}\log^p\left(\frac{1}{\alpha}(\theta_2+|x|)+\theta_1+e\right)-\frac{1}{p\sigma^p}(\log|x|-\mu)^p\right)\\
    &=\exp \left(\frac{1}{p(\sigma d)^p}\log^p\left(|x|\left(\frac{1}{\alpha}+\frac{\frac{\theta_2}{\alpha}+\theta_1+e}{|x|}\right)\right)-\frac{1}{p\sigma^p}\log^p|x|(1+o(1))^p\right)\\
    &=\exp \left(\frac{1}{p(\sigma d)^p}\log^p|x|\left(1+o(1)\right)^p-\frac{1}{p\sigma^p}\log^p|x|(1+o(1))^p\right)\\
    &=\exp \left(-\frac{1}{p\sigma^p}\log^p|x|\left(\frac{1}{d^p}(1+o(1))^p-(1+o(1))^p\right)\right).
\end{align*}
Therefore, the above expression has the same asymptotic as $\exp \left(-\frac{1}{p\sigma^p}\log^p|x|\left(\frac{1}{d^p}-1\right)\right)$. Since $d>1, p\geq 2\Rightarrow 1/d^p-1<0$, we have that \eqref{ocep_integral} is satisfied. 

Finally, to show that the divergence ball of $\phi_2$ contains distributions that have finite $d$-th moment. First, we examine the derivative of $\phi_2^*$, which is also equal to $(\phi'_2)^{-1}$:
\begin{align*}
    (\phi_2^*)'(s)=\begin{cases}
        c_1\exp\left\{\frac{1}{p(\sigma d)^p}\log^p(s+e)\right\}\left(\frac{1}{(\sigma d)^p}\log^{p-1}(s+e)+1\right)+c_2 &s>0\\
        e^s& s\leq 0.
    \end{cases}
\end{align*}
We define a function
\begin{align*}
    (\psi_2')^{-1}(s)=\begin{cases}
        c_1\exp\left\{\frac{1}{p(\sigma d)^p}\log^p(s+e)\right\}\left(\frac{1}{(\sigma d)^p}+1\right)+c_2 &s>0\\
        e^s& s\leq 0.
    \end{cases}
\end{align*}
We note that $(\phi_2^*)'(s)\geq (\psi_2')^{-1}(s)$, for all $s\in \mathbb{R}$. Moreover, $\psi_2'$ is given by
\begin{align*}
     \psi_2'(t)=\begin{cases}
         \exp\left\{\left(p\left(\sigma d\right)^p\log \left(\frac{t-c_2}{c_1\left(\frac{1}{(\sigma d)^p}+1\right)}\right)\right)^{\frac{1}{p}} \right\}-e& t\geq 1\\
         \log(t)& 0< t\leq 1.
     \end{cases} 
 \end{align*}
By Proposition \ref{prop: phi_bound}, we only have to show that the limit \eqref{moment_include} holds. We have that,
\begin{align*}
    \psi'_2\left(\frac{1}{f(x)}\right)&=\psi'_2\left(C\sigma |x|\exp\left\{\frac{1}{p\sigma^p}(\log|x|-\mu)^p\right\}\right)\\
    &=\exp \left\{\left(p(\sigma d)^p\cdot \left(\frac{1}{p\sigma^p}(\log|x|-\mu)^p+\log|x|+\mathcal{O}(1)\right)\right)^{\frac{1}{p}}\right\}\\
    &=\exp \left(\left(d^p\left(\log|x|-\mu\right)^p+\mathcal{O}(\log|x|)\right)^\frac{1}{p}\right)\\
    &=e^{\mathcal{O}(1)}\cdot |x|^d.
\end{align*}
Hence, $\psi_2'$ satisfies the limit \eqref{moment_include}

Similarly, we can also construct a function that bounds $(\phi_2^*)'$ from above. For any $\gamma>d$, we define 
\begin{align*}
    (\hat{\psi}_2^*)'(s)=\begin{cases}
        \max\left\{\exp\left\{\frac{1}{p(\gamma\sigma)^p}(\log^p(s+e)-1)\right\}, (\phi^*_2)'(s)\right\} &s>0\\
        e^s& s\leq 0.
    \end{cases}
\end{align*}
Then, we have that $(\hat{\psi}_2^*)'(s)\geq (\phi^*_2)'(s)$ for all $s\in \mathbb{R}$. Moreover, $\hat{\psi}_2^*(s):=\int^s_0(\hat{\psi}_2^*)'(y)\mathrm{d}y$ is an increasing convex function with $\hat{\psi}_2^*(0)=0, (\hat{\psi}_2^*)'(0)=1$. Hence, $\hat{\psi}_2:=\hat{\psi}_2^{**}$ is also a divergence function, and $(\hat{\psi}^*_2)'=(\hat{\psi}'_2)^{-1}$ by Proposition \ref{prop: construct_phi_conj}.

To show that the divergence $\hat{\psi}_2$ excludes all densities $g$ that satisfy $\liminf_{x\to -\infty}|x|^{t+1}g(x)>0$, for any $t\in (1,d)$. We use Proposition 3 of \citet{KCS19}, for which we have to verify the limit conditions $\limsup_{u\to\infty}\frac{u\hat{\psi}''_2(u)}{\hat{\psi}'_2(u)}<\infty$ and
    \begin{align}\label{moment_exclude}
        \liminf_{x\to -\infty}\frac{\hat{\psi}'_2\left(\frac{1}{f_0(x)}\cdot \frac{1}{|x|^{t+1}}\cdot c\right)}{|x|^t}>0,
    \end{align}
for all constants $c>0$. We note that since $\frac{1}{p(\gamma\sigma)^p}>\frac{1}{p(\sigma d)^p}$ due to $\gamma>d$, we have that $(\hat{\psi}'_2)^{-1}(s)=(\hat{\psi}_2^*)'(s)=\exp\left\{\frac{1}{p(\gamma\sigma)^p}(\log^p(s+e)-1)\right\}$, if $s$ is sufficiently large. This means that for $u$ sufficiently large, we have $\hat{\psi}'_2(u)=\exp\left\{(p(\gamma\sigma)^p\log(u)+1)^{\frac{1}{p}}\right\}-e$. Hence, we have that
\begin{align*}
    \limsup_{u\to\infty}\frac{u\hat{\psi}''_2(u)}{\hat{\psi}'_2(u)}=(\gamma \sigma)^p(p(\gamma\sigma)^p\log(u)+1)^{\frac{1}{p}-1}=0<\infty,
\end{align*}
since $1/p-1<0$ due to $p\geq 2$. We also have that,
for any constant $c>0$:
\begin{align*}
   \hat{\psi}'_2\left(\frac{1}{f_0(x)}\cdot \frac{1}{|x|^{t+1}}\cdot c\right)&=\exp \left(\left(\gamma^p(\log|x|-\mu)^p+\mathcal{O}(1)\cdot \log|x|+\mathcal{O}(1)\right)^{\frac{1}{p}}\right)\\
   &=|x|^\gamma \cdot e^{\mathcal{O}(1)}.
\end{align*}
Therefore, we have that $\lim_{x\to -\infty}|x|^{-t}\hat{\psi}'_2\left(\frac{1}{f_0(x)}\cdot \frac{1}{|x|^{t+1}}\cdot c\right)=\lim_{x\to -\infty}|x|^{\gamma -t}\mathcal{O}(1)>0$, since $\gamma >d>t$. Therefore, we have that $I_{\hat{\psi}_2}(g,f_0)=\infty$, which implies that $I_{\phi_2}(g,f_0)=\infty$, by Proposition \ref{prop: phi_bound}.
\end{proof}

\begin{proof}[\textbf{Proof of Lemma \ref{lem: exp_power}}]
    Following the construction in \eqref{def_phi^*}, we choose $\psi(s)=(s+1)e^{(s+1)^p}$, where $p=\frac{k}{d}$. The second derivative of $\psi$ is given by:
    \begin{align*}
        \psi''(s)=p(s+1)^{p-1}e^{(s+1)^p}(p+p(s+1)^p+1),
    \end{align*}
    which is positive for all $s\geq 0$. Hence, $\psi$ is an increasing, strictly convex and twice continuously differentiable function on $s\geq 0$. By Proposition \ref{prop: construct_phi_conj}, $\phi_2^*$ defined in \eqref{div:exp_power} induces a divergence function.
    We examine condition \eqref{ocep_integral}. Choose an arbitrary $\theta_2\in \mathbb{R}, \theta_1\geq 0$. For $x\to -\infty$, we have that
    \begin{align*}
        \phi_2^*(\phi_1^*(\theta_2-x)+\theta_1)f_0(x)&=\phi_2^*\left(\frac{1}{\alpha}(\theta_2-x)+\theta_1\right)\frac{k}{\lambda}\left(\frac{|x|}{\lambda}\right)^{k-1}e^{-\left(\frac{|x|}{\lambda}\right)^k}\\
        &=\mathcal{O}(|x|^k)\cdot \exp \left(\frac{1}{\alpha}|x|^{\frac{k}{d}}\cdot \mathcal{O}(1)-\left(\frac{|x|}{\lambda}^k\right)\right)+\mathcal{O}\left(|x|^ke^{-\left(\frac{|x|}{\lambda}\right)^k}\right).
    \end{align*}
Since $\frac{k}{d}<k$, we have that the term $\exp\left(\frac{1}{\alpha}|x|^{\frac{k}{d}}\cdot \mathcal{O}(1)-\left(\frac{|x|}{\lambda}^k\right)\right)$ is integrable. Hence, condition \eqref{ocep_integral} is satisfied.

To examine the content of the divergence ball, we first examine the derivative of $\phi_2^*$, which is equal to 
\begin{align*}
    (\phi_2^*)'(s)=\begin{cases}
        c_1e^{(s+1)^{\frac{k}{d}}}\left(1+\frac{k}{d}(s+1)^{\frac{k}{d}}\right)+c_2 &s>0\\
        e^s& s\leq 0.
    \end{cases}
\end{align*}
We then define the auxiliary function 
\begin{align*}
    (\psi_2')^{-1}(s)=\begin{cases}
        c_1e^{(s+1)^{\frac{k}{d}}}\left(1+\frac{k}{d}\right)+c_2 &s>0\\
        e^s& s\leq 0,
    \end{cases}
\end{align*}
which has inverse 
\begin{align*}
    \psi_2'(t)=\begin{cases}
        \log^{\frac{d}{k}}\left(\frac{t-c_2}{c_1(1+k/d)}\right)-1 &t>1\\
        \log(t)& 0<t\leq 1.
    \end{cases}
\end{align*}
We examine the asymptotic behavior:
\begin{align*}
    \psi_2'\left(\frac{1}{f_0(x)}\right)&=\log^{\frac{d}{k}}\left(|x|^{1-k}e^{\left(\frac{|x|}{\lambda}\right)^k}+\mathcal{O}(1)\right)-1\\
    &=\left(\left(\frac{|x|}{\lambda}\right)^k+(1-k)\log |x|+\mathcal{O}(1)\right)^{\frac{d}{k}}-1\\
    &=\mathcal{O}(|x|^d).
\end{align*}
Hence, we have $\limsup_{x\to-\infty}|x|^{-d}\psi_2'\left(\frac{1}{f_0(x)}\right)<\infty$. Proposition \ref{prop: phi_bound} implies that $I_{\phi_2}(g,f_0)<\infty$ for any continuous density $g$ with finite $d$-th moment.

On the other hand, for any $\gamma>1$, we can define the function 
\begin{align*}
    (\hat{\psi}_2^*)'(s)=\begin{cases}
        \max\left\{e^{\gamma(s+1)^{\frac{k}{d}}-\gamma}, (\phi^*_2)'(s)\right\} &s>0\\
        e^s& s\leq 0.
    \end{cases}
\end{align*}
Then, we have that $(\hat{\psi}_2^*)'(s)\geq (\phi^*_2)'(s)$ for all $s\in \mathbb{R}$, and that $\hat{\psi}_2^*(s):=\int^s_0(\hat{\psi}_2^*)'(y)\mathrm{d}y$ is an increasing convex function with $\hat{\psi}_2^*(0)=0, (\hat{\psi}_2^*)'(0)=1$. Hence, $\hat{\psi}_2:=\hat{\psi}_2^{**}$ is a divergence function by Proposition \ref{prop: construct_phi_conj}. Since the term $e^{\gamma(s+1)^{\frac{k}{d}}}$ dominates the term $e^{(s+1)^{\frac{k}{d}}}$ in $(\phi^*_2)'$ as $s\to \infty$, we have that for $s$ sufficiently large, that $(\hat{\psi}_2^*)'(s)=e^{\gamma(s+1)^{\frac{k}{d}}-\gamma}$. Hence, for $u$ sufficiently large, we have that 
\begin{align*}
    \hat{\psi}'_2(u)=\left(\frac{\log(u)+\gamma}{\gamma}\right)^{\frac{d}{k}}-1.
\end{align*}
Therefore, to apply Proposition 3 of \citet{KCS19}, we examine the limit conditions
\begin{align*}
\limsup_{u\to\infty}\frac{u(\hat{\psi}_2)''(u)}{(\hat{\psi}_2)'(u)}=\limsup_{u\to \infty}\frac{d}{k}\left(\frac{\log(u)+\gamma}{\gamma}\right)^{-1}-1<\infty,
\end{align*}
and for any $t\in (1,d)$, any constants $c>0$, we have
\begin{align*}
    \liminf_{x\to -\infty}\hat{\psi}'_2\left(\frac{1}{f_0(x)}\cdot \frac{1}{|x|^{t+1}}\cdot c\right)|x|^{-t}&=\liminf_{x\to -\infty}\left( \left(\left(\frac{|x|}{\lambda}\right)^k-(t+k)\log|x|\right)\cdot \mathcal{O}(1)\right)^{\frac{d}{k}}|x|^{-t}\\
    &=\mathcal{O}(|x|^{d-t})>0.
\end{align*}
Hence, by Proposition \eqref{prop: phi_bound}, we have that $I_{\phi_2}(g,f_0)=\infty$, for any $g$ such that $\liminf_{x\to -\infty}|x|^{t+1}g(x)>0$.
\end{proof}

\begin{proof}[\textbf{Proof of Lemma \ref{lem: entr_weibull}}]
    It follows from the proof of Lemma \ref{lem: div_gl_cvar} that $\phi_2^*$ as in \eqref{div_entr_weibull} is convex and twice differentiable and indeed defines a $\phi$-divergence. It remains to show that $\phi_2^*$ satisfies condition \eqref{ocep_integral} and \eqref{oceb_integral}. To verify this, we choose $\theta_1=\theta_2=0$. Then, we examine
    \begin{align*}
        &\phi_2^*\left(\gamma(\exp\left\{\gamma |x|\right\}-1)\right)\frac{k}{\lambda}\left(\frac{|x|}{\lambda}\right)^{k-1}\exp \left\{-\left(\frac{|x|}{\lambda}\right)^k\right\}\\
        &=c_1(\gamma(\exp\left\{\gamma |x|\right\}-1)+e)\exp\left\{\frac{1}{(2\gamma \lambda)^k}\log^k\left(\gamma(\exp\left\{\gamma |x|\right\}-1)\right)\right\}\\
        &\qquad\cdot \frac{k}{\lambda}\left(\frac{|x|}{\lambda}\right)^{k-1}\exp \left\{-\left(\frac{|x|}{\lambda}\right)^k\right\}+ \mathcal{O}\left(e^{\gamma |x|-\left(\frac{|x|}{\lambda}\right)^k}\right).
    \end{align*}
    Since $k>1$, $\mathcal{O}\left(e^{\gamma |x|-\left(\frac{|x|}{\lambda}\right)^k}\right)$ is integrable. Therefore, we only examine the dominant terms that determine the asymptotics:
    \begin{align*}
        &\exp \left\{\gamma |x|+\frac{1}{(2\gamma \lambda )^k}\log^k\left(\gamma(\exp\left\{\gamma |x|\right\}-1)\right)-\left(\frac{|x|}{\lambda}\right)^k\right\}\\
        &=\exp \left\{\gamma |x|+\left(\frac{|x|}{2\lambda }\right)^k(1+o(1))^k-\left(\frac{|x|}{\lambda}\right)^k\right\}\\
        &=\mathcal{O}\left(\exp\left\{\gamma |x|+\left(\frac{1}{(2\lambda )^k}-\frac{1}{\lambda^k}\right)|x|^k\right\}\right).
    \end{align*}
    Hence, this term is integrable and therefore \eqref{ocep_integral} is satisfied.

    We have that the derivative of $\phi_2^*$ is given by
    \begin{align*}
    (\phi_2^*)'(s)=\begin{cases}
        c_1\exp\left\{\frac{1}{(2\gamma \lambda )^k}\log^k(s+e)\right\}\left(\frac{k}{(2\gamma \lambda)^k}\log^{k-1}(s+e)+1\right)+c_2 &s>0\\
        e^s& s\leq 0.
    \end{cases}
\end{align*}
We define a function
\begin{align*}
    (\psi_2')^{-1}(s)=\begin{cases}
        c_1\exp\left\{\frac{1}{(2\gamma \lambda )^k}\log^k(s+e)\right\}\left(\frac{k}{(2\gamma \lambda )^k}+1\right)+c_2 &s>0\\
        e^s& s\leq 0.
    \end{cases}
\end{align*}
We note that $(\phi_2^*)'(s)\geq (\psi_2')^{-1}(s)$, for all $s\in \mathbb{R}$. Moreover, $\psi_2'$ is given by
\begin{align*}
     \psi_2'(t)=\begin{cases}
         \exp\left\{\left((2\gamma \lambda )^k\log \left(\frac{t-c_2}{c_1\left(\frac{k}{(2\gamma \lambda )^k}+1\right)}\right)\right)^{\frac{1}{k}} \right\}-e& t\geq 1\\
         \log(t)& 0< t\leq 1.
     \end{cases} 
 \end{align*}
 We consider a Weibull distribution $g$ with slightly heavier tail: $g(x)=\frac{l}{\lambda}\left(\frac{|x|}{\lambda}\right)^{l-1}e^{-\left(\frac{|x|}{\lambda}\right)^l}$, where $1<l<k$. Then, we show that $I_{\phi_2}(g,f_0)<\infty$. Following the proof of Proposition \ref{prop: phi_bound}, it is sufficient to show that the following function is integrable on $\{x:|x|>x_0\}$ for some $x_0\in \mathbb{R}$:
 \begin{align*}
     &g(x)\psi_2'\left(\frac{1}{f(x)}\right)\\
     &=\mathcal{O}\left(|x|^{l-1}\exp\left\{-\left(\frac{|x|}{\lambda}\right)^l+(2\gamma \lambda)\left(\left(\frac{|x|}{\lambda}\right)^k+(1-k)\log\left(\frac{|x|}{\lambda}\right)+\mathcal{O}(1)\right)^{\frac{1}{k}}\right\}\right)\\
     &=\mathcal{O}\left(|x|^{l-1}\exp\left\{-\left(\frac{|x|}{\lambda}\right)^l+2\gamma |x|\right\}\right).
 \end{align*}
 Since $l>1$, this is $g(x)\psi_2'\left(\frac{1}{f_0(x)}\right)$ is indeed integrable. Therefore, $I_{\phi_2}(g,f_0)<\infty$.
\end{proof}

\begin{proof}[\textbf{Proof of Theorem \ref{thm: uMean_b}}]
    We express the following optimization problem in terms of density functions:
\begin{align*}
 &\sup_{\lambda>0}\sup_{\substack{\{\mathbb{Q}: I_{\phi_2}(\mathbb{Q}, \mathbb{P}_0)\leq r\}}}\sup_{ \{\bar{\mathbb{Q}}:\bar{\mathbb{Q}}\ll \mathbb{Q}\}} \mathbb{E}_{\bar{\mathbb{Q}}}[-X]-I_{\phi_1}(\bar{\mathbb{Q}},\lambda \mathbb{Q}).\\
 &=-\inf_{\substack{\lambda>0, g, \bar{g}\geq 0,\\ \int_{\Omega} g(\omega)\mathrm{d}\mathbb{P}_0(\omega)=1,\int_{\Omega} \bar{g}(\omega)\mathrm{d}\mathbb{P}_0(\omega)=1\\
\int_{\Omega}\phi_2\left(g(\omega)\right)\mathrm{d}\mathbb{P}_0(\omega)\leq r}}\int_{\Omega}X(\omega)\bar{g}(\omega)+\lambda g(\omega)\phi_1\left(\frac{\bar{g}(\omega)}{\lambda g(\omega)}\right)\mathrm{d}\mathbb{P}_0(\omega).   
\end{align*}
Due to the product $\lambda g(\omega)$, this is a non-convex problem. However, using the substitution $\tilde{g}(\omega)=\lambda g(\omega)$, this is equivalent to the following convex problem
\begin{align*}
    \mu_0:=\inf_{\substack{\lambda>0, \tilde{g}, \bar{g}\geq 0,\\ \int_{\Omega} \tilde{g}(\omega)\mathrm{d}\mathbb{P}_0(\omega)=\lambda,\int_{\Omega} \bar{g}(\omega)\mathrm{d}\mathbb{P}_0(\omega)=1\\
    \int_{\Omega}\lambda \phi_2\left(\frac{\tilde{g}(\omega)}{\lambda}\right)\mathrm{d}\mathbb{P}_0(\omega)\leq \lambda r}}\int_{\Omega}X(\omega)\bar{g}(\omega)+ \tilde{g}(\omega)\phi_1\left(\frac{\bar{g}(\omega)}{\tilde{g}(\omega)}\right)\mathrm{d}\mathbb{P}_0(\omega).
\end{align*}
Assume without loss of generality that $\mu_0$ is finite. To apply Theorem \ref{thm:Luenberger}, we choose $\mathcal{X}=(0,\infty)\times \mathcal{L}^1(\Omega)\times \mathcal{L}^1(\Omega)$, $\mathcal{C}= (0,\infty)\times \mathcal{L}^1_+(\Omega)\times \mathcal{L}^1_+(\Omega)$, $H(\lambda,\tilde{g},\bar{g})=\left(\int_{\Omega}\tilde{g}(\omega)\mathrm{d}\mathbb{P}_0(\omega), \int_{\Omega}\bar{g}(\omega)\mathrm{d}\mathbb{P}_0(\omega)\right)-(\lambda,1)$, $G(\lambda, \tilde{g},\bar{g})= \int_{\Omega}\lambda \phi_2\left(\frac{\tilde{g}(\omega)}{\lambda }\right)\mathrm{d}\mathbb{P}_0(\omega)-\lambda r$. Then, we have that for $f_1\equiv 1$ the constant function, $(1,f_1,f_1)$ is a point such that $H(1,f_1,f_1)=(0,0)$ and $G(1,f_1,f_1)<0$. Furthermore, $(0,0)$ is an interior point of the image set of $H$. Hence, strong duality implies that
\begin{align*}
    \mu_0=\sup_{\eta\geq 0}-\eta r \lambda + \inf_{\substack{\lambda>0, \tilde{g}, \bar{g}\geq 0,\\ \int_{\Omega} \tilde{g}(\omega)\mathrm{d}\mathbb{P}_0(\omega)=\lambda\\\int_{\Omega} \bar{g}(\omega)\mathrm{d}\mathbb{P}_0(\omega)=1}}\int_{\Omega}X(\omega)\bar{g}(\omega)+ \tilde{g}(\omega)\phi_1\left(\frac{\bar{g}(\omega)}{\tilde{g}(\omega)}\right)+\eta \lambda \phi_2\left(\frac{\tilde{g}(\omega)}{\lambda}\right)\mathrm{d}\mathbb{P}_0(\omega).
\end{align*}
We first study the case when the supremum is attained at $\eta = 0$, in which case we have that
\begin{align*}
    \mu_0&=\inf_{\substack{\lambda>0, \tilde{g}, \bar{g}\geq 0,\\ \int_{\Omega} \tilde{g}(\omega)\mathrm{d}\mathbb{P}_0(\omega)=\lambda\\\int_{\Omega} \bar{g}(\omega)\mathrm{d}\mathbb{P}_0(\omega)=1}}\int_{\Omega}X(\omega)\bar{g}(\omega)+ \tilde{g}(\omega)\phi_1\left(\frac{\bar{g}(\omega)}{\tilde{g}(\omega)}\right)\mathrm{d}\mathbb{P}_0(\omega)\\
    &\geq \inf_{\substack{\lambda>0, \tilde{g}, \bar{g}\geq 0,\\ \int_{\Omega} \tilde{g}(\omega)\mathrm{d}\mathbb{P}_0(\omega)=\lambda\\\int_{\Omega} \bar{g}(\omega)\mathrm{d}\mathbb{P}_0(\omega)=1}}\int_{\Omega}X(\omega)\bar{g}(\omega)\mathrm{d}\mathbb{P}_0(\omega)=\essinf (X).
\end{align*}
On the other hand, taking $\lambda=1, \tilde{g}=\bar{g}$ shows that $\essinf(X)$ is also an upper bound of $\mu_0$. Hence, $\mu_0=\essinf(X)$. The dual formula \eqref{dual_sf_risk_measure} states that when $\eta =0$, $\essinf(X)$ is equal to
\begin{align*}
    \sup_{\theta_1,\theta_2\in \mathbb{R}}\left\{-\theta_2~\middle|~ \mathbb{E}_{\mathbb{P}_0}\left[0 \phi_2^*\left(\frac{\phi_1^*(-\theta_2-X)+\theta_1}{0}\right)\right]\leq \theta_1\right\},
\end{align*}
which implies that $0\leq \theta_1\leq -\phi_1^*(-\theta_2-\essinf(X))$. This inequality is satisfied if we take $-\theta_2= \essinf(X)$ and $\theta_1=0$, but violated if $-\theta_2>\essinf(X)$. Hence, the duality formula \eqref{dual_sf_risk_measure} holds for $\eta=0$, and we only have to examine the case where the supremum is attained at $\eta>0$, where strong duality implies
\begin{align*}
    \mu_0=\sup_{\substack{\eta> 0\\ \theta_1,\theta_2\in \mathbb{R}}}\inf_{\substack{\lambda>0\\\tilde{g},\bar{g}\geq 0}}(-\eta r+\theta_1)\lambda +\theta_2+\int_{\Omega}\beta_{\eta,\lambda}(\omega,\tilde{g}(\omega),\bar{g}(\omega))-\theta_1\tilde{g}(\omega)-\theta_2\bar{g}(\omega)\mathrm{d}\mathbb{P}_0(\omega),
\end{align*}
where 
\begin{align*}
    \beta_{\eta,\lambda}(\omega,s,t)=X(\omega)t+s\phi_1\left(\frac{t}{s}\right)+\lambda \eta \phi_2\left(\frac{s}{\lambda }\right).
\end{align*}
Note that when $\eta>0$, $\beta_{\eta,\lambda}(\omega,s,t)$ is only finite when $s,t\geq 0$. Analogous to the proof of and Theorem \ref{thm:J_beta}, we have
\begin{align*}
 &\inf_{\tilde{g},\bar{g}\geq 0}  \int_{\Omega}\beta_{\eta,\lambda}(\omega,\tilde{g}(\omega),\bar{g}(\omega))-\theta_1\tilde{g}(\omega)-\theta_2\bar{g}(\omega)\mathrm{d}\mathbb{P}_0(\omega) \\
 &=-\int_{\Omega}\beta^*_{\eta,\lambda}(\omega,\theta_1,\theta_2)\mathrm{d}\mathbb{P}_0(\omega),
\end{align*}
where $\int_{\Omega}\beta^*_{\eta,\lambda}(\omega,\theta_1,\theta_2)\mathrm{d}\mathbb{P}_0(\omega)=+\infty$ if $\eta=0$ and for $\eta>0$, we have
\begin{align*}
 \int_{\Omega}\beta^*_{\eta,\lambda}(\omega,\theta_1,\theta_2)\mathrm{d}\mathbb{P}_0(\omega)=\int_{\Omega}\lambda \eta \phi_2^*\left(\frac{\theta_1+\phi_1^*(\theta_2-X(\omega))}{\eta}\right)\mathrm{d}\mathbb{P}_0(\omega). 
\end{align*}
Hence, we have
\begin{align*}
    \mu_0&=\sup_{\substack{\eta> 0\\ \theta_1,\theta_2\in \mathbb{R}}}\theta_2+\inf_{\lambda>0}\left(\theta_1-\eta r-\int_{\Omega} \eta \phi_2^*\left(\frac{\theta_1+\phi_1^*(\theta_2-X(\omega))}{\eta}\right)\mathrm{d}\mathbb{P}_0(\omega)\right)\lambda\\
    &=\sup_{\substack{\eta> 0\\\theta_1,\theta_2\in \mathbb
    R}}\left\{\theta_2~\middle |~\theta_1-\eta r-\int_{\Omega}\eta \phi_2^*\left(\frac{\theta_1+\phi_1^*(\theta_2-X(\omega))}{\eta}\right)\mathrm{d}\mathbb{P}_0(\omega)\geq 0\right\}.
\end{align*}
\end{proof}

\begin{proof}[\textbf{Proof of Proposition \ref{prop:OCEMeasure}}]
By Proposition \ref{prop: compactOCE}, there exists a sufficiently large number $L>0$, such that 
\begin{align*}
  \rho^l_{\mathrm{oce},\mathbb{P}}(X_p)&=\inf_{\theta_1,\theta_2\in [-L,L]}-\theta_1-\theta_2+ p(\phi_2^*(\phi_1^*(\theta_2-x)+\theta_1)+(1-p)\phi_2^*(\phi_1^*(\theta_2)+\theta_1).
\end{align*}
We have that the minimizer of $\rho^l_{\mathrm{oce},\mathbb{P}}(X_p)$ is attained at $(0,0)$, if and only if, for all $ \theta_1,\theta_2\in[-L,L]$
\begin{align*}
     -\theta_1-\theta_2+p(\phi_2^*(\phi_1^*(\theta_2-x)+\theta_1))+(1-p)(\phi_2^*(\phi_1^*(\theta_2)+\theta_1))\geq p\phi_2^*(\phi_1^*(-x)),
\end{align*}
which is if and only if
\begin{align*}
    -\theta_1-\theta_2+\phi_2^*(\phi_1^*(\theta_2)+\theta_1)\geq p(\phi_2^*(\phi_1^*(-x))+\phi_2^*(\phi_1^*(\theta_2)+\theta_1)-\phi_2^*(\phi_1^*(\theta_2-x)+\theta_1)).
\end{align*}
We note the left-hand side is always strictly positive unless at $(\theta_1,\theta_2) =(0,0)$ and the right-hand side is bounded above by $pM$ for some constant $M>0$, uniformly for all $ \theta_1,\theta_2\in[-L,L]$, due to the continuity of the convex functions $\phi_1^*,\phi_2^*$ on their domain $\mathbb{R}$. Denote $F(\theta_1,\theta_2)=-\theta_1-\theta_2+\phi_2^*(\phi_1^*(\theta_2)+\theta_1)$. Then, for any $p$, the minimizer of $\rho^l_{\mathrm{oce},\mathbb{P}}(X_p)$ must be contained in the level set $L_p:=\{(\theta_1,\theta_2)\in [-L,L]^2: F(\theta_1,\theta_2)\leq pM\}$, since for any $(\theta_1,\theta_2)$ outside of the set $L_p$, the above inequality is satisfied and $(0,0)$ is a better solution. We argue that $\sup_{(\theta_1,\theta_2)\in L_p}\|(\theta_1,\theta_2)\|\to 0$ as $p\to 0$. Indeed, if not, then for any $p_n\to 0$, by the compactness of $[-L,L]^2$ and continuity of $F$, there exists a convergent sequence $(\theta_{1,n},\theta_{2,n})\in L_{p_n}$ with limit $(\theta_{1,L},\theta_{2,L})$, such that $\|(\theta_{1,n},\theta_{2,n})\|>\delta$ for some $\delta>0$ and all $n\geq 1$, but $F(\theta_{1,L},\theta_{2,L})=\lim_{n\to \infty}F(\theta_{1,n},\theta_{2,n})\leq \lim_{n\to \infty}p_nM=0$. This is contradictory since $F$ has the unique minimum at $(0,0)$. Hence, for any $p>0$, the minimizers $(\theta^*_{1,p},\theta^*_{2,p})$ of $\rho^l_{\mathrm{oce},\mathbb{P}}(X_p)$ converges to $(0,0)$. Thus, continuity implies that $\rho^l_{\mathrm{oce},\mathbb{P}}(X_p)- p\phi_2^*(\phi_1^*(-x))\to 0$.

To show \eqref{uMean_elicit}, we first note that by weak duality, we have $\rho^l_{\mathrm{oce},\mathbb{P}}(X_p)\leq \rho^l_{\mathrm{sf},\mathbb{P}}(X_p)$ for all $p>0$. Therefore, we already have
\begin{align*}
    \liminf_{p\downarrow 0}\frac{\rho^l_{\mathrm{sf},\mathbb{P}}(X_p)}{p}\geq \phi_2^*(\phi_1^*(-x)).
\end{align*}
We will now show the reverse. By definition, we have that
\begin{align*}
    \frac{1}{p}\rho^l_{\mathrm{sf},\mathbb{P}}(X_p)&=\inf_{\theta_1,\theta_2\in \mathbb{R}}\left\{\frac{1}{p}\theta_2~\middle |~p\phi_2^*(\phi_1^*(-\theta_2-x)+\theta_1)+(1-p)\phi_2^*(\phi_1^*(-\theta_2)+\theta_1)\leq \theta_1\right\}\\
    &=\inf_{\theta_1,\theta_2\in \mathbb{R}}\left\{\theta_2~\middle |~p\phi_2^*(\phi_1^*(-p\theta_2-x)+\theta_1)+(1-p)\phi_2^*(\phi_1^*(-p\theta_2)+\theta_1)\leq \theta_1\right\}.
\end{align*}
We bound the above infimum by consider only cases when $\theta_1=0$. Then, we show that as $p\downarrow 0$, the following inequality holds for any $\theta_2>\phi_2^*(\phi_1^*(-x))$:
\begin{align*}
    p\phi_2^*(\phi_1^*(-p\theta_2-x))+(1-p)\phi_2^*(\phi_1^*(-p\theta_2))\leq 0.
\end{align*}
The above inequality holds if and only if 
\begin{align*}
    \phi_2^*(\phi_1^*(-p\theta_2-x))-\phi_2^*(\phi_1^*(-p\theta_2))\leq \theta_2\cdot \frac{\phi_2^*(\phi_1^*(-p\theta_2))}{-p\theta_2}.
\end{align*}
As $p\downarrow 0$, the left-hand side converges to $\phi_2^*(\phi_1^*(-x))$, and the right-hand side converges to $\theta_2\cdot \lim_{t\to0}\frac{\phi_2^*(\phi_1^*(t))}{t}=\theta_2\cdot 1=\theta_2$. Since $\theta_2>\phi_2^*(\phi_1^*(-x))$, the above inequality is feasible for $p$ sufficiently small. Therefore, we also have $\limsup_{p\downarrow 0}\frac{\rho^l_{\mathrm{sf},\mathbb{P}}(X_p)}{p}\leq \phi_2^*(\phi_1^*(-x))$.
\end{proof}

\begin{proof}[\textbf{Proof of Theorem \ref{thm:GRO_ref_discrete}}]
We have that \eqref{primal_gro} is equivalent (up to a sign convention) to the following optimization problem over density functions:
\begin{align}\label{primal_gro}
\begin{split}
\mathrm{Minimize}&~\int_{\Omega}X(\omega)\bar{g}(\omega)+\tilde{g}(\omega)\phi_1\left(\frac{\bar{g}(\omega)}{\tilde{g}(\omega)}\right)+g(\omega)\phi_2\left(\frac{\tilde{g}(\omega)}{g(\omega)}\right)\mathrm{d}\mathbb{P}_0(\omega)\\
    \text{subject to}&~\int_{\Omega} g(\omega)\mathrm{d}\mathbb{P}_0(\omega)=\int_{\Omega} \tilde{g}(\omega)\mathrm{d}\mathbb{P}_0(\omega)=\int_{\Omega} \bar{g}(\omega)\mathrm{d}\mathbb{P}_0(\omega)=1\\
    &\int_\Omega \phi_3(g(\omega))\mathrm{d}\mathbb{P}_0(\omega)\leq r\\
    &g(\omega), \tilde{g}(\omega), \bar{g}(\omega)\geq 0.
\end{split}
\end{align}
We denote the optimal value of the primal problem \eqref{primal_gro} as $J_0$. Note that $J_0<\infty$ since the constant unity function $g,\tilde{g},\bar{g}\equiv 1$ is a feasible solution and we assumed $\int_\Omega X(\omega)\mathrm{d}\mathbb{P}_0(\omega)<\infty$. Strong duality then implies that
\begin{align*}
    J_0=&\sup_{\substack{\theta_1,\theta_2,\theta_3\in \mathbb{R}\\ \lambda \geq 0}}-\lambda r+\inf_{\substack{g,\tilde{g},\bar{g}\geq 0\\\int g=\int \tilde{g}=\int\bar{g}=1}}\int_{\Omega}X(\omega)\bar{g}(\omega)+\tilde{g}(\omega)\phi_1\left(\frac{\bar{g}(\omega)}{\tilde{g}(\omega)}\right)+g(\omega)\phi_2\left(\frac{\tilde{g}(\omega)}{g(\omega)}\right)\mathrm{d}\\
    &\qquad +\lambda\phi_3(g(\omega)) \mathbb{P}_0(\omega).
\end{align*}
Suppose this supremum is attained at $\lambda =0$, then we would have $J_0=\essinf(X)$. The dual formula in \eqref{dual_GRO_discrete} would also yield
\begin{align*}
    &\sup_{\theta_1,\theta_2,\theta_3\in \mathbb{R}}\theta_1+\theta_2+\theta_3-\mathbb{E}_{\mathbb{P}_0}\left[0\phi_3^*\left(\frac{\phi_2^*(\phi_1^*(\theta_3-X(\omega))+\theta_2)+\theta_1}{0}\right)\right]\\
    &=\sup_{\theta_2,\theta_3\in \mathbb{R}}\theta_2+\theta_3-\phi_2^*(\phi_1^*(\theta_2+\theta_3-\essinf(X)))\\
    &=\essinf(X).
\end{align*}
Hence, the duality formula holds if the supremum is attained at $\lambda=0$. We may now assume that the supremum is at $\lambda>0$.
Applying Theorem \ref{thm:Luenberger} and Theorem 14.60 of \citet{Rockafellar&Wets}, gives the duality
\begin{align*}
    J_0
    &=\sup_{\substack{\theta_1,\theta_2,\theta_3\in \mathbb{R}\\ \lambda > 0}}-\lambda r +\theta_1+\theta_2+\theta_3-\int_\Omega \beta_{\lambda}^*(\omega,\theta_1,\theta_2,\theta_3)\mathrm{d}\mathbb{P}_0(\omega),
\end{align*}
where 
\begin{align*}
\beta_\lambda(\omega,t_1,t_2,t_3)=X(\omega)t_3+t_2\phi_1\left(\frac{t_3}{t_2}\right)+t_1\phi_2\left(\frac{t_2}{t_1}\right)+\lambda \phi_3(t_1).
\end{align*}
We note that $\beta_\lambda(\omega,t_1,t_2,t_3)$ is only finite if $t_1,t_2,t_3\geq 0$, since $\lambda>0$. Therefore, we have
\begin{align*}
    \beta_{0}^*(\omega,\theta_1,\theta_2,\theta_3)&=\sup_{t_1,t_2,t_3\geq 0}\theta_1t_1+\theta_2t_2+\theta_3t_3-X(\omega)t_3-t_2\phi_1\left(\frac{t_3}{t_2}\right)-t_1\phi_2\left(\frac{t_2}{t_1}\right)-\lambda \phi_3(t_1)\\
    &=\sup_{t_1\geq 0, t_2\geq 0}\theta_1t_1+\theta_2t_2+t_2\phi_1^*(\theta_3-X(\omega))-t_1\phi_2\left(\frac{t_2}{t_1}\right)-\lambda \phi_3(t_1)\\
    &=\sup_{t_1\geq 0}\theta_1t_1+t_1\phi_2^*(\phi_1^*(\theta_3-X(\omega))+\theta_2)-\lambda \phi_3(t_1)\\
    &=\lambda \phi_3^*\left(\frac{\phi_2^*(\phi_1^*(\theta_3-X(\omega))+\theta_2)+\theta_1}{\lambda}\right).
\end{align*}

\end{proof}

\begin{proof}[\textbf{Proof of Proposition \ref{prop: compactOCE}}]
    Denote
    \begin{align*}
        K(\theta_1,\theta_2):=-\theta_1-\theta_2+\mathbb{E}_{\mathbb{P}_0}\left[\phi_2^*(\phi_1^*(\theta_2-X)+\theta_1)\right].
    \end{align*}
    Since $\phi_2(1)=0$, we have that $\phi_2^*(s)\geq s,\forall s\in \mathbb{R}$. Therefore, for any $\theta_1,\theta_2\in \mathbb{R}$:
    \begin{align*}
         K(\theta_1,\theta_2)\geq -\theta_2+\mathbb{E}_{\mathbb{P}_0}\left[\phi_1^*(\theta_2-X)\right].
    \end{align*}
By assumption the of $\phi_2\in \Phi_0$, we have that $\mathrm{dom}(\phi_2)$ has a non-empty interior around $1$. Therefore, let $x_0<1<y_0$ be such that $\phi_2(x_0),\phi_2(y_0)<\infty$. Then, we have that 
\begin{align*}
    \phi_1^*(\theta_2-X)&\geq \max\{x_0(\theta_2-X)-\phi_1(x_0),\\
    &\qquad y_0(\theta_2-X)-\phi_1(y_0)\}.
\end{align*}
Hence,
\begin{align*}
    K(\theta_1,\theta_2)&\geq \max\{(x_0-1)\theta_2-\phi_1(x_0)-x_0\mathbb{E}_{\mathbb{P}_0}[X],\\
    &\qquad (y_0-1)\theta_2-\phi_1(y_0)-x_0\mathbb{E}_{\mathbb{P}_0}[X]\}.
\end{align*}
By assumption, there exists a $(\theta_1^*,\theta_2^*)$ such that $K(\theta_1^*,\theta_2^*)<\infty$. We now claim that $\theta_2$ can be restricted on the set
\begin{align}\label{interval}
    \left[\frac{1}{x_0-1}(\phi_1(x_0)+x_0\mathbb{E}_{\mathbb{P}_0}[X]+K(\theta_1^*,\theta_2^*)),~ \frac{1}{y_0-1}(\phi_1(y_0)+y_0\mathbb{E}_{\mathbb{P}_0}[X]+K(\theta_1^*,\theta_2^*))\right].
\end{align}
Indeed, if $\theta_2<\frac{1}{x_0-1}(\phi_1(x_0)+x_0\mathbb{E}_{\mathbb{P}_0}[X]+K(\theta_1^*,\theta_2^*))$, then since $(x_0-1)<0$, we have that 
\begin{align*}
    K(\theta_1,\theta_2)&\geq (x_0-1)\theta_2-\phi_1(x_0)-x_0\mathbb{E}_{\mathbb{P}_0}[X]\\
    &>K(\theta_1^*,\theta_2^*).
\end{align*}
Similarly, if $\theta_2>\frac{1}{y_0-1}(\phi_1(y_0)+y_0\mathbb{E}_{\mathbb{P}_0}[X]+K(\theta_1^*,\theta_2^*))$, then $K(\theta_1,\theta_2)>K(\theta_1^*,\theta_2^*)$.

Now, let $L_2>0$ be a number such that $[-L_2,L_2]$ contains \eqref{interval}. We may thus restrict $\theta_2$ on $[-L_2,L_2]$. Then, by the non-decreasing property of $\phi_1^*, \phi_2^*$, we have that
\begin{align}
    \notag K(\theta_1,\theta_2)&\geq -\theta_1-L_2+\mathbb{E}_{\mathbb{P}_0}\left[\phi_2^*(\phi_1^*(-L_2-X)+\theta_1)\right]\\
    \label{bound_t1}&\geq -\theta_1-L_2+\mathbb{E}_{\mathbb{P}_0}\left[\phi_2^*(-L_2-X+\theta_1)\right]
\end{align}
Again, let $\tilde{x}_0<1<\tilde{y}_0$ be such that $\phi_2^*(\tilde{x}_0),\phi_2^*(\tilde{y}_0)<\infty$. Then, with exactly the same argument, we may restrict $\theta_1$ on the interval
\begin{align*}
    &\left[\frac{1}{\tilde{x}_0-1}\left((\tilde{x}_0+1)L_2+\tilde{x}_0\mathbb{E}_{\mathbb{P}_0}[X]+\phi_2(\tilde{x}_0)+K(\theta_1^*,\theta_2^*)\right)\right.\\
    &\qquad\left.,~\frac{1}{\tilde{y}_0-1}\left((\tilde{y}_0+1)L_2+\tilde{y}_0\mathbb{E}_{\mathbb{P}_0}[X]+\phi_2(\tilde{y}_0)+K(\theta_1^*,\theta_2^*)\right)\right]
\end{align*}
Hence, $\theta_1,\theta_2$ can be both restricted on compact sets, without changing the optimum. We now examine 
\begin{align*}
    V(\theta_1,\theta_2,\lambda):=-\theta_1-\theta_2+\lambda r +\mathbb{E}_{\mathbb{P}_0}\left[\lambda \phi_2^*\left(\frac{\phi_1^*(\theta_2-X)+\theta_1}{\lambda}\right)\right].
\end{align*}
We note that since $\lambda \phi_2^*\left(\frac{s}{\lambda}\right)=\sup_{t\geq 0}\{st-\lambda \phi_2(s)\}=(\lambda \phi_2)^*(s)$ and $\phi_2\geq 0$, we have that $\lambda \phi_2^*\left(\frac{s}{\lambda}\right)$ is non-increasing in $\lambda$, and that $(\lambda \phi_2)^*(s)\geq s,~\forall s\in \mathbb{R}$. Therefore, for any $\theta_1,\theta_2\in \mathbb{R}$ and $\lambda> 0$, we have
\begin{align*}
    V(\theta_1,\theta_2,\lambda)&= -\theta_1-\theta_2+\lambda r +\mathbb{E}_{\mathbb{P}_0}\left[(\lambda\phi_2)^*( \phi_1^*(\theta_2-X)+\theta_1)\right]\\
    &\geq -\theta_1-\theta_2+\lambda r +\mathbb{E}_{\mathbb{P}_0}\left[\phi_1^*(\theta_2-X)+\theta_1\right]\\
    &\geq \lambda r-\mathbb{E}_{\mathbb{P}_0}[X].
\end{align*}
Since $V(\theta_1^*,\theta_2^*,1)=K(\theta_1^*,\theta_2^*)<\infty$, we have that $\lambda$ may be restricted on the set
\begin{align*}
    \left[0, \frac{V(\theta_1^*,\theta_2^*,1)+\mathbb{E}_{\mathbb{P}_0}[X]}{r}\right].
\end{align*}
Let $L_{\lambda}$ denotes this upper bound of $\lambda$. Then, we have that for all $\lambda \in [0, L_{\lambda}]$,
\begin{align*}
    V(\theta_1,\theta_2,\lambda)&\geq -\theta_1-\theta_2 +\mathbb{E}_{\mathbb{P}_0}\left[L_{\lambda} \phi_2^*\left(\frac{\phi_1^*(\theta_2-X)+\theta_1}{L_{\lambda}}\right)\right]\\
    &\geq -\theta_2 +\mathbb{E}_{\mathbb{P}_0}\left[\phi_1^*(\theta_2-X)\right].
\end{align*}
Hence, $\theta_2$ can be bounded by the same interval as in \eqref{interval}, thus we may assume that $\theta_2\in [-L_2,L_2]$ for some $L_2>0$. Then, we also have
\begin{align*}
    V(\theta_1,\theta_2,\lambda)&\geq -\theta_1-L_2 +\mathbb{E}_{\mathbb{P}_0}\left[L_{\lambda} \phi_2^*\left(\frac{\phi_1^*(-L_2-X)+\theta_1}{L_{\lambda}}\right)\right]\\
    &\geq -\theta_1-L_2 +\mathbb{E}_{\mathbb{P}_0}\left[L_{\lambda} \phi_2^*\left(\frac{-L_2-X+\theta_1}{L_{\lambda}}\right)\right]\\
    &=-\theta_1-L_2 +\mathbb{E}_{\mathbb{P}_0}\left[ (L_{\lambda}\phi_2)^*\left(-L_2-X+\theta_1\right)\right].
\end{align*}
Since $L\phi_2$ is also a divergence function belonging to $\Phi_0$. It follows from the analysis of \eqref{bound_t1} that $\theta_1$ can also be bounded in the interval
\begin{align*}
    &\left[\frac{1}{\tilde{x}_0-1}\left((\tilde{x}_0+1)L_2+\tilde{x}_0\mathbb{E}_{\mathbb{P}_0}[X]+L_{\lambda}\phi_2(\tilde{x}_0)+K(\theta_1^*,\theta_2^*)\right)\right.\\
    &\qquad\left.,~\frac{1}{\tilde{y}_0-1}\left((\tilde{y}_0+1)L_2+\tilde{y}_0\mathbb{E}_{\mathbb{P}_0}[X]+L_{\lambda}\phi_2(\tilde{y}_0)+K(\theta_1^*,\theta_2^*)\right)\right],
\end{align*}
where $\tilde{x}_0<1<\tilde{y}_0$ are such that $\tilde{x}_0,\tilde{y}_0 \in \mathrm{dom}(\phi_2)$.
\end{proof}

\begin{proof}[\textbf{Proof of Theorem \ref{thm: ConvergenceSAA}}]
We verify the six conditions of Theorem 5.4 in \citet{LectureSP}. We have that \eqref{dual_OCE_p} is equal to 
\begin{align*}
    \inf_{\theta_1,\theta_2\in \mathbb{R}}\mathbb{E}_{\mathbb{P}_0}\left[\phi_2^*(\phi_1^*(\theta_2-X)+\theta_1)-\theta_1-\theta_2\right]=:\inf_{\boldsymbol{\theta}\in \mathbb{R}^2}\mathbb{E}_{\mathbb{P}_0}[F(\boldsymbol{\theta},X)],
\end{align*}
and \eqref{dual_OCE_b} is equal to
\begin{align*}
    \inf_{\lambda \geq 0,\theta_1,\theta_2\in \mathbb{R}}\mathbb{E}_{\mathbb{P}_0}\left[\lambda \phi_2^*\left(\frac{\phi_1^*(\theta_2-X)+\theta_1}{\lambda}\right)-\theta_1-\theta_2+\lambda r\right]=:\inf_{\lambda \geq 0,\boldsymbol{\theta}\in \mathbb{R}^2}\mathbb{E}_{\mathbb{P}_0}[G(\boldsymbol{\theta},\lambda, X)],
\end{align*}
where we replaced $\lambda>0$ with $\lambda\geq 0$, which does not change the optimization problem.

By our assumptions $\mathrm{dom}(\phi_2^*)=\mathrm{dom}(\phi_1^*)=\mathbb{R}$, it follows from convexity that $\phi_1^*,\phi_2^*$ are continuous on $\mathbb{R}$. Therefore, the functions $F(\boldsymbol{\theta},X(\omega))$ and $G((\boldsymbol{\theta},\lambda,X(\omega))$ are convex, lower-semicontinuous functions in $\boldsymbol{\theta}$ or $(\boldsymbol{\theta},\lambda)$ with non-empty interior effective domain, for each $\omega\in \Omega$. Clearly, for each $\boldsymbol{\theta}$ or $(\boldsymbol{\theta},\lambda)$, they are also measurable functions in $\omega$. Hence, it follows from Proposition 14.39 of \citet{Rockafellar&Wets} that they are normal convex integrand. Therefore, conditions (i) and (ii) of Theorem 5.4 are satisfied. The sets $\mathbb{R}^2$ and $\mathbb{R}^2\times [0,\infty)$ are also closed and convex. Thus, (iii) of Theorem 5.4 is satisfied. We also have that for any $\boldsymbol{\theta}\in \mathbb{R}^2$ and $\omega\in \Omega$, that $F(\boldsymbol{\theta},X(\omega))\geq -X(\omega)$ and for any $\boldsymbol{\theta}\in \mathbb{R}^2, \lambda\geq 0$ and $\omega\in \Omega$, that $G(\boldsymbol{\theta},\lambda,X(\omega))\geq -X(\omega)+\lambda r\geq -X(\omega)$. Since $X$ is assumed to be $\mathbb{P}_0$-integrable, it follows from Theorem 7.42 of \citet{LectureSP} that condition (iv) of Theorem 5.4 holds. Condition (v), which requires boundedness and non-emptiness of the set of solutions follows from Proposition \ref{prop: compactOCE} and the assumption. Finally, condition (vi), requires that the law of large number must hold pointwise for $F(\boldsymbol{\theta},X)$ and $G(\boldsymbol{\theta},\lambda,X)$. For cases, where $\mathbb{E}_{\mathbb{P}_0}[F(\boldsymbol{\theta},X)]<\infty$ or $\mathbb{E}_{\mathbb{P}_0}[G(\boldsymbol{\theta},\lambda,X)]<\infty$, this is trivial. For other cases, this holds by Theorem 2.4.5 of \citet{Durrett}, due to the fact that $F(\boldsymbol{\theta},X)\geq -X$ and  $G(\boldsymbol{\theta},\lambda,X)\geq -X$, with $\mathbb{E}_{\mathbb{P}_0}[-X]<\infty$.
\end{proof}

\section{Minimizing Convex Integral}
Let $(\Omega,\mathcal{F},\mathbb{P})$ be a probability space. Consider the convex integral minimization problem over integrable functions.
\begin{align}\label{J_beta}
    J_\beta(\mathbf{a}):=\inf_{\mathbf{g}\in \mathcal{L}^1_+(\mathbf{a})}\int_{\Omega}\beta(\omega,\mathbf{g}(\omega))\mathrm{d}\mathbb{P}(\omega),
\end{align}
where $\mathbf{a}\in \mathbb{R}^d$ and the set of $d$-dimensional non-negative $\mathbb{P}$-integrable functions with value $\mathbf{a}$, denoted as: 
\begin{align*}
\mathcal{L}^1_+(\mathbf{a})=\left\{\mathbf{g}=(g_1,\ldots,g_d)~\middle|~g_i(\omega)\geq 0,~\int_{\Omega}g_i(\omega)\mathrm{d}\mathbb{P}(\omega)=a_i, \forall i=1,\ldots,d\right\}.    
\end{align*}
Here, $\beta(\omega,\mathbf{y}): \mathbb{R}^d\to \mathbb{R}$ is a convex function in $\mathbf{y}$ for $\mathbb{P}$-almost every $\omega$ and a normal convex integrand in the sense of \citet{Rockafellar&Wets}. Furthermore, we assume that $\mathrm{dom}(\beta(\omega,\mathbf{y}))\subset [0,\infty)^d$ for all $\omega\in \Omega$. Therefore, \eqref{J_beta} is also equivalent to 
\begin{align}\label{J_beta1}
J_\beta(\mathbf{a})=\inf_{\mathbf{g}\in\mathcal{L}^1(\mathbf{a})}\int_{\Omega}\beta(\omega,\mathbf{g}(\omega))\mathrm{d}\mathbb{P}(\omega),
\end{align}
where $\mathcal{L}^1(\mathbf{a})$ is the set of all 
$\mathbb{P}$-integrable functions with integral equals to $\mathbf{a}$:
\begin{align*}
  \mathcal{L}^1(\mathbf{a})=\left\{\mathbf{g}=(g_1,\ldots,g_d)~\middle|~\int_{\Omega}g_i(\omega)\mathrm{d}\mathbb{P}(\omega)=a_i, \forall i=1,\ldots,d\right\}.   
\end{align*}
In the sequel, we let $\mathcal{L}^1:=\bigcup_{\mathbf{a}\in \mathbb{R}^d}\mathcal{L}^1(\mathbf{a})$ denote the set of all $\mathbb{P}$-integrable functions.
Problem \eqref{J_beta1} can be reformulated using Lagrangian duality. Since we are optimizing over density functions, which belong to an infinite-dimensional vector space of integrable functions, we invoke the following Lagrangian duality theorem given in \citet{Luenberger}, which we state here for completeness.
\begin{mytheorem}[\citealt{Luenberger}, exercise 8.8.7]\label{thm:Luenberger}
Let $\mathcal{X}$ be a real vector space, $\mathcal{C}$ a convex subset of $\mathcal{X}$. Let $f$ be a real-valued convex functional on $\mathcal{C}$ and $G:\mathcal{C}\to \mathbb{R}^d$ a convex mapping, i.e., $G(\lambda x_1+(1-\lambda x_2))\leq_{\mathbb{R}^d} \lambda G(x_1)+(1-\lambda)G(x_2), \forall x_1,x_2\in \mathcal{C}, \lambda\in (0,1)$ ($x\leq_{\mathbb{R}^d}y \Leftrightarrow x_i\leq y_i,\forall i$). Let $H:\mathcal{X}\to \mathbb{R}^m$ be an affine mapping, i.e, $H(x)=Ax+b$ for some linear mapping $A$. Consider the optimization problem,
\begin{align*}
    \mu_0=\inf_{x\in \mathcal{X}}\{f(x)~|~ x\in \mathcal{C}, G(x)\leq \mathbf{0}, H(x)=\mathbf{0}\}, 
\end{align*}
where $\mathbf{0}$ denotes the zero-vector. Assume the existence of a $x_1\in \mathcal{C}$ such that $G(x_1)<_{\mathbb{R}^d}\mathbf{0}$ and $H(x_1)=\mathbf{0}$. Suppose also that $\mathbf{0}$ is an interior point of the image set $\{y\in \mathbb{R}^m~|~ H(x)=y~\text{for some}~x\in \mathcal{C}\}$. If $\mu_0$ is finite, then there exists a $\boldsymbol{\mu}^*\in \mathbb{R}^d$ with $\mu\geq_{\mathbb{R}^d} \mathbf{0}$ and $\boldsymbol{\nu}^*\in \mathbb{R}^m$, such that
\begin{align*}
    \mu_0&=\inf_{x\in \mathcal{C}}\{f(x)+\boldsymbol{\mu}^{*T}G(x)+\boldsymbol{\nu}^{*T}H(x)\}\\
    &=\sup_{\substack{\boldsymbol{\mu}\geq_{\mathbb{R}^d} \mathbf{0}\\\boldsymbol{\nu}\in \mathbb{R}^m}}\inf_{x\in \mathcal{C}}\{f(x)+\boldsymbol{\mu}^{T}G(x)+\boldsymbol{\nu}^{T}H(x)\}.
\end{align*}
\end{mytheorem}
To apply Theorem \ref{thm:Luenberger}, we set $\mathcal{X}=(\mathcal{L}^1)^d$ and $H(\mathbf{g})=\int_{\Omega}\mathbf{g}(\omega)\mathrm{d}\mathbb{P}(\omega)-\mathbf{a}$. Note that the zero vector $\mathbf{0}\in \mathbb{R}^d$ is an interior point of the image set of $H$, since any vector $\boldsymbol{\epsilon}\in \mathbb{R}^d$ sufficiently close to $\mathbf{0}$ is the image of $H(\mathbf{a}+\boldsymbol{\epsilon})$ of the constant function $\mathbf{g}(\omega)=\mathbf{a}+\boldsymbol{\epsilon}$. Hence, Theorem \ref{thm:Luenberger} can be applied to $J_\beta(\mathbf{a})$ if $J_\beta(\mathbf{a})<+\infty$.
\begin{mytheorem}\label{thm:J_beta}
    Assume $J_\beta(\mathbf{a})<+\infty$. Then, we have the duality
    \begin{align}\label{J_beta_dual}
        J_\beta(\mathbf{a})=\sup_{\boldsymbol{\theta}\in \mathbb{R}^d}\boldsymbol{\theta}^T\mathbf{a}-K_\beta(\boldsymbol{\theta}),
    \end{align}
    where $K_\beta(\boldsymbol{\theta})=\int_{\Omega}\beta^*(\omega,\boldsymbol{\theta})\mathrm{d}\mathbb{P}(\omega)$ and $\beta^*(\omega,\boldsymbol{\theta})$ is the convex conjugate of $\beta(\omega,\mathbf{y})$. Moreover, if $J_\beta(\mathbf{a})>-\infty$, then there exists $\boldsymbol{\theta}^*\in \mathbb{R}^d$ such that the supremum in \eqref{J_beta_dual} is attained.
\end{mytheorem}
\begin{proof}
    Applying Theorem \ref{thm:Luenberger} to \eqref{J_beta1} yields that
    \begin{align*}
        J_\beta(\mathbf{a})=\sup_{\boldsymbol{\theta}\in \mathbb{R}^d}\boldsymbol{\theta}^T\mathbf{a}+\inf_{\mathbf{g}\in \mathcal{L}^1}\int_\Omega\beta(\omega,\mathbf{g}(\omega))-\boldsymbol{\theta}^T\mathbf{g}(\omega)\mathrm{d}\mathbb{P}(\omega).
    \end{align*}
We have,
\begin{align*}
    &\inf_{\mathbf{g}\in \mathcal{L}^1}\int_{\Omega}\beta(\omega,\mathbf{g}(\omega))-\boldsymbol{\theta}^T\mathbf{g}(\omega) \mathrm{d}\mathbb{P}(\omega)\\
    &\stackrel{(*)}{=}\int_{\Omega}\inf_{\mathbf{y}\in \mathbb{R}^d}(\beta(\omega,\mathbf{y})-\boldsymbol{\theta}^T\mathbf{y})\mathrm{d}\mathbb{P}(\omega)\\
    &\stackrel{(**)}{=}\int_{\Omega}\inf_{\mathbf{y}\geq \mathbf{0}}(\beta(\omega,\mathbf{y})-\boldsymbol{\theta}^T\mathbf{y})\mathrm{d}\mathbb{P}(\omega)\\
    &=-\int_{\Omega}\beta^*(\omega,\boldsymbol{\theta})\mathrm{d}\mathbb{P}(\omega),
\end{align*}
where the interchange of infimum and integral in $(*)$ follows from [\citet{Rockafellar&Wets}, Theorem 14.60], for which we may apply the theorem since  $J_\beta(\mathbf{a})<\infty$ and $\mathcal{L}^1$ is decomposable (see \citet{Rockafellar&Wets}, Definition 14.59). In $(**)$ we used that $\beta(\omega,\mathbf{y})=+\infty$ if $\mathbf{y}$ has negative component. Finally, the existence of $\boldsymbol{\theta}^*$ is guaranteed by Theorem \ref{thm:Luenberger}.
\end{proof}

Hence, the \textit{primal problem} $J_{\beta}(\mathbf{a})$ can be reformulated into a finite-dimensional \textit{dual problem} as in \eqref{J_beta_dual}. We call the solutions of the corresponding problems the \textit{primal} and \textit{dual solution} respectively. The following Theorem establishes the condition in which a convex function of the form $\boldsymbol{\theta}\mapsto  K_\beta(\boldsymbol{\theta})=\int_\Omega \beta^*(\omega,\boldsymbol{\theta})\mathrm{d}\mathbb{P}(\omega)$ is differentiable.
\begin{mytheorem}\label{thm:diff_kbeta}
    Let $\boldsymbol{\theta}\in \mathbb{R}^d$. If $\mathrm{int}(\mathrm{dom}(K_\beta))\neq \emptyset$ and $\beta^*(\omega,\cdot)$ is differentiable on its interior effective domain for $\mathbb{P}$-a.e. $\omega\in \Omega$. Then, $K_\beta$ is differentiable on $\mathrm{int}(\mathrm{dom}(K_\beta))$ and we have
    \begin{align*}
        \nabla K_\beta(\boldsymbol{\theta})=\int_\Omega \nabla \beta^*(\omega,\boldsymbol{\theta})\mathrm{d}\mathbb{P}(\omega).
    \end{align*}
\end{mytheorem}
\begin{proof} Since $K_\beta$ is a convex function, we can show its differentiability by showing its directional derivative function is linear (\citet{Rockafellar}, Theoreom 25.2). Let $\boldsymbol{\theta}_0\in \mathrm{int}(\mathrm{dom}(K_\beta))$. For any other $\boldsymbol{\theta}\in \mathrm{int}(\mathrm{dom}(K_\beta))$, we examine its directional derivative
\begin{align*}
    \lim_{t\downarrow 0}\frac{1}{t}\left(K_\beta(\boldsymbol{\theta}_0+t(\boldsymbol{\theta}-\boldsymbol{\theta}_0))-K_\beta(\boldsymbol{\theta}_0)\right)= \lim_{t\downarrow 0}\int_\Omega\frac{1}{t}\left(\beta^*(\omega,\boldsymbol{\theta}_0+t(\boldsymbol{\theta}-\boldsymbol{\theta}_0))-\beta^*(\omega,\boldsymbol{\theta}_0)\right)\mathrm{d}\mathbb{P}(\omega).
\end{align*}
Since, $\boldsymbol{\theta}_0\in \mathrm{int}(\mathrm{dom}(K_\beta))$, we have that $\boldsymbol{\theta}_0\in \mathrm{int}(\mathrm{dom}(\beta^*(\omega,\cdot)))$ for $\mathbb{P}$-a.e. $\omega\in \Omega$. For these $\omega$, the difference quotient 
\begin{align*}
    q_\omega(t):=\frac{1}{t}(\beta^*(\omega,\boldsymbol{\theta}_0+t(\boldsymbol{\theta}-\boldsymbol{\theta}_0))-\beta^*(\omega,\boldsymbol{\theta}_0))
\end{align*}
is non-increasing as $t\downarrow 0$ (Theorem 23.1, \citet{Rockafellar}). Therefore, the quotient $q_\omega(t)-q_\omega(\frac{1}{2})$ is negative and non-increasing for $t\leq 1/2$. Hence, applying the monotone convergence theorem and removing the term $q_\omega(\frac{1}{2})$ afterward (note that $q_\omega(\frac{1}{2})$ is a finite number since $\boldsymbol{\theta},\boldsymbol{\theta}_0\in \mathrm{dom}(\beta^*(\omega,\cdot))$, for $\mathbb{P}$-a.e. $\omega\in \Omega$) shows that $K_\beta$ is differentiable at $\boldsymbol{\theta}_0$. Indeed,
\begin{align*}
    \lim_{t\downarrow 0}\frac{1}{t}\left(K_\beta(\boldsymbol{\theta}_0+t(\boldsymbol{\theta}-\boldsymbol{\theta}_0))-K_\beta(\boldsymbol{\theta}_0)\right)&=\int_\Omega \lim_{t\downarrow 0}\frac{1}{t}(\beta^*(\omega,\boldsymbol{\theta}_0+t(\boldsymbol{\theta}-\boldsymbol{\theta}_0))-\beta^*(\omega,\boldsymbol{\theta}_0))\mathrm{d}\mathbb{P}(\omega)\\
    &\stackrel{(*)}{=}\int_\Omega \langle \boldsymbol{\theta}-\boldsymbol{\theta}_0 ,\nabla \beta^*(\omega,\boldsymbol{\theta}_0)\rangle \mathrm{d}\mathbb{P}(\omega)\\
    &\stackrel{(**)}{=}\left\langle \boldsymbol{\theta}-\boldsymbol{\theta}_0 ,\int_\Omega\nabla \beta^*(\omega,\boldsymbol{\theta}_0)\mathrm{d}\mathbb{P}(\omega)\right\rangle,
\end{align*}
where in $(*)$ we used the differentiability of $\beta^*(\omega,\cdot)$ and in $(**)$ we used the continuity of inner product.
\end{proof}
\begin{mytheorem}\label{thm:wc_density}
Assume $J_\beta(\mathbf{a})<+\infty$. Let $\boldsymbol{\theta}^*$ be such that $K_\beta(\boldsymbol{\theta}^*)<\infty$.  If the partial derivatives $\mathbf{g}_{\boldsymbol{\theta}^*}(\omega)=\nabla \beta^*(\omega,\boldsymbol{\theta}^*)$ exist for $\mathbb{P}$-almost every $\omega$ and satisfy:
\begin{align}\label{optimal_condition}
    \int_\Omega \mathbf{g}_{\boldsymbol{\theta}^*}(\omega)\mathrm{d}\mathbb{P}(\omega)=\mathbf{a}.
\end{align}
Then, $\boldsymbol{\theta}^*$ is a dual solution and $\mathbf{g}_{\boldsymbol{\theta}^*}$ is a primal solution.
On the other hand, if a primal solution $\mathbf{g}^*$ exists, then $\mathbf{g}^*=\mathbf{g}_{\boldsymbol{\theta}^*}$ for a particular dual solution $\boldsymbol{\theta}^*$, if the partial derivatives $\nabla \beta^*(\omega,\boldsymbol{\theta}^*)$ at this particular $\boldsymbol{\theta}^*$ exists, for $\mathbb{P}$-almost every $\omega\in \Omega$.
\end{mytheorem}
\begin{proof}
By the definition of the conjugate, we have for any $\boldsymbol{\theta}\in \mathbb{R}^d$ and functions $\mathbf{g}$, the inequality
\begin{align}\label{conj_ineq}
    \beta(\omega,\mathbf{g}(\omega))+\beta^*(\omega,\boldsymbol{\theta})\geq \boldsymbol{\theta}^T\mathbf{g}(\omega).
\end{align}
Inequality \eqref{conj_ineq} is tight if and only if the gradient of $\beta^*(\omega,\cdot)$ at $\boldsymbol{\theta}$ exists and is equal to $\mathbf{g}(\omega)$  (\citep{Rockafellar}, Theorem 23.5). 
Integrating both sides of \eqref{conj_ineq} with respect to $\omega$, for any function $\mathbf{g}$ such that $\int_\Omega \mathbf{g}(\omega)\mathrm{d}\mathbb{P}(\omega)=\mathbf{a}$ yields
\begin{align}\label{conj_ineq_int}
\begin{split}
\int_{\Omega}\beta(\omega,\mathbf{g}(\omega))\mathrm{d}\mathbb{P}(\omega)+\int_{\Omega}\beta^*(\omega,\boldsymbol{\theta})\mathrm{d}\mathbb{P}(\omega)&\geq \boldsymbol{\theta}^T\mathbf{a}.
\end{split}
\end{align}
In particular, inequality \eqref{conj_ineq} is an equality for $\mathbf{g}_{\boldsymbol{\theta}^*}(\omega)=\nabla \beta^*(\omega,\boldsymbol{\theta}^*)$ with $\boldsymbol{\theta}^*$ that satisfies condition \eqref{optimal_condition}. Hence, we have that \eqref{conj_ineq_int} is also an equality and thus
\begin{align*}
\int_{\Omega}\beta(\omega,\mathbf{g}_{\boldsymbol{\theta}^*}(\omega))\mathrm{d}\mathbb{P}(\omega)= \mathbf{a}^T\boldsymbol{\theta}^*-\int_{\Omega}\beta^*(\omega,\boldsymbol{\theta}^*)\mathrm{d}\mathbb{P}(\omega).
\end{align*}
Therefore, by Theorem \ref{thm:J_beta} we conclude that $\boldsymbol{\theta}^*$ is a dual solution and $\mathbf{g}_{\boldsymbol{\theta}^*}$ is a primal solution.

Conversely, let $\mathbf{g}^*$ be a primal solution. Then, $J_{\beta}(\mathbf{a})$ is finite, and thus a dual solution $\boldsymbol{\theta}^*$ exists and we must have $K_\beta(\boldsymbol{\theta}^*)<\infty$. Furthermore, \eqref{conj_ineq_int} holds for $\mathbf{g}^*$ and $\boldsymbol{\theta}^*$ and is an equality. Therefore, we have that the following integral with positive integrand is zero:
\begin{align*}
    \int_\Omega \beta(\omega,\mathbf{g}^*(\omega))+\beta^*(\omega,\boldsymbol{\theta}^*)-\mathbf{a}^T\boldsymbol{\theta}^*~\mathrm{d}\mathbb{P}(\omega)=0.
\end{align*}
Hence, \eqref{conj_ineq} is an equality for $\mathbf{g}^*$ and $\boldsymbol{\theta}^*$ for $\mathbb{P}$-a.e. $\omega\in \Omega$. Therefore, $\mathbf{g}^*=\mathbf{g}_{\boldsymbol{\theta}^*}$, by \citet{Rockafellar}, Theorem 23.5.
\end{proof}

\section{Characterization of Risk Aversion of Robust OCE and its Preservation of Convex Order}
The formulation of the robust OCE and robust shortfall risk measures in terms of $(\phi_1^*,\phi_2^*)$ also allows us to characterize risk aversion by simple convexity conditions imposed on $\phi_1^*,\phi_2^*$. We say that a robust risk measure $\rho\in \{\rho^l_{\mathrm{oce},\mathbb{P}}, \rho^l_{\mathrm{sf},\mathbb{P}}\}$ is risk-averse, if and only if
\begin{equation}\label{Risk_Aversion}
    \rho(X)\geq \mathbb{E}_{\mathbb{P}}[-X], \forall X.
\end{equation}
We show that for the robust OCE and shortfall risk measures, the risk aversion property can be easily characterized by some mild conditions on $\phi_1^*,\phi_2^*$.
\begin{proposition}\label{prop: Risk_Aversion}
For any non-decreasing convex function $\phi_1^*,\phi_2^*$ with $\phi_1^*(0)=\phi_2^*(0)=0$, we have that $\rho^l_{\mathrm{oce},\mathbb{P}}$ satisfies \eqref{Risk_Aversion}, if and only if 
    \begin{equation*}
        \phi_2^*(x)\geq x, \phi_1^*(x)\geq x, \forall x\in \mathbb{R}.
    \end{equation*}
    If we further assume that $\phi_1^*(x)>0, \forall x>0$, then we have that $\rho^l_{\mathrm{sf},\mathbb{P}}$ satisfies \eqref{Risk_Aversion}, if and only if
    \begin{equation*}
        \phi_2^*(x)\geq x, \forall x\in \mathbb{R}.
    \end{equation*}
\end{proposition}
\begin{proof}[\textbf{Proof of Proposition \ref{prop: Risk_Aversion}}]
For $\rho^l_{\mathrm{oce},\mathbb{P}}$, we have that \eqref{Risk_Aversion} implies that for all $x\in \mathbb{R}$, that 
\begin{align*}
    \inf_{\theta_1,\theta_2\in \mathbb{R}}\left\{-\theta_1-\theta_2+\phi_2^*(\phi_1^*(\theta_2-x)+\theta_1)\right\}\geq -x,
\end{align*}
which implies in particular that for $\theta_2=x$:
\begin{align*}
    \phi_2^*(\theta_1)\geq \theta_1, \forall \theta_1\in \mathbb{R}.
\end{align*}
It also implies that for $\theta_1=-\phi_1^*(x),\theta_2=2x$, that
\begin{align*}
    \phi_1^*(x)\geq x, \forall x\in \mathbb{R}.
\end{align*}
Conversely, if $\phi_2^*(x)\geq x, \phi_1^*(x)\geq x, \forall x\in \mathbb{R}$, then by the monotonicity of $\phi_1^*,\phi_2^*$, we have that for any $\theta_1,\theta_2\in\mathbb{R}$ and any $X$:
\begin{align*}
    -\theta_1-\theta_2+\mathbb{E}_{\mathbb{P}}[\phi_2^*(\phi_1^*(\theta_2-X)+\theta_1)]\geq -\theta_1-\theta_2+\mathbb{E}_{\mathbb{P}}[\theta_2-X+\theta_1]\geq \mathbb{E}_{\mathbb{P}}[-X]. 
\end{align*}
For $\rho^l_{\mathrm{sf},\mathbb{P}}$, we have that \eqref{Risk_Aversion} implies that for all $x\in \mathbb{R}$,
\begin{align*}
    \inf_{\theta_1,\theta_2\in \mathbb{R}}\left\{\theta_2~\middle |~\phi_2^*(\phi_1^*(-\theta_2-x)+\theta_1)\leq \theta_1\right\}\geq -x
\end{align*}
This implies that for all $x,\theta_1\in \mathbb{R}$, 
\begin{align*}
    \phi_2^*(\phi_1^*(-(-x)-x)+\theta_1)=\phi_2^*(\theta_1)\geq \theta_1.
\end{align*}
On the other hand, if $\phi_2^*(x)\geq x$ for all $x\in \mathbb{R}$, then for any $\theta_1\in \mathbb{R}$ and any random variable $X$,
\begin{align*}
    \mathbb{E}_{\mathbb{P}}[\phi_2^*(\phi_1^*(-\theta_2-X)+\theta_1)]\leq \theta_1&\Rightarrow \mathbb{E}_{\mathbb{P}}[\phi_1^*(-\theta_2-X)]\leq 0\\
    &\Rightarrow \phi_1^*(-\theta_2+\mathbb{E}[-X])\leq 0
\end{align*}
By the non-decreasing property of $\phi_1^*$ and the assumption that $\phi_1^*(x)>0, \forall x>0$, we can then conclude that the above inequality implies that $\theta_2\geq \mathbb{E}[-X]$. Therefore, $\rho^l_{\mathrm{sf},\mathbb{P}}(X)\geq \mathbb{E}_{\mathbb{P}}[-X]$.
\end{proof}

Furthermore, it is also easy to show that the robust OCE and shortfall risk measure preserves convex order. We say that $X$ is convex less order than $Y$, i.e. $X\leq_{\mathrm{cv}} Y$  if for all convex function $f$, we have $\mathbb{E}_\mathbb{P}[f(X)]\leq \mathbb{E}_{\mathbb{P}}[f(Y)]$.
\begin{proposition}\label{prop: CV_Order}
    Let $\phi_1^*,\phi_2^*$ be a non-decreasing convex function. Then 
    \begin{align*}
        X\leq_{\mathrm{cv}} Y &\Rightarrow  \rho^l_{\mathrm{oce},\mathbb{P}}(X)\leq \rho^l_{\mathrm{oce},\mathbb{P}}(Y) \\
        &\qquad\mathrm{and}~ \rho^l_{\mathrm{sf},\mathbb{P}}(X)\leq \rho^l_{\mathrm{sf},\mathbb{P}}(Y).
    \end{align*}
\end{proposition}

\begin{proof}[\textbf{Proof of Proposition \ref{prop: CV_Order}}]
    Suppose $X\leq_{\mathrm{cv}} Y$, then, we have that for all $\theta_1,\theta_2\in \mathbb{R}$,
    \begin{align*}
        \mathbb{E}_{\mathbb{P}}[\phi_2^*(\phi_1^*(\theta_2-X)+\theta_1)]\geq \mathbb{E}_{\mathbb{P}}[\phi_2^*(\phi_1^*(\theta_2-Y+\theta_1)], 
    \end{align*}
    since the function $\phi_2^*(\phi_1^*(\theta_2-x)+\theta_1)$ is convex in $x$. Hence, $\rho^l_{\mathrm{oce},\mathbb{P}}(X)\leq \rho^l_{\mathrm{oce},\mathbb{P}}(Y)$.

    Similarly, the function $\phi_2^*(\phi_1^*(-\theta_2-x)+\theta_1)$ is also convex in $x$. Therefore, for any $\theta_1,\theta_2\in \mathbb{R}$, if $X\leq_{\mathrm{cv}} Y$, then
    \begin{align*}
        \mathbb{E}[\phi_2^*(\phi_1^*(-\theta_2-Y)+\theta_1)]\leq \theta_1 \Rightarrow \mathbb{E}[\phi_2^*(\phi_1^*(-\theta_2-X)+\theta_1)]\leq \theta_1.
    \end{align*}
    Hence, $\rho^l_{\mathrm{sf},\mathbb{P}}(X)\leq \rho^l_{\mathrm{sf},\mathbb{P}}(Y)$.
\end{proof}

\section{Additional Details on Tail Analysis }\label{app:Finiteness_expOCE}
\subsection{Details in Table \texorpdfstring{\ref{tab:exp_OCE_finite}}{TEXT}}
For exponential OCE with $\psi^*(s)=e^{s}-1$,
\begin{itemize}
    \item Weibull versus KL, consider $f(x)\sim \left(\frac{|x|}{\lambda}\right)^{k-1}e^{-\left(\frac{|x|}{\lambda}\right)^k}$ for $x\leq 0$ and $f(x)=0$ for $x>0$, with $k,\lambda>0$. Examine the integral
    \begin{align*}
        \int^0_{-\infty}\left(\frac{|x|}{\lambda}\right)^{k-1}e^{-\left(\frac{|x|}{\lambda}\right)^k}\left(e^{\theta_1+e^{\theta_2-x}-1}-1\right)dx
    \end{align*}
    Since for any $k,\lambda>0$, $\lim_{x\to -\infty}|x|^{k-1}e^{-\left(\frac{|x|}{\lambda}\right)^k+e^{\theta_2-x}}\to \infty$, the integral diverges.
    \item Weibull versus Polynomial $p>1$, The tail behavior of $\psi^*(\theta_1+e^{\theta_2-x}-1)$ is governed by $e^{-x\frac{p}{p-1}}$ for $x\to-\infty$. For $k<1$, we have that
    \begin{align*}
        \lim_{x\to -\infty}\left(\frac{|x|}{\lambda}\right)^{k-1}e^{-\left(\frac{|x|}{\lambda}\right)^k}e^{-x\frac{p}{p-1}}&=\lim_{x\to -\infty}e^{-\left(\frac{|x|}{\lambda}\right)^k+\frac{p}{p-1}|x|+(k-1)\log\left(\frac{|x|}{\lambda}\right)}\\
        &=\infty,
    \end{align*}
    since $|x|$ outgrows $|x|^k+\log(|x|)$ for $k<1$. Hence, the integral will also diverge.
    
    For $k=1$, we have $f(x)\sim e^{-\frac{|x|}{\lambda}}$. The tail behavior we have to examine is: 
    \begin{align*}
        e^{-\frac{|x|}{\lambda}}e^{|x|\frac{p}{p-1}}=e^{\left(\frac{p}{p-1}-\frac{1}{\lambda}\right)|x|},
    \end{align*}
    which is integrable if and only if $\frac{p}{p-1}<\frac{1}{\lambda}$.
    
    For $k>1$, we need to examine the tail behavior of 
    \begin{align*}
        e^{-\left(\frac{|x|}{\lambda}\right)^k+\frac{p}{p-1}|x|+(k-1)\log\left(\frac{|x|}{\lambda}\right)}=e^{-|x|^k\left(a+b|x|^{1-k}+c\frac{\log(|x|)}{|x|^k}\right)},
    \end{align*}
    for $a>0$ and $b,c$ constants. Since $k>1$, the term $\left(a+b|x|^{1-k}+c\frac{\log(|x|)}{|x|^k}\right)$ converges to $a$ and can thus be bounded by a constant. Therefore, the exponential term is integrable. 
    \item Log-normal versus KL, consider $f(x)\sim \frac{1}{|x|}e^{-\frac{(\log(|x|)-\mu)^2}{2\sigma^2}}$ for $x<0$ and $f(x)=0$ for $x\geq 0$. We have that:
    \begin{align*}
        \lim_{x\to -\infty}\frac{1}{|x|}e^{-\frac{(\log(|x|)-\mu)^2}{2\sigma^2}}e^{e^{\theta_2-x}}=\infty.
    \end{align*}
    Therefore, the resulting integral diverges.
    \item Log-normal versus Polynomial $p>1$. The tail behavior of $\psi^*(\theta_1+e^{\theta_2-x}-1)$ is governed by $e^{-x\frac{p}{p-1}}$. We have
    $\frac{1}{|x|}e^{-\frac{(\log(|x|)-\mu)^2}{2\sigma^2}}=e^{-\frac{(\log(|x|)-\mu)^2}{2\sigma^2}-\log(|x|)}$. Since $-x$ approaches faster to $\infty$ than any polynomial of $\log(|x|)$ as $x\to -\infty$, it diverges.
    \item Pareto versus KL, $f(x)\sim \frac{1}{|x|^{\alpha+1}}$, for $-x\leq -\beta$, $\alpha, \beta>0$ and zero otherwise. Clearly, $\lim_{x\to-\infty}\frac{1}{|x|^{\alpha+1}}e^{e^{\theta_2-x}}=\infty$.
    \item Pareto versus Polynomial $p>1$, as $x\to -\infty$, we have that $\psi^*(\theta_1+e^{\theta_2-x})\to e^{-x\cdot \frac{p}{p-1}}$ and $\lim_{x\to -\infty}\frac{1}{|x|^{\alpha+1}}e^{-x\cdot \frac{p}{p-1}}=\infty$. Hence, it diverges.
    \item Pareto, Log-normal, Weibull vs Polynomial $p<1$, $\psi^*$ becomes $+\infty$ for $x$ sufficiently negative, hence they all diverges.
\end{itemize}
\subsection{Details in Table \texorpdfstring{\ref{tab:CVaR_finite}}{TEXT}}\label{app:Finiteness_CVAR}
\begin{itemize}
    \item Pareto versus KL-divergence. For x sufficiently negative, we study the tail behavior of $\frac{1}{|x|^{\alpha_0+1}}e^{\frac{1}{1-\alpha}|x|}$, which diverges as $x\to -\infty$.
    \item Log-normal vs KL-divergence, we study the tail behavior of:
    \begin{align*}
    \frac{1}{|x|}e^{-\frac{(\log(|x|)-\mu)^2}{2\sigma^2}}e^{\frac{1}{1-\alpha}|x|},    
    \end{align*}
    which also diverges for $x\to -\infty$.
    \item Weibull vs KL-divergence, we study the tail behavior of
    \begin{align*}
        e^{-|x|^k+\frac{1}{1-\alpha}|x|+(k-1)\log|x|}
    \end{align*}
    For $k>1$, this is integrable, for $0<k<1$, this diverges. For $k=1$, it depends on the remaining constants which determines whether there is $e^{|x|}$ or $e^{-|x|}$ term.
    \item Polynomial $p>1$ versus Gaussian, we study the tail behavior of
    \begin{align*}
        e^{-\frac{(x-\mu)^2}{2\sigma^2}}|x|^{\frac{p}{p-1}},
    \end{align*}
    which is clearly integrable.
    \item Polynomial $p>1$ versus Log-normal, we study the tail behavior of:
    \begin{align*}
        |x|^{\frac{p}{p-1}}\frac{1}{|x|}e^{-\frac{(\log |x|-\mu)^2}{2\sigma^2}}=e^{(p-1)\log|x|-\frac{(\log |x|-\mu)^2}{2\sigma^2}},
    \end{align*}
    which is also integrable.
    \item Polynomial $p>1$ vs Weibull, the tail behavior is
    \begin{align*}
        |x|^{k-1}e^{-|x|^k}|x|^{\frac{p}{p-1}}=e^{-|x|^k+(k-1+\frac{p}{p-1})\log|x|},
    \end{align*}
    which is also integrable for all $k>0$.
    \item Polynomial $p>1$ vs Pareto, the tail behavior is
    \begin{align*}
        |x|^{-\alpha_0-1}|x|^{\frac{p}{p-1}}=|x|^{\frac{p}{p-1}-(\alpha_0+1)}
    \end{align*}
    The integrability depends clearly on whether $\frac{p}{p-1}-(\alpha_0+1)<-1$ or not.
     \item Pareto, Log-normal, Weibull vs Polynomial $p<1$, $\psi^*$ becomes $+\infty$ for $x$ sufficiently negative, hence they all diverges.
\end{itemize}
\end{document}